\documentclass[11pt]{amsart}

\usepackage{amsmath}
\usepackage{amsthm}
\usepackage{amssymb}
\usepackage{enumerate}
\usepackage{fullpage}
\usepackage[all,cmtip]{xy}
\usepackage[hidelinks]{hyperref}
\usepackage{url}
\usepackage{color}
\definecolor{purple}{rgb}{1,0,1}

\usepackage{tikz}
\usetikzlibrary{cd,angles,quotes,decorations.pathreplacing,math}

\usepackage{enumitem}
\usepackage{mathtools} 

\usepackage{bm}
\usepackage{stackengine}

\usepackage{graphicx}
\usepackage{rotating}
\usepackage{pdflscape}

\usepackage{float}
\floatstyle{boxed}
\restylefloat{figure}

\theoremstyle{plain}
\newtheorem{thm}{Theorem}[section]

\newtheorem*{thm*}{Theorem}
\newtheorem*{guess*}{Educated guess}
\newtheorem{lemma}[thm]{Lemma}
\newtheorem*{lemma*}{Lemma}
\newtheorem{cor}[thm]{Corollary}
\newtheorem*{cor*}{Corollary}
\newtheorem{prop}[thm]{Proposition}
\newtheorem*{prop*}{Proposition}

\newtheorem{fact}[thm]{Fact}

\newtheorem{dimensionTheorem}[thm]{Dimension Theorem}
\newtheorem{smallSetsTheorem}[thm]{Small Sets Theorem}
\newtheorem{kindaSmallSetsTheorem}[thm]{Kinda Small Sets Theorem}
\newtheorem{verySmallSetsTheorem}[thm]{Very Small Sets Theorem}

\theoremstyle{definition}
\newtheorem{definition}[thm]{Definition}
\newtheorem*{def*}{Definition}
\newtheorem{remark}[thm]{Remark}
\newtheorem*{remark*}{Remark}

\newtheorem{example}[thm]{Example}
\newtheorem*{example*}{Example}

\newtheorem*{question*}{Question}
\newtheorem*{openQuestion*}{Open Question}

\numberwithin{equation}{section}

\newcommand{\fr}{\mathfrak}

\newcommand\N{\mathbb{N}}
\newcommand\Q{\mathbb{Q}}
\newcommand\R{\mathbb{R}}
\newcommand\Z{\mathbb{Z}}

\newcommand\calL{\mathcal{L}}

\newcommand\calS{\mathcal{S}}

\DeclareMathOperator\lspan{span}
\DeclareMathOperator\Qspan{\lspan_{\Q}}

\newcommand\0{\varnothing}

\DeclareMathOperator\image{image}

\DeclareMathOperator\acl{acl}

\DeclareMathOperator\cl{cl}
\DeclareMathOperator\rk{rk}

\renewcommand\L{\mathcal{L}}

\DeclareMathOperator\Th{Th}

\DeclareMathOperator\dpr{dp}

\renewcommand\d{\operatorname{d}}

\renewcommand\k{\boldsymbol{k}}

\DeclareMathOperator\res{res}

\DeclareFontFamily{U}{fsy}{}
\DeclareFontShape{U}{fsy}{m}{n}{<->s*[.9]psyr}{}
\DeclareSymbolFont{der@m}{U}{fsy}{m}{n}
\DeclareMathSymbol{\der}{\mathord}{der@m}{182}

\DeclareSymbolFont{der@m}{U}{fsy}{m}{n}
\DeclareMathSymbol{\derdelta}{\mathord}{der@m}{100}

\DeclareFontFamily{OMS}{smallo}{}
\DeclareFontShape{OMS}{smallo}{m}{n}{<->s*[.65]cmsy10}{}
\DeclareSymbolFont{smallo@m}{OMS}{smallo}{m}{n}
\DeclareMathSymbol{\smallo}{\mathord}{smallo@m}{79}

\DeclareSymbolFont{imag@m}{OT1}{cmr}{m}{ui}
\DeclareMathSymbol{\imag}{\mathord}{imag@m}{105}

\newcommand{\overbar}[1]{\mkern 2mu\overline{\mkern-2mu#1\mkern-2mu}\mkern 2mu}
\DeclareFontFamily{U}{mathb}{\hyphenchar\font45}
\DeclareFontShape{U}{mathb}{m}{n}{
      <5> <6> <7> <8> <9> <10> gen * mathb
      <10.95> matha10 <12> <14.4> <17.28> <20.74> <24.88> matha12
      }{}
\DeclareSymbolFont{mathb}{U}{mathb}{m}{n}
\DeclareMathSymbol{\monus}{2}{mathb}{"01}

\def\Ind#1#2{#1\setbox0=\hbox{$#1x$}\kern\wd0\hbox to 0pt{\hss$#1\mid$\hss}\lower.9\ht0\hbox to 0pt{\hss$#1\smile$\hss}\kern\wd0}
\def\Notind#1#2{#1\setbox0=\hbox{$#1x$}\kern\wd0\hbox to 0pt{\mathchardef
	\nn=12854\hss$#1\nn$\kern1.4\wd0\hss}\hbox to
	0pt{\hss$#1\mid$\hss}\lower.9\ht0 \hbox to
	0pt{\hss$#1\smile$\hss}\kern\wd0}

\newbox\gnBoxA
\newdimen\gnCornerHgt
\setbox\gnBoxA=\hbox{$\ulcorner$}
\global\gnCornerHgt=\ht\gnBoxA
\newdimen\gnArgHgt
\def\Code #1{%
\setbox\gnBoxA=\hbox{$#1$}%
\gnArgHgt=\ht\gnBoxA%
\ifnum     \gnArgHgt<\gnCornerHgt \gnArgHgt=0pt%
\else \advance \gnArgHgt by -\gnCornerHgt%
\fi \raise\gnArgHgt\hbox{$\ulcorner$} \box\gnBoxA %
\raise\gnArgHgt\hbox{$\urcorner$}}

\author{Allen Gehret}
\address{Czech Technical University in Prague, Artificial Intelligence Center, Charles Square 13, Prague 2, Czech Republic}
\address{Universit\"{a}t Wien, Institut f\"{u}r Mathematik, Kurt G\"{o}del Research Center, Kolingasse 14--16, 1090 Wien, Austria}
\email{gehreall@fel.cvut.cz}
\email{allen.gehret@univie.ac.at}

\author{Elliot Kaplan}
\address{Max-Planck-Institut f\"{u}r Mathematik, Vivatsgasse 7, 53111 Bonn, Germany}
\address{Universit\'{e} de Mons, D\'{e}partement de Math\'{e}matique, Avenue Maistriau 15, 7000 Mons, Belgium}
\email{ekaplan@mpim-bonn.mpg.de}
\email{elliot.kaplan@umons.ac.be}

\author{Nigel Pynn-Coates}
\address{Universit\"{a}t Wien, Institut f\"{u}r Mathematik, Kurt G\"{o}del Research Center, Kolingasse 14--16, 1090 Wien, Austria}
\email{nigel.pynn-coates@univie.ac.at}

\title{Dimension theory for the asymptotic couple of the field of logarithmic transseries}
\begin{document}

\begin{abstract}
In this paper we completely characterize all dimension functions on all models of the theory $T_{\log}$ of the asymptotic couple of the field of logarithmic transseries (Dimension Theorem). This is done by characterizing the ``small'' $1$-variable definable sets (Small Sets Theorem). As a byproduct, we show that $T_{\log}$ is d-minimal and does not eliminate imaginaries. Separately, we provide an abstract criterion for d-minimality, which we use to observe some new examples of d-minimal expansions of valued fields.
\end{abstract}
\maketitle

\section{Introduction}\label{section_introduction}

\noindent
The differential field $\mathbb{T}_{\log}$ of logarithmic transseries is conjectured to have nice model-theoretic properties~\cite{gehret2017towards}. As $\mathbb{T}_{\log}$ is a so-called \emph{$H$-field}~\cite{ADAMTT}, it is an expansion of a \emph{valued} differential field; as such, in the Ax--Kochen--Ershov (AKE) tradition we view $\mathbb{T}_{\log}$ in terms of the ``decomposition'': 
\[
\begin{tikzcd}
&\arrow[ld]\mathbb{T}_{\log} \arrow[rd]  &\\
\R&& (\Gamma_{\log},\psi)
\end{tikzcd}
\]
Here, $\R$ is simultaneously the residue field of the valuation and the constant field of the derivation (conjectured to have semialgebraic induced structure). The object $\Gamma_{\log}$ is the value group of the valuation, further equipped with a map $\psi\colon\Gamma_{\log}\to\Gamma_{\log}\cup\{\infty\}$ induced by the logarithmic derivative of $\mathbb{T}_{\log}$.
Collectively, the pair $(\Gamma_{\log},\psi)$ is called the \emph{asymptotic couple} of $\mathbb{T}_{\log}$.

\medskip\noindent
Incidentally, the object $(\Gamma_{\log},\psi)$ is also the asymptotic couple of the Hardy field $\R(\fr{L})$, i.e., the Hardy field over $\R$ which is generated by all power-products of the form $x^{r_0}(\log x)^{r_1}(\log^{\circ 2}x)^{r_2}\cdots(\log^{\circ n}x)^{r_n}$ for varying $n$, with each $r_i\in\R$.
As pointed out in~\cite{elliott2024analytic_book}: ``many functions in number theory are comparable, or are conjectured to be comparable'' to a function in $\R(\fr{L})$. We contend that the asymptotic couple $(\Gamma_{\log},\psi)$ is an appropriate universal domain for capturing the generic asymptotic behavior of the functions in $\R(\fr{L})$, as well as their interaction with the (logarithmic) derivative.

\medskip\noindent
This paper is the fourth in a series~\cite{GehretACTlog,GehretACTlogNIP,GehretKaplan-ACTlogdistal} about the theory $T_{\log}$ (defined in Section~\ref{section_overview_Tlog}) of the asymptotic couple $(\Gamma_{\log},\psi)$. To summarize, here are the most relevant things already known:
\begin{itemize}
\item $T_{\log}$ has quantifier elimination (QE) and a universal axiomatization (UA) in a natural language~\cite{GehretACTlog}; in particular, the quantifier-free definable sets enjoy a Tarski--Seidenberg-type theorem, and definable functions are given piecewise by terms in the language.
\item $T_{\log}$ has the non-independence property (NIP)~\cite{GehretACTlogNIP}, i.e., definable families have finite \emph{Vapnik-Chervonenkis (VC) dimension}; in particular, this implies definable hypothesis spaces are subject to the so-called \emph{Fundamental Theorem of Statistical Learning}, which tells us (modulo some measurability assumptions) they are always \emph{PAC learnable}~(see, e.g.,~\cite{chase2019model,krapp2024measurability} for more on this connection). 
Moreover, $T_{\log}$ is not strongly dependent and so it is not dp-minimal nor does it have finite dp-rank~\cite{GehretKaplan-ACTlogdistal}.
\item $T_{\log}$ is distal~\cite{GehretKaplan-ACTlogdistal}; in particular, definable relations satisfy strong combinatorial bounds including a definable version of the \emph{strong Erdős--Hajnal property}~\cite{chernikov2018regularity,chernikov2020cutting}. 
\item The model-theoretic algebraic closure $\operatorname{acl}$ is in general \emph{not} a pregeometry~\cite{GehretACTlogNIP}; in particular, there is no immediate off-the-shelf \emph{dimension theory} we can use similar to the likes of \emph{vector spaces},  \emph{algebraically closed fields}, and \emph{o-minimal structures}.
\end{itemize}

\noindent
In this paper, we further examine the nature of definable sets in models of $T_{\log}$ from topological, algebraic, and model-theoretic perspectives.
Specifically, we answer the following questions:

\begin{question*}
What are the dimension functions on models of $T_{\log}$?
\end{question*}
\begin{proof}[Answer]\renewcommand{\qedsymbol}{}
The Dimension Theorem~\ref{dimensionTheorem_label} completely characterizes all dimension functions on models of $T_{\log}$, where \emph{dimension function} is meant in the axiomatic sense of~\cite{vdD89}; see Definition~\ref{def_dimfn}.
Corollary~\ref{agreement_topological_dimension} shows that each such dimension coincides with an appropriate \emph{topological dimension}.
\end{proof}

\begin{question*}
What is a more precise description of the $1$-variable definable subsets of models of $T_{\log}$?
\end{question*}
\begin{proof}[Answer]\renewcommand{\qedsymbol}{}
The Small Sets Theorem~\ref{smallSetsTheorem_label} characterizes the ideal of ``small'' unary definable sets (i.e., sets of dimension $\leq 0$) for each dimension function, and forms the technical core of this paper.
\end{proof}

\begin{question*}
What is a more precise description of the $n$-variable definable subsets of models of $T_{\log}$? 
\end{question*}
\begin{proof}[Answer]\renewcommand{\qedsymbol}{}
This is the Kinda Small Sets Theorem~\ref{smallishsets} and the Very Small Sets Theorem~\ref{corverysmallsets}. When $n>1$, the conditions in the Small Sets Theorem~\ref{smallSetsTheorem_label} generalize in two ways: they either characterize sets of dimension $<n$ (``kinda small sets'') or sets of dimension $\leq 0$ (``very small sets'').
\end{proof}

\begin{question*}
What is a more precise description of definable functions in models of $T_{\log}$?
\end{question*}
\begin{proof}[Answer]\renewcommand{\qedsymbol}{}
Corollary~\ref{cor_locallin} asserts that every definable function $f\colon\Gamma^n\to\Gamma_{\infty}$ is either locally affine or locally constant outside of a ``kinda small set''.
Conversely, Proposition~\ref{characterization_of_functions_from_Psin} characterizes functions defined on the ``typical very small set'' $\Psi^n$. Understanding these two extremes is sufficient for many questions about definable functions due to the inductive nature of dimension.
\end{proof}

\begin{question*}
Is the theory $T_{\log}$ d-minimal?
\end{question*}
\begin{proof}[Answer]\renewcommand{\qedsymbol}{}
Yes, i.e., every $1$-variable definable set in any model either has nonempty interior or is a finite union of discrete sets (Definition~\ref{d_minimality_def_new}). First, this is apparent from \eqref{smallequiv2}$\Leftrightarrow$\eqref{smallequiv3} of the Small Sets Theorem~\ref{smallSetsTheorem_label}. A second proof follows readily from a general \emph{d-minimality criterion}, Proposition~\ref{d-minimality_criterion_ms}. Finally, Corollary~\ref{corfornasierodmin} shows that $T_{\log}$ is d-minimal in the stronger sense of~\cite[Definition 9.1]{fornasiero-dimmatroiddense}, which places additional topological requirements on the definable sets in $n$-variables; this relies on the Kinda Small and Very Small Sets Theorems~\ref{smallishsets} and~\ref{corverysmallsets}. 
\end{proof}

\begin{question*}
Does the $1$-sorted theory $T_{\log}$ have elimination of imaginaries (EI)?
\end{question*}
\begin{proof}[Answer]\renewcommand{\qedsymbol}{}
No. The analogue of the \emph{RV sort} cannot be eliminated (Corollary~\ref{RV_sort_not_EI}). This uses properties of dimension and a cardinality argument in what we call \emph{the standard model}.
\end{proof}

\subsection*{Overview and main ideas}
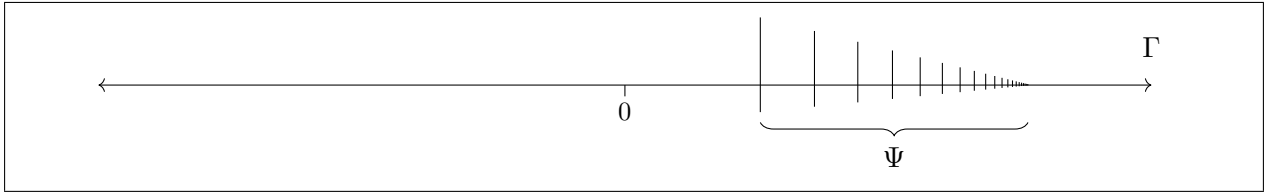
\begin{figure}[h!]
\caption{The standard model of $T_{\log}$. }
\label{standardmodel}
\begin{center}
\begin{tikzpicture}
\tikzmath{\x = .8; \y = .9; \z =2;
\w =1/(1-\y);}
\draw[<->] (-7,0)--(7,0); 
\tikzmath{\x = .8; \y = .9; \z =2;
\w =1/(1-\y);}
\draw (0,-.15)--(0,0);
\node at (0,-.35) {\small$0$}; 

\node at (7,.5) {$\Gamma$};

\draw [decorate,decoration={brace,amplitude=5pt,mirror,raise=3ex}]
  (3*\y*\z+\w*\y*\y*\z*2-2*\z*\w*\y*\y-2*\z*\y*\x^0,0) -- (3*\y*\z+\w*\y*\y*\z*2-2*\z*\w*\y*\y-2*\z*\y*\x^20,0) node[midway,yshift=-2.5em]{$\Psi$};

\foreach \a in {0,1,2,3,4,5,6,7,8,9,10,11,12,13,14,15,16,17,18,19,20} {
\begin{scope}
\draw (3*\y*\z+\w*\y*\y*\z*2-2*\z*\w*\y*\y-2*\z*\y*\x^\a, -.4*\y*\x^\a)--(3*\y*\z+\w*\y*\y*\z*2-2*\z*\w*\y*\y-2*\z*\y*\x^\a, 1/2*\z*\y*\x^\a);
\end{scope}}
\end{tikzpicture}
\end{center}
\end{figure}

Throughout, $m$ and $n$ range over $\N=\{0,1,2,\ldots\}$. Let $\bigoplus_n\R e_n$ be a vector space over $\R$ with basis $(e_n)$. Then $\bigoplus_n\R e_n$ can be made into an ordered group using the usual lexicographical order, i.e., by requiring for nonzero $\sum_ir_ie_i$, that
\[
\textstyle\sum r_ie_i > 0 \quad \iff \quad r_n > 0 \ \text{ for the least $n$ such that $r_n\neq 0$.}
\]
Let $\Gamma_{\log}$ be the above ordered abelian group $\bigoplus_n\R e_n$. It is convenient to think of an element $\sum r_ie_i$ as the vector $(r_0,r_1,r_2,\ldots)$. We follow Rosenlicht
~\cite{rosen-dvalgp2}
in taking the function:
\[
\psi \colon \Gamma_{\log}\setminus\{0\}\to\Gamma_{\log},\quad (\underbrace{0,\ldots,0}_{n},\underbrace{r_n}_{\neq 0},r_{n+1},\ldots)\ \mapsto \ (\underbrace{1,\ldots,1}_{n+1},0,0,\ldots)
\]
as a new primitive, calling the pair $(\Gamma_{\log},\psi)$ an \emph{asymptotic couple}. 
Throughout, we refer to this specific asymptotic couple, depicted in Figure~\ref{standardmodel}, as \emph{the standard model}. 

\medskip\noindent
The first key observation about the function $\psi$ is that it is a convex valuation on the ordered abelian group $\Gamma_{\log}$. Moreover, the \emph{value set} $\Psi$ of $\psi$ is a very important definable subset of $\Gamma_{\log}$:
\[
\Psi \ \coloneqq \ \psi(\Gamma_{\log}\setminus\{0\}) \ = \ \{(\underbrace{1,\ldots,1}_n,0,0,\ldots):n\geq 1\}.
\]
Note that this situation is a bit atypical in valuation theory, as usually a value set (or value group) of a valuation lives on its own more primitive sort, and is not a subset of the domain of the valuation. Furthermore, $\Psi$ introduces a discrete set into the otherwise ``continuous'' object $\Gamma_{\log}$, which itself is o-minimal as an ordered divisible abelian group. In general, these seem to be the main sources of complications when dealing with $(\Gamma_{\log},\psi)$.

\medskip\noindent
The story of dimension now begins with the following observation:

\begin{quote}
Although $\operatorname{acl}$ is not a pregeometry in general, the \emph{relativization} of $\operatorname{acl}$ to the definable set $\Psi$ is a pregeometry. Moreover, this relativization $X\mapsto \operatorname{acl}(X\cup\Psi)$ is essentially the same as the ``linear'' pregeometry $X\mapsto \Qspan(X\cup\Psi)$ coming from the underlying divisible abelian group structure of $\Gamma_{\log}$. 
\end{quote}

\noindent
This pregeometry gives rise to a \emph{dimension function}. However, the following issue still remains:
\begin{quote}
\emph{How many} dimension functions does a model of $T_{\log}$ have? Since there is no obvious definable field structure, uniqueness results such as~\cite[Theorem 3.48]{fornasiero-dimmatroiddense} do not apply.
\end{quote}

\noindent
The general role of \emph{coarsening} in the analysis of $H$-fields~\cite{ADAMTT} provides a natural guess at what the other dimension functions might be:
they are parametrized by a certain ``scale'' that is uniformly indexed by the $\Psi$-set of the asymptotic couple, which we now explain.

\begin{figure}[h!]\label{figquotient}
\caption{Quotienting by $\Delta_\phi$. }
\label{Quotient}
\begin{center}
\begin{tikzpicture}
\draw[<->] (-7,0)--(7,0); 
\tikzmath{\x = .8; \y = .9; \z =2;
\w =1/(1-\y); \u = -3;}
\draw (0,-.15)--(0,0);
\node at (3*\y*\z+\w*\y*\y*\z*2-2*\z*\w*\y*\y-2*\z*\y*\x^2, .9) {\small$\phi$}; 

\node at (7,.5) {$\Gamma$};

\node at (-1.3,0) {(};
\node at (1.3,0) {)};
\draw [decorate,decoration={brace,amplitude=5pt,mirror,raise=2ex}]
  (-1,0) -- (1,0) node[midway,yshift=-2em]{$\Delta_\phi$};

\foreach \a in {0,1,2,3,4,5,6,7,8,9,10,11,12,13,14,15,16,17,18,19,20} {
\begin{scope}
\draw (3*\y*\z+\w*\y*\y*\z*2-2*\z*\w*\y*\y-2*\z*\y*\x^\a, -.4*\y*\x^\a)--(3*\y*\z+\w*\y*\y*\z*2-2*\z*\w*\y*\y-2*\z*\y*\x^\a, 1/2*\z*\y*\x^\a);
\end{scope}}

\draw[<->] (-7,\u)--(7,\u); 

\node at (3*\y*\z+\w*\y*\y*\z*2-2*\z*\w*\y*\y-2*\z*\y*\x^2+.3, .9+\u) {\small$\phi+\Delta_\phi$}; 

\draw (0,-.15+\u)--(0,\u);
\foreach \a in {0,1,2} {
\begin{scope}
\draw (3*\y*\z+\w*\y*\y*\z*2-2*\z*\w*\y*\y-2*\z*\y*\x^\a, -.4*\y*\x^\a+\u)--(3*\y*\z+\w*\y*\y*\z*2-2*\z*\w*\y*\y-2*\z*\y*\x^\a, 1/2*\z*\y*\x^\a+\u);
\end{scope}}

\node at (7,.5+\u) {$\Gamma/\Delta_\phi$};

\node at (0,-1.8) {\rotatebox[origin=c]{90}{$\Longleftarrow$}};
\end{tikzpicture}
\end{center}
\end{figure}
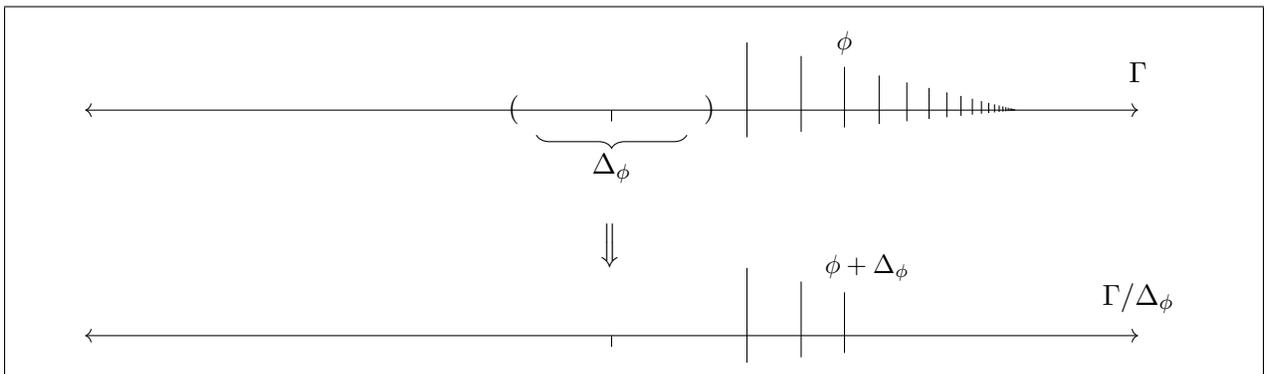

\medskip\noindent
For an arbitrary model $(\Gamma,\psi)$ of $T_{\log}$ we set $\Psi\coloneqq \psi(\Gamma^{\neq})$, where $\Gamma^{\neq}\coloneqq \Gamma\setminus\{0\}$.
Then for $\phi\in\Psi$ we define a proper convex subgroup of $\Gamma$:
\[
\Delta_{\phi} \ \coloneqq \ \{x\in\Gamma^{\neq}:\psi(x)>\phi\}\cup\{0\}.
\]
We extend this to $\phi\in\Psi\cup\{\infty\}$ by setting $\Delta_{\infty}\coloneqq\{0\}$. We may then proceed to study a definable set $X\subseteq\Gamma$ in terms of its image in the quotient $\Gamma/\Delta_{\phi}$, which we visualize in Figure~\ref{Quotient}.

\medskip\noindent
To do this, we consider the further relativization $X\mapsto \operatorname{acl}(X\cup\Psi\cup\Delta_{\phi})$, which is also a pregeometry and yields a notion of \emph{dimension} $\dim_{\phi}$ on definable sets. 
However, it is not obvious that $\dim_{\phi}$ is a \emph{dimension function} in the axiomatic sense of~\cite{vdD89} (indeed, not every pregeometric dimension is: the theory of $(\N;<)$ is pregeometric, but the corresponding dimension is not definable).
Our main theorem establishes that the $\dim_{\phi}$ are dimension functions, and they are the only ones:

\begin{dimensionTheorem}\label{dimensionTheorem_label}
Suppose $(\Gamma,\psi)$ is a model of $T_{\log}$. Then:
\begin{enumerate}
\item for each $\phi\in\Psi\cup\{\infty\}$, there exists a unique dimension function $\dim_{\phi}$ on $(\Gamma,\psi)$ such that $\dim_{\phi}\Delta_{\phi}=0$ and $\dim_{\phi}\Delta_{\xi}=1$ for all $\xi\in\Psi^{<\phi}$, and
\item if $\operatorname{d}$ is an arbitrary dimension function on $(\Gamma,\psi)$, then $\operatorname{d}=\dim_{\phi}$ for some $\phi\in\Psi \cup\{\infty\}$.
\end{enumerate}
\end{dimensionTheorem}

\noindent
This theorem follows from the Small Sets Theorem~\ref{smallSetsTheorem_label}, which characterizes the definable sets $X\subseteq\Gamma$ such that $\dim_{\phi}X\leq 0$. Part of this characterization includes the fact that $T_{\log}$ is d-minimal.

\medskip\noindent
At an earlier stage of this work, we also produced a direct proof that $T_{\log}$ is d-minimal by extracting and applying a d-minimality criterion, Proposition~\ref{d-minimality_criterion_ms}, from other proofs of d-minimality. This proof was later superseded by the finer analysis needed for the Small Sets Theorem~\ref{smallSetsTheorem_label}. Nevertheless, we have included this criterion because it may be useful for other topological theories. For example, we use it to show in Subsection~\ref{subsec_valueddmin} that henselian valued fields of equicharacteristic zero, equipped with a section of the valuation and a lift of the residue field, are d-minimal. See Section~\ref{sec_dminimality_criterion} for a self-contained treatment of this criterion and other applications.

\subsection{Outline of paper}
In Section~\ref{section_overview_Tlog} we review the theory $T_{\log}$ and establish some basic facts we need.
Section~\ref{sec_dimensions_SST} contains the main results of the paper; see the introduction of that section for an overview.
Section~\ref{sec_sstproof} is dedicated to the proof of the Small Sets Theorem;
in particular, Subsection~\ref{subsec_roadmap} provides a roadmap of the proof and the rest of the section consists of independent subsections contributing some part of the proof.
Section~\ref{sec_dminimality_criterion} provides an abstract criterion for d-minimality for topological theories expanded by unary functions, along with several applications.
In Section~\ref{section_final_comments} we make some final comments and observations.
Finally, in Appendix~\ref{sec_appendix} we collect a few basic topological facts and definitions that we use and prove some facts needed for Section~\ref{sec_dminimality_criterion}.

\subsection{Conventions}

\subsubsection*{Set theory conventions}
Given a boolean algebra $\mathcal{C}$ of subsets of a set $X$, recall that an \textbf{ideal} of $\mathcal{C}$ is a collection $\mathcal{I}\subseteq\mathcal{C}$ such that $\varnothing\in\mathcal{I}$, if $A,B\in\mathcal{I}$, then $A\cup B\in\mathcal{I}$, and if $A\subseteq B$ and $B\in\mathcal{I}$, then $A\in\mathcal{I}$.
Given a set $X\subseteq A\times B$ and $a\in A$, $X_a\coloneqq\{b\in B:(a,b)\in X\}$ denotes the \emph{fiber} of $X$ over $a$. We may use $\sqcup$ and $\bigsqcup$ instead of $\cup$ and $\bigcup$ to emphasize that a given union is a disjoint union.

\subsubsection*{Ordered set conventions} By ``ordered set'' we mean ``totally ordered set''. 
Let $S$ be an ordered set. Below, the ordering on $S$ will be denoted by $\leq$, and a subset of $S$ is viewed as ordered by the induced ordering. 
Suppose that $B$ is a subset of $S$. We put $S^{>B}\coloneqq\{s\in S:\text{$s>b$ for every $b\in B$}\}$ and $S^{>a} \coloneqq S^{>\{a\}}$; similarly for $\geq$, $<$, and $\leq$ instead of $>$.
For $a,b\in S$ we put
\[
[a,b]_B \ \coloneqq \ \{x\in B:a\leq x\leq b\}.
\]
If $B=S$, then we usually write $[a,b]$ instead of $[a,b]_S$. A subset $C$ of $B$ is said to be \textbf{convex} in $B$ if for all $a,b\in C$ we have $[a,b]_B\subseteq C$. For $a<b$ we also set
\[
(a,b) \ \coloneqq \ \{x\in X:a<x<b\}.
\]
The sets of the form $(a,b)$ are called \textbf{intervals} in $S$. The intervals form a basis for a hausdorff topology on $S$, the \textbf{order topology} on $S$.

\subsubsection*{Ordered abelian group conventions} Suppose that $G$ is an ordered abelian group. Then we set $G^{\neq}\coloneqq G\setminus\{0\}$. Also $G^{<}\coloneqq G^{<0}$; similarly for $\geq$, $\leq$, and $>$ instead of $<$. We define $|g|\coloneqq\max(g,-g)$ for $g\in G$. For $a\in G$, the \textbf{archimedean class} of $a$ is defined by
\[
[a] \ \coloneqq \ \{g\in G:\text{$|a|\leq n|g|$ and $|g|\leq n|a|$ for some $n\geq 1$}\}.
\]
Basic facts about ordered abelian groups and archimedean classes may be found in~\cite[2.4]{ADAMTT}.

\subsubsection*{Topology conventions} Given a subset $A$ of a topological space $X$, we let $A'$ denote the \textbf{derived set} of $A$ (in $X$), i.e., the set of limit points of $A$ in $X$.
For $n>0$, we set $A^{(n)} \coloneqq (A^{(n-1)})'$, where $A^{(0)} \coloneqq A$. We say that $A$ is \textbf{d-finite} (in $X$) if $A^{(n)}=\varnothing$ for some $n$. 
Given a function $f\colon X\to Y$ between topological spaces $X,Y$, we let $\mathsf{Discont}(f)$ denote the set of points of $X$ at which $f$ is discontinuous.
We let $\operatorname{br}A \coloneqq A \setminus \operatorname{int}A$ denote the \textbf{border} of $A$.  
See Appendix~\ref{sec_appendix} for more on d-finite sets, the border of $A$, and other topological facts we use. 

\subsubsection*{Model theory conventions} 
In general we adopt the conventions of~\cite[Appendix B]{ADAMTT}. In Sections~\ref{section_overview_Tlog},~\ref{sec_dimensions_SST} and~\ref{sec_sstproof} we work in a $1$-sorted setting, whereas in Section~\ref{sec_dminimality_criterion} we work in a possibly many-sorted setting.

\medskip\noindent
Throughout, \emph{$A$-definable} has its usual meaning, whereas \emph{definable} means \emph{definable with arbitrary parameters} (as opposed to meaning \emph{$\varnothing$-definable}). 
In general we choose not to consider a monster model (except in Lemma~\ref{lem_existentialmat}) and instead opt to work in arbitrary models of a (complete) theory, taking (sufficiently saturated) elementary extensions as needed. 

\medskip\noindent
Given an elementary extension $\bm{M}\preccurlyeq\bm{M}^*$, and a definable set $X$ in $\bm{M}$, it is understood that $X^*$ denotes its realization in $\bm{M}^*$; moreover, given a tuple $\alpha$ in $M^*$, we regard $\operatorname{tp}(\alpha|M)$ as an ultrafilter on a boolean algebra of $M$-definable subsets, as opposed to a collection of formulas. 

\medskip\noindent
As is standard when working with pregeometries, we often use concatenation to denote a union of sets, e.g., $A\Psi\Delta_{\phi}\coloneqq A\cup\Psi\cup\Delta_{\phi}$. We refer to~\cite[C.1]{tent2012course} for the basic properties of pregeometries.

\section{Overview of \texorpdfstring{$T_{\log}$}{Tₗₒ₉}}\label{section_overview_Tlog}

\noindent
In this section we provide an overview of the theory $T_{\log}$ and establish a few basic facts and conventions needed for the remainder of the paper. Since every asymptotic couple we will encounter will be a model of $T_{\log}$, we forgo a systematic development and instead define $T_{\log}$ directly from the standard model introduced in Section~\ref{section_introduction}. See~\cite[Sections~6.5, 9.2]{ADAMTT} for a general treatment of asymptotic couples (including the definition), and~\cite{GehretACTlog} for an axiomatic treatment of $T_{\log}$.

\subsection{More functions on the standard model} First, since the ordered abelian group $\Gamma_{\log}$ is divisible, for $n\geq 1$ we may define:
\[
\delta_n \colon \Gamma_{\log}\to\Gamma_{\log},\quad \alpha \ \mapsto \ \delta_n(\alpha) \ \coloneqq \ \textstyle\frac{1}{n}\alpha.
\]
Next, observe that $(\Gamma_{\log},\psi)$ has \emph{asymptotic integration}, i.e., for every $\alpha\in\Gamma_{\log}$ there exists a unique $\beta\in\Gamma_{\log}^{\neq}$ such that $\beta+\psi(\beta)=\alpha$. 
Thus, we may define the \emph{asymptotic integral}:
\[
\textstyle \int \colon \Gamma_{\log}\to \Gamma_{\log}^{\neq},\quad \alpha \ \mapsto \ \int\alpha \ \coloneqq \ \text{the unique $\beta$ such that $\beta+\psi(\beta)=\alpha$}.
\]
Note that $\int$ is an increasing bijection.
Next we define the \emph{successor function}:
\[
s \colon \Gamma_{\log}\to\Psi,\quad \alpha \ \mapsto \ s(\alpha) \ \coloneqq \ \textstyle \psi(\int\alpha).
\]
Note that $s0$ is the smallest element of $\Psi$. For each $\phi\in\Psi$, the element $s\phi$ is the immediate successor of $\phi$ in the ordered set $(\Psi;<)$, so $s$ restricts to an increasing bijection $s\colon\Psi\to\Psi^{>s0}$. Thus, we may define the \emph{predecessor function}:
\[
p \colon \Psi^{>s0}\to\Psi,\quad \phi \ \mapsto \ p(\phi) \ \coloneqq \ \text{the unique $\xi\in\Psi$ such that $s\xi=\phi$}.
\]
Finally, we adjoin to the underlying set $\Gamma_{\log}$ a new element $\infty$ and extend the ordering on $\Gamma_{\log}$ to $\Gamma_{\log}\cup\{\infty\}$ by declaring $\Gamma_{\log}<\infty$. Likewise, we extend the domains of the functions  defined so far for $\Gamma_{\log}$ by declaring $\infty$ to be a default value, i.e., for every $n\geq 1$, $\alpha\in\Gamma_{\log}\cup\{\infty\}$, and $\beta\in(\Gamma_{\log}\cup\{\infty\})\setminus\Psi^{>s0}$
we define:
\[
\textstyle
-\infty \ = \ \alpha+\infty \ = \ \infty+\alpha \ = \ \psi(0) \ = \ \psi(\infty) \ = \ \delta_n(\infty) \ = \ \int\infty \ = \ s(\infty) \ = \ p(\beta) \ \coloneqq \ \infty.
\]

\subsection{The \texorpdfstring{$\calL_{\log}$}{Lₗₒ₉}-theory \texorpdfstring{$T_{\log}$}{Tₗₒ₉}} We construe the standard model $(\Gamma_{\log},\psi)$ as an $\calL_{\log}$-structure, where
\[
\calL_{\log} \ \coloneqq \ \{0,-,+,<,\psi,\infty,(\delta_n)_{n\geq 1},s,p\}
\]
and define $T_{\log}\coloneqq\operatorname{Th}_{\calL_{\log}}(\Gamma_{\log},\psi)$.
By convention, we will always denote a model of $T_{\log}$ by a pair $(\Gamma,\psi)$, where it is understood that:
\begin{itemize}
\item $\Gamma=(\Gamma;0,-,+,<,(\delta_n)_{n\geq 1})$ is an ordered divisible abelian group, which we also regard as an ordered $\Q$-vector space;
\item $\psi$ is a function $\Gamma^{\neq}\to\Gamma$;
\item the primitives $s,p$ are left implicit as they can be defined in terms of $\psi$;
\item the underlying set of the $\calL_{\log}$-structure is $\Gamma_{\infty}\coloneqq\Gamma\cup\{\infty\}$, where $\Gamma<\infty$ and all primitives may be regarded as total functions;
\item although not part of the language, we can make use of the definable function $\int\colon\Gamma\to\Gamma^{\neq}$, defined in the same way as for $(\Gamma_{\log},\psi)$. \end{itemize}

\begin{remark}[Disclaimer about $\infty$] As $\infty$ is an element of our underlying structure, we treat it as such in model-theoretic statements. However, most of our statements are made only considering $\Gamma^n$ (without $\infty$). In this way we can use $\Q$-vector space arguments without having to accommodate $\infty$. Such statements can be readily adapted to  $\Gamma_{\infty}^n$. 
\end{remark}

\noindent
\emph{For the rest of Section~\ref{section_overview_Tlog}, we assume $(\Gamma,\psi)$ is a model of $T_{\log}$.}

\medskip\noindent
We need the following identities, established in \cite[Lemmas~3.2, 3.7, and 3.4]{GehretACTlog}.
The first connects $s$ and $\int$ and implies we may regard $\int$ as an $\calL_{\log}$-term:

\begin{fact}\label{fact_identities}
For all $\alpha,\beta\in\Gamma$:
\begin{enumerate}[itemindent=8em]
    \item[Integral Identity.] $\textstyle{\int}\alpha=\alpha-s\alpha$.
    \item[Fixed Point Identity.] $\beta=\psi(\alpha-\beta)$ if and only if $\beta=s\alpha$.
    \item[Successor Identity.] If $s\alpha<s\beta$, then $\psi(\beta-\alpha)=s\alpha$.
\end{enumerate}
\end{fact}

\begin{fact}[{QE~\cite[Theorem~5.2]{GehretACTlog} and UA~\cite[Lemma~5.1]{GehretACTlog}}]
$T_{\log}$ has quantifier elimination (QE) and a universal axiomatization (UA).
\end{fact}

\noindent
Here is an immediate consequence of QE and UA:

\begin{lemma}\label{easy_consequences_QE_UA}
Suppose $(\Gamma_1,\psi_1)\models T_{\log}$ and $(\Gamma_0,\psi_0)\subseteq(\Gamma_1,\psi_1)$ is an $\calL_{\log}$-substructure. Then $(\Gamma_0,\psi_0)\preccurlyeq(\Gamma_1,\psi_1)$, 
and thus $\psi_0(\Gamma_0^{\neq})=\Gamma_0\cap\psi_1(\Gamma^{\neq}_1)$.
\end{lemma}

\noindent
Given a tuple $\alpha$ in an extension $(\Gamma_1,\psi_1)\succcurlyeq (\Gamma,\psi)$, we let $\Gamma\langle \alpha\rangle$ denote the $\calL_{\log}$-substructure of $\Gamma_1$ generated by $\Gamma$ and $\alpha$. Then $(\Gamma,\psi)\preccurlyeq (\Gamma\langle \alpha\rangle,\psi_1|_{\Gamma\langle \alpha\rangle})\preccurlyeq (\Gamma_1,\psi_1)$ by the previous lemma.

\subsection{The \texorpdfstring{$\Psi$}{Ψ}-set of a model}
Given a model $(\Gamma_0,\psi_0)$ of $T_{\log}$, we set $\Psi_{\Gamma_0}\coloneqq\psi(\Gamma^{\neq}_0)$, which we refer to as the \emph{$\Psi$-set} of $(\Gamma_0,\psi_0)$. Thus by Lemma~\ref{easy_consequences_QE_UA}:
\[
\Psi_{\Gamma_0} \ = \ \Gamma_0\cap\Psi_{\Gamma_1} \quad \text{for any model $(\Gamma_1,\psi_1)$ of $T_{\log}$ extending $(\Gamma_0,\psi_0)$}.
\]
When a distinguished model $(\Gamma,\psi)$ is clear from context (as is currently the case), we will denote $\Psi_{\Gamma}$ as just $\Psi$. Finally, we set $\Psi_{\infty}\coloneqq\Psi\cup\{\infty\}$. 
\begin{lemma}\label{lem_elemsub}
Let $\Gamma_1$ be a $\Q$-linear subspace of $\Gamma$ with $\Psi\subseteq\Gamma_1$.
Then $(\Gamma_1,\psi|_{\Gamma_1})$ is an elementary substructure of $(\Gamma,\psi)$ with $\Psi_{\Gamma_1}=\Psi$.
\end{lemma}
\begin{proof}
From $\Psi\subseteq\Gamma_1$, we get that $\Gamma_{1,\infty}$ is closed under the functions $\psi$, $s$, and $p$.
Hence $(\Gamma_1,\psi|_{\Gamma_1})$ is an $\calL_{\log}$-substructure of $(\Gamma,\psi)$. The claim now follows by Lemma~\ref{easy_consequences_QE_UA}.
\end{proof}

\noindent
Observe that $(\Psi;<)\equiv(\N;<)$. In fact, this entirely characterizes the structure of the $\Psi$-set:

\begin{fact}[{\cite[Corollary~7.2]{GehretACTlog}}]\label{fact_Psistabembed}
If $(\Gamma,\psi)\models T_{\log}$, the structure $(\Psi;<)$ is purely stably embedded in $(\Gamma,\psi)$ in the sense that the structure induced on $\Psi$ by $(\Gamma,\psi)$ is just its structure as an ordered set.
\end{fact}

\begin{lemma}\label{definable_subset_Psi_set_fact}
If $C\subseteq\Psi$ is a definable convex subset of $\Psi$, then $C$ is either of the form $[\alpha,\beta]_{\Psi}$ or $\Psi^{\geq\alpha}$ for some $\alpha,\beta\in\Psi$.
\end{lemma}
\begin{proof}
This follows from Fact~\ref{fact_Psistabembed}, and can also be seen directly by noting that ``proper $s$-cuts'' (see~\cite[2.7]{GehretACTlogNIP}) are not definable since this would violate the universal property of the extension lemma~\cite[4.12]{GehretACTlog}.
\end{proof}

\begin{lemma}[Locally constant primitives]\label{loc_const_primitives}
When $\Gamma$ is equipped with the order topology we have:
\begin{enumerate}
\item $p$ takes constant value $\infty$ on $\Gamma\setminus\Psi^{>s0}$,
\item $\psi$ is locally constant on $\Gamma^{\neq}$,
\item $s$ is locally constant on $\Gamma$.
\end{enumerate}
\end{lemma}
\begin{proof}
(1) Clear from the definition. (2) Suppose $x\neq 0$, say $x>0$. Then $\psi(x/2)=\psi(2x)$, and so $\psi$ takes constant value $\psi(x)$ on the interval $(x/2,2x)$ since $\psi$ is a convex valuation. (3) Note that $\int\colon\Gamma\to\Gamma^{\neq}$ is a homeomorphism since it is a strictly increasing bijection. Thus, the composition $s=\psi\circ\int$ is locally constant by (2).
\end{proof}

\subsection{Quotienting by \texorpdfstring{$\Delta_{\phi}$}{∆ᵩ}}\label{subsec_quotienting}
\emph{In this subsection we let $\phi$ range over $\Psi_{\infty}$.} Recall:

\begin{definition}
For $\phi\in\Psi$ define the proper convex subgroup of $\Gamma$:
\[
\Delta_{\phi} \ \coloneqq \ \{x\in\Gamma:\psi(x)>\phi\}.
\]
We extend this to $\phi\in\Psi_{\infty}$ by setting $\Delta_{\infty}\coloneqq\{0\}$.
\end{definition}

\begin{lemma}\label{lem_char_def_convex_subgroups}
The set of definable convex subgroups of $\Gamma$ is $\{\Delta_{\phi}:\phi\in\Psi_{\infty}\}\cup\{\Gamma\}$.
\end{lemma}
\begin{proof}
Suppose $\Delta\subseteq\Gamma$ is a definable convex subgroup and consider $X\coloneqq\{\phi\in\Psi:\Delta_{\phi}\supseteq\Delta\}$, a definable initial segment of $\Psi$; by Lemma~\ref{definable_subset_Psi_set_fact} there are three cases to consider. If $X=\Psi$, then $\Delta=\{0\}=\Delta_{\infty}$. 
If $X\neq\Psi$ is nonempty, then for $\phi\coloneqq\max X$, we have $\Delta_{\phi}\supseteq \Delta\supsetneq\Delta_{s\phi}$. We claim $\Delta_{\phi}=\Delta$. Otherwise, since $\Delta$ is definable, this would contradict the universal property in the extension lemma~\cite[4.6]{GehretACTlog} since there would be two nonisomorphic ways of adding an element at the cut ``$\Delta^+$'': one that adds a new archimedean class to the extension of the definable set $\Delta$, and one that does not.
Finally, if $X=\varnothing$, then we have $\Delta=\Gamma$ by a similar argument.
\end{proof}

\noindent
Finally, we establish the following notation and conventions with regard to quotienting by $\Delta_{\phi}$:

\begin{itemize}
\item Given the subgroup $\Delta_{\phi}\subseteq\Gamma$, let $\pi_{\phi}\colon\Gamma\to\Gamma/\Delta_{\phi}$ denote the projection map. When $\phi$ is understood from context, we let $\bar{\Gamma}\coloneqq\Gamma/\Delta_{\phi}$. Since $\Delta_{\phi}$ is a convex $\Q$-subspace, we may construe $\bar{\Gamma}$ as an ordered $\Q$-vector space, with ordering induced by the ordering on $\Gamma$.
\item Given the subgroup $\Delta_{\phi}^n\subseteq\Gamma^n$, we identify $\Gamma^n/\Delta^n_{\phi}$ with $(\Gamma/\Delta_{\phi})^n$, which we also denote by $\bar{\Gamma}^n$ when $\phi$ is understood from context. We also have the natural projection map
\[
\pi_{\phi} \colon \Gamma^n\to\bar{\Gamma}^n,\quad (\alpha_1,\ldots,\alpha_n) \ \mapsto \ (\alpha_1+\Delta_{\phi},\ldots,\alpha_n+\Delta_{\phi}) \ = \ (\alpha_1,\ldots,\alpha_n)+\Delta_{\phi}^n.
\]
\item Given $\alpha\in\Gamma^n$ and $X\subseteq\Gamma^n$, we denote $\pi_{\phi}(\alpha)$ and $\pi_{\phi}(X)$ by just $\bar{\alpha}$ and $\bar{X}$ when $\phi$ is understood from context.
\end{itemize}

\begin{remark}\label{remark_phitopdef}
We define the \emph{$\phi$-topology} to be the coarsest topology on $\Gamma$ making $\pi_{\phi} \colon \Gamma \to \bar{\Gamma}$ continuous, where $\bar{\Gamma}$ is equipped with the order topology.
A basis for this topology is given by the sets of the form $\pi^{-1}_{\phi}((\bar{\alpha},\bar{\beta}))$, where $\alpha,\beta \in \Gamma$ with $\beta-\alpha>\Delta_{\phi}$.
Note that the $\infty$-topology equals the order topology on $\Gamma$, while for $\phi<\infty$, the $\phi$-topology is strictly coarser than the order topology on $\Gamma$ and does not satisfy the separation axiom $T_0$.
We define analogously the \emph{$\phi$-topology} on $\Gamma^n$, which equals the product topology on $\Gamma^n$ coming from the $\phi$-topology on $\Gamma$.

To reduce confusion, we confine the explicit use of the $\phi$-topologies to Remark~\ref{remark_phitop}, which indicates the statements that can be rewritten to incorporate the $\phi$-topologies, and Corollary~\ref{dim_equals_dim_closure}.
\end{remark}

\section{Dimension, the Small Sets Theorem, and related results}\label{sec_dimensions_SST}

\noindent
This section contains the main results of this paper.
First, Subsection~\ref{subsec_dimfun} reviews the definition and main properties of a \emph{dimension function} on a structure, and makes some basic observations about how dimension functions on models of $T_{\log}$ must behave.
Next, Subsection~\ref{subsec_pregeo} introduces the family of pregeometries on models of $T_{\log}$ that will later give rise to the family of dimension functions.
Subsection~\ref{subsec_sststatement} contains the statement of the Small Sets Theorem~\ref{smallSetsTheorem_label} (proof deferred to Section~\ref{sec_sstproof}).
Subsection~\ref{subsec_dimthmproof} shows how the Dimension Theorem~\ref{dimensionTheorem_label} follows from (part of) the Small Sets Theorem.
In Subsection~\ref{KSST_VSST_subsection} we see how the Small Sets Theorem splits into the Kinda Small Sets Theorem~\ref{smallishsets} and Very Small Sets Theorem~\ref{corverysmallsets} when considering definable sets of higher arity. We also observe some consequences, including local linearity (Corollary~\ref{cor_locallin}) and d-minimality (Corollary~\ref{corfornasierodmin}).
Finally, Subsection~\ref{subsec_EIfailure} shows the failure of elimination of imaginaries, as a consequence of dimension theory.

\medskip\noindent
\emph{In this section $(\Gamma,\psi)$ ranges over models of $T_{\log}$, and $\phi$ ranges over $\Psi_{\infty}$, where $\Psi=\psi(\Gamma^{\neq})$. We equip $\bar{\Gamma}=\Gamma/\Delta_{\phi}$ with the order topology and $\bar{\Gamma}^n=(\Gamma/\Delta_{\phi})^n=\Gamma^n/\Delta_{\phi}^n$ with the product topology. Given a set $X\subseteq\Gamma^n$, we construe $\bar{X}=\pi_{\phi}(X)$ as a subset of the ambient topological space $\bar{\Gamma}^n$.}

\subsection{Dimension functions}\label{subsec_dimfun}
We adopt the following definition of a \emph{dimension function} from \cite{vdD89}; below we declare $-\infty<\N$ and set $n+(-\infty)\coloneqq -\infty$ for every $n\in\N$.

\begin{definition}\label{def_dimfn}
Let $\bm{M}$ be a $1$-sorted structure.
A \textbf{dimension function} on $\bm{M}$ is a function $\d$ from the definable subsets of $M^n$ ($n$ varying) to $\N\cup\{-\infty\}$ such that:
\begin{enumerate}[label=(D\arabic*)]
    \item\label{D1}
        \begin{enumerate}[label=(\alph*)]
            \item\label{D1a} $\d(S)=-\infty \Leftrightarrow S=\0$ for definable $S \subseteq M^n$;
            \item\label{D1b} $\d(\{a\})=0$ for all $a \in M$;
            \item\label{D1c} $\d(M)=1$;
        \end{enumerate}
    \item\label{D2} $\d(S_1\cup S_2)=\max\{\d(S_1), \d(S_2)\}$ for definable $S_1,S_2 \subseteq M^n$;
    \item\label{D3} $\d$ is preserved under permutation of coordinates;
    \item\label{D4} if $S \subseteq M^{n+1}$ is definable, then $B_i\coloneqq \{ a\in M^n : \d(S_a)=i \}$ is definable and
    \[
    \d\big(\{(a,b)\in S : a\in B_i \}\big)\ =\ i + \d(B_i) \quad\text{for $i=0,1$.}
    \]
\end{enumerate}
\end{definition}

\noindent
These four axioms have several natural consequences. Here we collect a few that we need.

\begin{fact}[{\cite[1.1, 1.3, 1.5]{vdD89}}, {\cite[2.2]{angel2016bounded}}]\label{fact_vdD89}
Let $\bm{M}$ be a structure and $\d$ be a dimension function on $\bm{M}$.
\begin{enumerate}
    \item\label{dim_fact_subset} If $S_1,S_2\subseteq M^n$ are definable and $S_1\subseteq S_2$, then $\d S_1 \leq \d S_2$.
    \item If $S\subseteq M^n$ is finite and nonempty, then $\d S=0$.
    \item $\d M^n=n$.
    \item If $\d'$ is a dimension function on $\bm{M}$ such that $\d' S=\d S$ for all definable $S\subseteq M$, then $\d'=\d$.
    \item
    If $S \subseteq M^m$ and $f\colon S \to M^n$ are definable, then:
    \begin{enumerate}
        \item\label{dim_funct_fact_product} $\d(S\times T) = \d S+\d T$ for any definable $T \subseteq M^n$;
        \item\label{dim_funct_fact_def_bij} $\d S \geq \d f(S)$ (in particular, $\d S = \d f(S)$ for injective $f$);
        \item\label{dim_funct_fact_map_fibers} $B_i \coloneqq \{ a \in M^n : \d f^{-1}(a)=i \}$ is definable and $\d f^{-1}(B_i)=i+\d B_i$ for $i=0,\dots,m$.
    \end{enumerate}
    \item\label{dim_funct_fact_proj} For $S\subseteq M^n$ definable and $d \leq n$, we have $\d S\geq d$ if and only if $\d\pi(S) = d$ for some coordinate projection $\pi\colon M^n\to M^d$.
\end{enumerate}
\end{fact}

\begin{lemma}\label{dim_union_fact_linear_order}
Let $\operatorname{d}$ be a dimension function on an expansion of a linear order $(M;<,\ldots)$ and $X\subseteq M$ be a definable subset such that for every $a\in X$, we have $\operatorname{d}X^{\leq a}=0$. Then $\operatorname{d}X=0$.
\end{lemma}
\begin{proof}
Clearly, $\operatorname{d}X \in \{0,1\}$, so $\operatorname{d}X^2 \in \{0,2\}$ by Fact~\ref{fact_vdD89}\eqref{dim_funct_fact_product}. However, by \ref{D4} the ``triangle'' $\{(a,b)\in X^2:b\leq a\}$ has dimension $\d X$ since the vertical fibers $X_a=X^{\leq a}$ have dimension $0$. Since the square $X^2$ can be covered by two definably bijective copies of this triangle, it follows from Fact~\ref{fact_vdD89}\eqref{dim_funct_fact_def_bij} and~\ref{D2} that $\operatorname{d}X^2 \in \{0,1\}$, hence $\d X^2=0$. Thus $\operatorname{d}X=0$.
\end{proof}

\begin{cor}\label{cor_dPsi0}
Let $\operatorname{d}$ be a dimension function on $(\Gamma,\psi)$. Then $\operatorname{d}\Psi=0$.
\end{cor}
\begin{proof}
Consider $X\coloneqq\{\alpha \in \Psi:\d[s0,\alpha]_{\Psi}=0\}$, which is a definable initial segment of $\Psi$ by \ref{D4} and Fact~\ref{fact_vdD89}\eqref{dim_fact_subset}.
By Lemma~\ref{dim_union_fact_linear_order}, it suffices to show that $X=\Psi$.
By \ref{D1} and \ref{D2}, $s0 \in X$ and if $\alpha \in X$, then $s\alpha \in X$.
Hence $X=\Psi$ by Lemma~\ref{definable_subset_Psi_set_fact}, as desired.
\end{proof}

\noindent
In an ordered abelian group, there is always a largest definable convex subgroup with dimension~$0$:

\begin{lemma}\label{lem_odag_convex_dim_fact}
Let $\d$ be a dimension function on an expansion $(G;<,+,\ldots)$ of an ordered  abelian group.
Then $G$ has a largest definable convex subgroup $\Delta$ such that $\d\Delta=0$.
For an interval $(a,b)\subseteq G$, we have $\d(a,b)=1$ if and only if $b-a>\Delta$.
\end{lemma}
\begin{proof}
Define $\Delta\coloneqq\{ x \in G : \d[-|x|,|x|]=0 \}$.
The set $\Delta$ is convex by Fact~\ref{fact_vdD89}\eqref{dim_fact_subset} and definable by \ref{D4}.
It follows by Fact~\ref{fact_vdD89}\eqref{dim_funct_fact_def_bij} and \ref{D2} that if $x \in \Delta$, then $2x\in\Delta$, so $\Delta$ is a subgroup.
Applying Lemma~\ref{dim_union_fact_linear_order} to $\Delta^{\geq}$ yields $\d\Delta=0$.
The final statement follows by similar calculations.
\end{proof}

\noindent
The following consequence of Lemmas~\ref{lem_char_def_convex_subgroups} and~\ref{lem_odag_convex_dim_fact} limits the possible dimension functions:

\begin{cor}\label{cor_dphiwide}
If $\operatorname{d}$ is a dimension function on $(\Gamma,\psi)$, then there exists a unique $\phi$ such that $\operatorname{d}\Delta_{\phi}=0$ and $\operatorname{d}\Delta_{\xi}=1$ for every $\xi\in\Psi^{<\phi}$; moreover, for an interval $(a,b)\subseteq\Gamma$, we have $\operatorname{d}(a,b)=0$ if and only if $b-a\in\Delta_{\phi}$.
\end{cor}

\noindent
Corollary~\ref{cor_dphiwide} notwithstanding, the following two things are still not clear at this point:
\begin{itemize}
\item (Uniqueness) For each $\phi$, there is \emph{at most one} dimension function $\operatorname{d}$ such that $\operatorname{d}\Delta_{\phi}=0$ and $\operatorname{d}\Delta_{\xi}=1$ for $\xi\in\Psi^{<\phi}$.
\item (Existence) For each $\phi$, there is \emph{at least one} such dimension function.
\end{itemize}

\begin{remark}
If $\bm{R}=(R;<,+,\ldots)$ is an o-minimal expansion of an ordered group, then Lemma~\ref{lem_odag_convex_dim_fact} already provides the \emph{uniqueness} of a dimension function as an easy consequence of the o-minimality axiom (it implies all intervals must have dimension $1$), although \emph{existence} of a dimension function requires more work (via the Cell Decomposition Theorem, or by proving $\operatorname{acl}$ is a pregeometry as a consequence of the Monotonicity Theorem).
\end{remark}

\subsection{Pregeometries}\label{subsec_pregeo}
For each $\phi$ define a closure operator $\cl_{\phi}$ on $\Gamma_{\infty}$ by setting for each $A\subseteq \Gamma_{\infty}$:
\[
\cl_{\phi}(A)\ \coloneqq\ \acl(A\Psi\Delta_{\phi}).
\]
Although the model-theoretic algebraic closure $\acl$ in $T_{\log}$ is not a pregeometry, it follows from QE and UA that $\cl_{\phi}(A) \setminus\{\infty\}=\Qspan(A\Psi\Delta_{\phi})=\Qspan(A\Psi)+\Delta_{\phi}$ for $A\subseteq\Gamma$, so $\cl_{\phi}$ is a pregeometry; in particular, $\cl_{\phi}(A)=\cl_{\infty}(A)+\Delta_{\phi}$.
Moreover, $(\cl_{\phi}(A),\psi|_{\cl_{\phi}(A)})$ is an elementary substructure of $(\Gamma,\psi)$ with $\Psi_{\cl_{\phi}(A)}=\Psi$ by Lemma~\ref{lem_elemsub}.

\medskip\noindent
Each pregeometry yields a notion of dimension defined as follows.
Let $\rk_{\phi}(B|A)$ be the size of a basis of $\cl_{\phi}(AB)$ over $\cl_{\phi}(A)$. Then for a definable $X\subseteq\Gamma_{\infty}^n$, set 
\[
\dim_{\phi}(X)\ \coloneqq\ \sup\{\rk_{\phi}(\{x_0,\dots,x_{n-1}\}|\Gamma) : (x_0,\dots,x_{n-1}) \in X^* \}  \ \in \ \{-\infty,0,1,\ldots,n\},
\]
where $(\Gamma^*,\psi^*) \succcurlyeq (\Gamma,\psi)$ is $|\Gamma|^+$-saturated and $\rk_{\phi}$ is computed using the pregeometry of $(\Gamma^*,\psi^*)$.
This is independent of the choice of $|\Gamma|^+$-saturated $(\Gamma^*,\psi^*) \succcurlyeq (\Gamma,\psi)$ (see \cite[Section~2]{angel2016bounded} for a general statement of this kind).
Thus, we have a family of dimensions such that if $\phi\leq\xi\in \Psi_{\infty}$, then $\dim_{\phi}(X)\leq\dim_{\xi}(X)$ for all definable $X \subseteq\Gamma_{\infty}^n$.

\medskip\noindent
These dimensions also fit into the framework of \cite{angel2016bounded}: In every model of $T_{\log}$, the pregeometry $\cl_{\infty}$ is defined by the collection of $\calL_{\log}$-formulas of the form
\[
\exists x_0\cdots x_{m-1}\ \big(\sum_{i=0}^{m-1} q_i\psi(x_i) + \sum_{j=0}^{n-1} r_jy_j = u\big),
\]
where $q_i, r_j \in \Q$ for $i=0,\dots,m-1$ and $j=0,\dots,n-1$.
For $\phi \in \Psi$, the pregeometry $\cl_{\phi}$ is defined in every model of $\Th_{\calL_{\log}\cup\{\phi\}}(\Gamma,\psi)$ by the collection of $\calL_{\log}\cup\{\phi\}$-formulas of the form
\[
\exists x_0\cdots x_{m-1} \exists z\ \big(\psi(z)>\phi \wedge \big(\sum_{i=0}^{m-1} q_i\psi(x_i) + \sum_{j=0}^{n-1} r_jy_j + z = u\big)\big).
\]

\medskip\noindent
It is not obvious that each $\dim_{\phi}$ is a dimension function on $(\Gamma,\psi)$ in the sense of Definition~\ref{def_dimfn}, although \ref{D1}\ref{D1a}, \ref{D1}\ref{D1b}, \ref{D2}, and \ref{D3} are easy. 

\begin{lemma}\label{lem_Gammanotsmall}
For any $\phi \in \Psi_{\infty}$, we have $\dim_\phi(\Gamma) = 1$, i.e., $\dim_{\phi}$ satisfies \ref{D1}\ref{D1c}.
\end{lemma}
\begin{proof}
Let $(\Gamma^*,\psi^*)$ be an elementary extension of $(\Gamma,\psi)$ containing an element $\alpha> \Gamma$. We claim that $\Gamma\langle \alpha \rangle = \Gamma \oplus \Q\alpha$. It is enough to show that $\Gamma \oplus \Q\alpha$ is closed under $\psi^*$, $s$, and $p$. Consider an element $\gamma+q\alpha$, where $\gamma \in \Gamma$ and $q \in \Q^{\neq}$. Since $[\alpha]>[\Gamma^{\neq}]$, we have $\psi^*(\gamma+q\alpha) = \psi^*(\alpha) = s0$. In particular, $s0 = \psi^*(\gamma+q\alpha-s0)$, so the Fixed Point Identity (Fact~\ref{fact_identities}) gives $s(\gamma+q\alpha) = s0$. We also see that $\gamma+q\alpha$ is not in $\psi^*((\Gamma \oplus \Q\alpha)^{\neq})$, so $p(\gamma+q\alpha) = \infty$.

Having established this claim, we see that
\[
\Psi_{\Gamma\langle \alpha \rangle}= \Psi,\qquad \{\beta \in \Gamma\langle \alpha \rangle:\psi^*(\beta)>\phi\} = \Delta_\phi,
\]
so $\cl_\phi(\varnothing)$, computed in $\Gamma\langle \alpha \rangle$, is contained in $\Gamma$. Thus $\rk_\phi(\alpha|\Gamma) = 1$, so $\dim_\phi(\Gamma) = 1$. 
\end{proof}

\medskip\noindent
In light of Lemma~\ref{lem_Gammanotsmall},
to get that $\dim_{\phi}$ is a dimension function it remains to establish \ref{D4}, which we do in Corollary~\ref{dimension_theorem_existence} using part of the Small Sets Theorem.
As a first step, compactness and \ref{D2} yield \eqref{smallequiv1}$\Leftrightarrow$\eqref{smallequiv5} of the Small Sets Theorem (or see \cite[Lemma~2.3]{angel2016bounded}).
\begin{lemma}\label{lem_dphi0equivdphiset}
Let $X\subseteq\Gamma$ be definable.
Then $\dim_{\phi} X \leq 0$ if and only if $X$ is covered by finitely many affine maps $\Psi^n\times\Delta_{\phi} \to \Gamma$.
\end{lemma}

\noindent
For the remainder of the paper, we call a definable subset $X \subseteq \Gamma$ \textbf{$\phi$-small} if $\dim_\phi X \leq 0$; equivalently, $X$ is $\phi$-small if it  is covered by finitely many affine maps $\Psi^n\times\Delta_{\phi} \to \Gamma$.

\subsubsection*{Connection to existential matroids}\label{existential_matroids_remark}
When working in a monster model, the pregeometries $\cl_{\phi}$ fit into the framework of \emph{existential matroids} from \cite{fornasiero-dimmatroiddense}.
We need existential matroids in the proof of Corollary~\ref{dim_directed_union_fact}, which is used in the proofs of \eqref{smallequiv8}$\Rightarrow$\eqref{smallequiv1} of the Small Sets Theorem and \eqref{ksmallequiv6}$\Rightarrow$\eqref{ksmallequiv1} of the Kinda Small Sets Theorem.
We also use them to prove Corollary~\ref{cor_localprop} and Corollary~\ref{cor_dimfibre}, which are included for their own sakes and not used later.

\begin{lemma}\label{lem_existentialmat}
Suppose $(\Gamma,\psi)$ is a monster model of $T_{\log}$, expanded to the language $\calL_{\log}\cup\{\phi\}$. Then $\cl_{\phi}$ is an existential matroid in the sense of~\cite{fornasiero-dimmatroiddense}.
\end{lemma}

\subsection{The Small Sets Theorem}\label{subsec_sststatement} 
We first consider an example of a typical ($\infty$-)small set.

\begin{figure}[h!]
\caption{A set $X$ with $X,X',X''\neq \varnothing$ and $X^{(3)} = \varnothing$. }
\label{CB3}
\begin{center}
\begin{tikzpicture}
\tikzmath{\x = .8; \y =.9; \z =28;
\w =1/(1-\y); \u = -3.5; \v = 1-\x;}

\draw[<->] (-7,0)--(7,0); 
\draw[<->] (-7,\u)--(7,\u); 

\draw (0,-.15)--(0,0);
\draw[thick] (0,-.15+\u)--(0,.2+\u);

\node at (0,-.35) {\small$0$}; 
\node at (0,-.35+\u) {\small$0$}; 

\node at (7,.5) {$\Gamma$};
\node at (7,.5+\u) {$\Gamma$};

\draw [decorate,decoration={brace,amplitude=5pt,mirror,raise=4ex}]
(-.05,0) -- (\z*\y*\v+\z*\y*\x^2*\v*\v+.05,0) node[midway,yshift=-3em]{$X = \{(sx-x)+(sy-y):x\neq y \in \Psi\}$};

\draw [decorate,decoration={brace,amplitude=5pt,mirror,raise=4ex}]
(-.05,\u) -- (\z*\y*\v+.05,\u) node[midway,yshift=-3em]{$X' = \{sx-x:x \in \Psi\}\cup\{0\}$};

\foreach \a in {0,1,2,3,4,5,6,7,8,9,10,11,12,13,14,15,16,17,18,19,20} {
\begin{scope}
\draw[thick] (\z*\y*\x^\a*\v, -.4*\y*\x^\a+\u)--(\z*\y*\x^\a*\v, \y*\x^\a+\u);
\foreach \b in {2,3,4,5,6,7,8,9,10,11,12,13} {
\begin{scope}
\draw (\z*\y*\x^\a*\v+\z*\y*\x^\b*\x^\a*\v*\v, -.4*\y*\x^\b*\x^\a)--(\z*\y*\x^\a*\v+\z*\y*\x^\b*\x^\a*\v*\v, \y*\x^\b*\x^\a);
\end{scope}
}
\end{scope}
}

\node at (-1.5,-4.5) {$X'' = \{0\}$};
\node at (-.6,-4) {\rotatebox[origin=c]{40}{$\longrightarrow$}};

\node at (0,-2) {\rotatebox[origin=c]{90}{$\Longleftarrow$}};

\end{tikzpicture}
\end{center}
\end{figure}

\begin{figure}[h!]
\caption{The image in the quotient $\Gamma/\Delta_\phi$ of the set in Figure~\ref{CB3} for the value $\phi=s^50$.} 
\label{CBquot}
\begin{center}
\begin{tikzpicture}
\tikzmath{\x = .8; \y =.9; \z =28;
\w =1/(1-\y); \v = 1-\x;}

\draw[<->] (-7,0)--(7,0); 

\draw[thick] (0,-.15)--(0,.2);

\node at (0,-.35) {\small$0$}; 

\node at (7,.5) {$\Gamma/\Delta_\phi$};

\draw [decorate,decoration={brace,amplitude=5pt,mirror,raise=4ex}]
(-.05,0) -- (\z*\y*\v+\z*\y*\v*\v*\x+.05,0)  node[midway,yshift=-3em]{$\bar{X}$};

\foreach \a in {0,1,2,3} {
\begin{scope}
\draw[thick]  (\z*\y*\x^\a*\v+\z*\y*\v*\v*\x^4,-.4*\y*\x^4)--(\z*\y*\x^\a*\v+\z*\y*\v*\v*\x^4, \y*\x^4);
\foreach \b in {0,1,2,3} {
 \ifthenelse{\a < \b}{
\begin{scope}
\draw (\z*\y*\x^\a*\v+\z*\y*\x^\a*\v*\v*\x^4/\x^\b, -.4*\y*\x^\a*\x^4/\x^\b)--(\z*\y*\x^\a*\v+\z*\y*\x^\a*\v*\v*\x^4/\x^\b, \y*\x^\a*\x^4/\x^\b);
\end{scope}
}}
\end{scope}
}

\end{tikzpicture}
\end{center}
\end{figure}
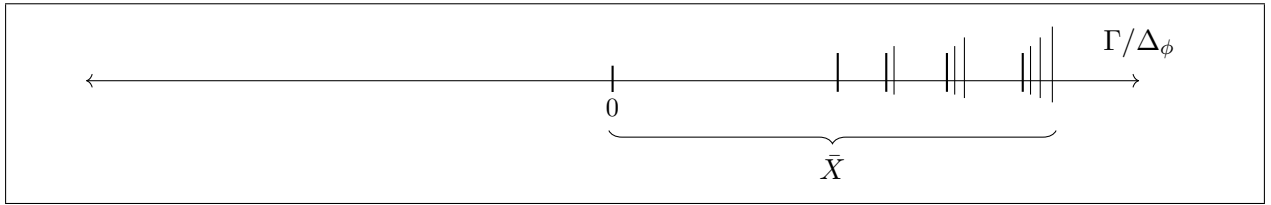

\medskip\noindent
Consider the definable set (pictured in Figure~\ref{CB3}):
\[
X \ \coloneqq \ \{(sx-x)+(sy-y):x\neq y\in\Psi\} \ \subseteq \ \Gamma.
\]
In the standard model, the set $X$ is countable and has the following explicit description:
\[
X \ = \ \{(\underbrace{0,\ldots,0}_{m},1,\underbrace{0,\ldots,0}_{n},1,0,0,\ldots):m\geq 1,n\geq 0\}.
\]
To see that $X$ is $\infty$-small using Lemma~\ref{lem_dphi0equivdphiset}, consider the affine map:
\[
F \colon \Psi^4\to\Gamma,\quad (x_0,x_1,x_2,x_3,) \ \mapsto \ (x_0-x_1)+(x_2-x_3)
\]
and note that $X\subseteq\operatorname{image}(F)$. Moreover, for the set:
\[
W \ \coloneqq \ \{(x_0,x_1,x_2,x_3):x_0=sx_1,x_2=sx_3,x_1< x_3\} \ \subseteq \ \Psi^4,
\]
we have a bijection $F|_W\colon W\to X$, where $W$ is definable in the structure $(\Psi;<)$. Next, observe that after taking the derived set finitely many times, we will arrive at $\varnothing$. Indeed:
\begin{align*}
X' \ &= \ \{sx-x:x\in\Psi\}\cup\{0\} \\
&= \ \{(\underbrace{0,\ldots,0}_{m},1,0,0,\ldots):m\geq 1\}\cup \{0\} \quad \text{(in the standard model)}
\end{align*}
and so $X''=\{0\}$ and thus $X^{(3)}=\varnothing$.
Next, consider the image $\pi_{s^n0}(X)$ in the quotient $\Gamma/\Delta_{s^n0}$, which we picture in Figure~\ref{CBquot}. Since in the standard model
\[
\Delta_{s^n0} \ = \ \{(\underbrace{0,\ldots,0}_n,r_n,r_{n+1},\ldots)\}
\]
we have:
\[
\Gamma/\Delta_{s^n0} \ = \ \{(r_0,\ldots,r_{n-1})\} \ \cong \ \R^n
\]
and thus $\pi_{s^n0}(X)$ is the set of $0/1$-vectors of length $n$ that begin with $0$ and contain $\leq 2$-many $1$'s. Observe that in this case $\pi_{\phi}(X)$ is finite, which indeed always happens when $\phi=s^n0$ for some $n$ (Corollary~\ref{phi_small_is_finite_sk0}). 

\medskip\noindent
Here is the general case.

\begin{smallSetsTheorem}\label{smallSetsTheorem_label}
For $\phi\in\Psi_{\infty}$ and $X\subseteq\Gamma$ definable in $(\Gamma,\psi)$, the following are equivalent:
\begin{enumerate}
\item\label{smallequiv1} $\dim_{\phi} X\leq 0$.
\item\label{smallequiv2} $\bar{X}$ is a finite union of discrete sets.
\item\label{smallequiv3} $\bar{X}$ has empty interior.
\item\label{smallequiv4} $X$ does not contain an interval $(a,b)$ such that $b-a>\Delta_{\phi}$.
\item\label{smallequiv5} $X$ can be covered by finitely many affine maps $\Psi^n\times\Delta_{\phi}\to\Gamma$.
\item\label{smallequiv6} $\bar{X}$ is d-finite.
\item\label{smallequiv7} $\bar{X}$ is nowhere dense \item\label{smallequiv8} $\bar{X}$ is definably meager in $\bar{\Gamma}$.
\item\label{smallequiv9} $X$ is $\Psi\Delta_{\phi}$-internal.
\item\label{smallequiv10} $X$ is $\Psi\Delta_{\phi}$-coanalysable.
\end{enumerate}
Furthermore, if $(\Gamma,\psi)$ is the standard model, then the above are additionally equivalent to:
\begin{enumerate}
\setcounter{enumi}{10}
\item\label{smallequiv11} $|\bar{X}|\leq\aleph_0$.
\item\label{smallequiv12} $|\bar{X}|<2^{\aleph_0}$.
\end{enumerate}
If $\phi<\infty$, then the above are additionally equivalent to:
\begin{enumerate}\setcounter{enumi}{12}
\item\label{smallequiv13} $\bar{X}$ is closed and discrete, and thus $\bar{X}'=\varnothing$.
\end{enumerate}
Finally, if $\phi=\infty$, then the above are equivalent to:
\begin{enumerate}\setcounter{enumi}{13}
\item\label{smallequiv14} $X$ is dp-finite.
\item\label{smallequiv15} The dp-rank of $X$ is less than $\aleph_0$.
\end{enumerate}
\end{smallSetsTheorem}

\noindent
See Subsections~\ref{subsec_roadmap},~\ref{interality_coanalysability_subsection} and~\ref{dp_rank_subsection} for precise definitions of \emph{definably meager}, \emph{internality}, \emph{coanalysability}, and \emph{dp-rank}.

\subsection{Proof of the Dimension Theorem from \texorpdfstring{\eqref{smallequiv1}$\Leftrightarrow$\eqref{smallequiv4}$\Leftrightarrow$\eqref{smallequiv5}}{(1)↔(4)↔(5)} of the Small Sets Theorem}\label{subsec_dimthmproof}

\begin{cor}[Uniqueness]\label{dimension_theorem_uniqueness}
Let $\operatorname{d}$ be a dimension function on $(\Gamma,\psi)$. Then $\operatorname{d}=\dim_{\phi}$ for some $\phi\in\Psi_{\infty}$.
\end{cor}
\begin{proof}
By Corollary~\ref{cor_dphiwide}, let $\phi$ be such that $\operatorname{d}\Delta_{\phi}=0$ and $\operatorname{d}\Delta_{\xi}=1$ for every $\xi\in\Psi^{<\phi}$. By Fact~\ref{fact_vdD89} it suffices to prove: for every definable $X\subseteq\Gamma_{\infty}$,
\begin{enumerate}[label=(\alph*)]
\item if $\dim_{\phi}X\leq 0$, then $\operatorname{d}X\leq 0$, and
\item if $\dim_{\phi}X=1$, then $\operatorname{d}X=1$.
\end{enumerate}
Let $X\subseteq\Gamma_{\infty}$ be an arbitrary definable set. We may assume $X\subseteq\Gamma$ since removing $\infty$ does not change whether dimension is $\leq0$ or $=1$. We have two cases:
\begin{itemize}
\item Suppose $\dim_{\phi}X\leq 0$. By~\ref{D2} and Lemma~\ref{lem_dphi0equivdphiset}, we further reduce to the case $X=\operatorname{image}(F)$, where $F\colon\Psi^n\times\Delta_{\phi}\to\Gamma$ is an affine map. Since $\operatorname{d}\Psi=0$ by Corollary~\ref{cor_dPsi0}, it follows from  Fact~\ref{fact_vdD89}(\ref{dim_funct_fact_product},\ref{dim_funct_fact_def_bij}) that $\operatorname{d}X\leq0$.
\item Suppose $\dim_{\phi}X=1$. Then by \eqref{smallequiv4}$\Rightarrow$\eqref{smallequiv1}, $X$ contains an interval $(a,b)$ where $b-a>\Delta_{\phi}$. Thus $\operatorname{d}(a,b)=1$ by Corollary~\ref{cor_dphiwide}, and thus $\operatorname{d}X=1$ by Fact~\ref{fact_vdD89}\eqref{dim_fact_subset}. \qedhere
\end{itemize}
\end{proof}

\begin{cor}[Existence]\label{dimension_theorem_existence}
For each $\phi\in\Psi_{\infty}$, the function $\dim_{\phi}$ is a dimension function.
\end{cor}
\begin{proof}
It remains to show that each $\dim_{\phi}$ satisfies \ref{D4}. 
Suppose $S\subseteq\Gamma_{\infty}^{n+1}$ is definable and set $B_i\coloneqq\{a\in\Gamma_{\infty}^n:\dim_{\phi}S_a=i\}$ for $i=0,1$. Then using \eqref{smallequiv1}$\Leftrightarrow$\eqref{smallequiv4} we see that $B_1$ is definable:
\begin{align*}
B_1 \ &= \ \{a\in\Gamma_{\infty}^n:\dim_{\phi}(S_a\setminus\{\infty\})=1\} \\
&= \ \{a\in\Gamma_{\infty}^n:\text{there exists an interval $(b_1,b_2)\subseteq S_a\setminus\{\infty\}$ such that $b_2-b_1>\Delta_{\phi}$}\}.
\end{align*}
Thus $B_0$ is definable as well. 
It now follows from standard pregeometry arguments that $\dim_{\phi}S_i=\dim_{\phi}B_i+i$, where $S_i\coloneqq\{(x,y)\in S:x\in B_i\}$. This is asserted 
in~\cite[Proposition 2.8(2)]{angel2016bounded}.
\end{proof}

\begin{cor}\label{dim_directed_union_fact}
Fix $d\in\{-\infty,0,\ldots,n\}$ and suppose $(I;\leq)$ is a definable directed set and $X\subseteq I\times\Gamma^n_{\infty}$ is a definable family such that $X_a\subseteq X_b$ for every $a\leq b$ in $I$. If $\dim_{\phi}X_a\leq d$ for every $a\in I$, then $\dim_{\phi}(\bigcup_{a\in I}X_a)\leq d$.
\end{cor}
\begin{proof}
Let $(\Gamma^*,\psi^*)$ be a monster model extending $(\Gamma,\psi)$. Then Corollary~\ref{dimension_theorem_existence} gives $\dim_{\phi}(X^*)_a\leq d$ for every $a\in I^*$. By Lemma~\ref{lem_existentialmat} and ~\cite[Lemma 3.71]{fornasiero-dimmatroiddense}, it follows that $\dim_{\phi}(\bigcup_{a\in I^*}(X^*)_a)\leq d$. Thus $\dim_{\phi}(\bigcup_{a\in I}X_a)\leq d$ since $(\bigcup_{a\in I}X_a)^*=\bigcup_{a\in I^*}(X^*)_a$. \end{proof}

\subsection{Higher arity definable sets}\label{KSST_VSST_subsection}

For definable subsets of $\Gamma^n$ with $n>1$, the conditions in the Small Sets Theorem are clearly not equivalent. However, each turns out to be equivalent to either having $\phi$-dimension $\leq 0$ (being ``very small'') or having $\phi$-dimension $< n$ (being ``kinda small''). In this subsection, we freely use both the Small Sets Theorem and the Dimension Theorem.

\begin{kindaSmallSetsTheorem}\label{smallishsets}
For $\phi\in\Psi_{\infty}$ and $X\subseteq\Gamma^n$ definable in $(\Gamma,\psi)$ the following are equivalent:
\begin{enumerate}
\item\label{ksmallequiv1} $\dim_{\phi} X<n$.
\item\label{ksmallequiv2} $\bar{X}$ has empty interior.
\item\label{ksmallequiv3} $X$ does not contain an open box $\prod_{i=1}^n(a_i,b_i)$ where $b_i-a_i>\Delta_{\phi}$ for each $i$.
\item\label{ksmallequiv4} $X$ can be covered by finitely many sets of the form 
\[
\{\alpha \in \Gamma^n:q\cdot \alpha \in Y\}
\]
where $q \in \Q^n$ is not the zero tuple and $Y \subseteq \Gamma$ is a definable $\phi$-small set \item\label{ksmallequiv5} $\bar{X}$ is nowhere dense.
\item\label{ksmallequiv6} $\bar{X}$ is definably meager in $\bar{\Gamma}^n$.
\end{enumerate}
\end{kindaSmallSetsTheorem}
\begin{proof}
Suppose $\dim_\phi X < n$, and let $\alpha$ be a tuple in an elementary extension $(\Gamma^*, \psi^*)$ of $(\Gamma,\psi)$ with $\alpha \in X^*$. Then $\rk_\phi(\alpha|\Gamma)<n$, so $\alpha$ is $\Q$-linearly dependent over
\[
\cl_{\phi}(\Gamma)\ =\ \Gamma + \Qspan(\Psi^*)+ \Delta_\phi^*.
\]
By a standard compactness argument, $X$ is covered by finitely many sets of the form 
\[
\{\alpha \in \Gamma^n:q\cdot \alpha \in Y\}
\]
where $q \in \Q^n$ is not the zero tuple and where $Y$ is the image of an affine map $\Psi^{m}\times\Delta_{\phi}\to\Gamma$. This gives \eqref{ksmallequiv1}$\Rightarrow$\eqref{ksmallequiv4}. 

We now show that \eqref{ksmallequiv4}$\Rightarrow$\eqref{ksmallequiv5}. Since nowhere dense sets form an ideal, we may assume that $X$ is of the form $\{\alpha \in \Gamma^n:q\cdot \alpha \in Y\}$. Then $\bar{X} \subseteq \{\bar{\alpha} \in \bar{\Gamma}^n:q\cdot \bar{\alpha} \in \bar{Y}\}$. 
We may further assume that $\bar{Y}$ is closed---when  $\phi < \infty$, this holds by \eqref{smallequiv13} of the Small Sets Theorem, and when $\phi = \infty$, then $\bar{Y} = Y$ is nowhere dense, so its closure is as well. Then $\bar{X}$ is closed, so we need only show that $\bar{X}$ has empty interior. Since the map $\bar{\alpha}\mapsto q\cdot \bar{\alpha}\colon \bar{\Gamma}^n\to \bar{\Gamma}$ is open, this follows from the fact that $\bar{Y}$ has empty interior. 

The implications \eqref{ksmallequiv5}$\Rightarrow$\eqref{ksmallequiv2}$\Rightarrow$\eqref{ksmallequiv3} are immediate. We establish the implication \eqref{ksmallequiv3}$\Rightarrow$\eqref{ksmallequiv1} in Proposition~\ref{prop_simpleextinterior} below.
Finally, having established \eqref{ksmallequiv5}$\Rightarrow$\eqref{ksmallequiv1}, the implications \eqref{ksmallequiv5}$\Rightarrow$\eqref{ksmallequiv6}$\Rightarrow$\eqref{ksmallequiv1} follow as in the $n=1$ case done in Stage (IV) of the proof of the Small Sets Theorem below, also using Corollary~\ref{dim_directed_union_fact}.
\end{proof}

\noindent
Using Fact~\ref{fact_vdD89}\eqref{dim_funct_fact_proj} and the equivalence \eqref{ksmallequiv1}$\Leftrightarrow$\eqref{ksmallequiv2} of the Kinda Small Sets Theorem, we obtain:
\begin{cor}[Coincidence with topological dimension]\label{agreement_topological_dimension}
The dimension function $\dim_\phi$ coincides with naive topological dimension in the quotient $\bar{\Gamma}$. That is, for $X \subseteq \Gamma^n$ and $d \leq n$, we have $\dim_\phi X \geq d$ if and only if $\pi(\bar{X})$ has nonempty interior for some coordinate projection $\pi\colon \bar{\Gamma}^n\to \bar{\Gamma}^d$.
\end{cor}

\begin{cor}[Local linearity]\label{cor_locallin}
Suppose $f\colon\Gamma^n\to\Gamma_{\infty}$ is definable. Then there exists a definable dense open set $U\subseteq\Gamma^n$ such that $f$ is locally affine at all $x \in U$. 
Moreover, there is a finite set $Q\in\Q^n$ that the ``slopes'' of the affine functions come from.
\end{cor}
\begin{proof}
As $T_{\log}$ has quantifier elimination and a universal axiomatization in the language $\calL_{\log}$, definable functions are given piecewise by terms. Thus, we may assume that $f$ is a term $\tau$. We proceed by induction on complexity of terms, assuming that for any term $\sigma$ less complex than $\tau$, there is a definable dense open set $U_\sigma$ on which $\sigma$ is locally affine (with slopes coming from a finite subset $Q_\sigma \subseteq \Q^n$). If $\tau$ is a sum $\sigma_1+\sigma_2$, then this still holds for $\tau$ on the dense open set $U_{\sigma_1}\cap U_{\sigma_2}$, and if   $\tau = -\sigma$ or $\delta_n(\sigma)$ for some $n\geq 1$, then this still holds for $\tau$ on $U_\sigma$, so we may assume that $\tau$ is either of the form $\psi \circ \sigma$, $s\circ \sigma$, or $p \circ \sigma$.

Let $V\subseteq U_\sigma$ be the set of $x \in \Gamma^n$ at which $\sigma$ is locally constant, so $\tau$ is locally constant on $V$ as well. Let $W \coloneqq U_\sigma \setminus V$, so $\sigma$ is locally affine and nonconstant on $W$. We will find a dense open subset of $W$ on which $\tau$ is locally constant. By Lemma~\ref{loc_const_primitives}, $s$ is locally constant on $\Gamma$, $\psi$ is locally constant on $\Gamma \setminus \{0\}$, and $p$ is locally constant (with constant value $\infty$) on $\Gamma \setminus \Psi^{>s0}$. Thus, it is enough to show that 
\[
X\ \coloneqq\ \{x \in W: \sigma(x) \in \{0\} \cup \Psi^{>s0}\}
\]
is nowhere dense.

For $q \in Q_{\sigma}$ and $\beta \in \Gamma$, let $W_{q,\beta}$ be the set of $x \in W$ such that $\sigma(y) = q\cdot y+\beta$ for $y$ in  an open  neighbourhood of $x$. Then $(W_{q,\beta})_{q \in Q_\sigma,\beta \in \Gamma}$ is a definable family of disjoint open subsets of $W$ with $W = \bigcup_{q,\beta}W_{q,\beta}$. Let 
\[
X_{q,\beta}\ \coloneqq\ \{x \in W_{q,\beta}:q\cdot x +\beta\in \{0\}\cup \Psi^{>s0}\},
\]
so $X = \bigcup_{q,\beta}X_{q,\beta}$ and each $X_{q,\beta}$ is nowhere dense by the equivalence \eqref{ksmallequiv4}$\Leftrightarrow$\eqref{ksmallequiv5} of the Kinda Small Sets Theorem. Since the sets $W_{q,\beta}$ are disjoint, we see that $X$ is also nowhere dense, as desired.
\end{proof}

\begin{cor}[Uniform definability of dimension]\label{uniform_definability_of_dimension}
Suppose $S\subseteq \Gamma^n$ and $f\colon S\to\Gamma^m$ are definable and fix $d\in\{-\infty,0,\ldots,n\}$. Then
\[
\{(\phi,b):\dim_{\phi}f^{-1}(b)=d\} \ \subseteq \ \Psi_{\infty}\times\Gamma^m
\]
is definable.
\end{cor}
\begin{proof}
For any $\phi$, we have $\dim_\phi f^{-1}(b) = -\infty$ if and only if $b \not\in f(S)$, so we may assume $d \geq 0$. It suffices to show that $\{(\phi,b):\dim_{\phi}f^{-1}(b)\geq d\}$ is definable. By Fact~\ref{fact_vdD89}\eqref{dim_funct_fact_proj} and the equivalence \eqref{ksmallequiv1}$\Leftrightarrow$\eqref{ksmallequiv3} of the Kinda Small Sets Theorem, we have $\dim_{\phi}f^{-1}(b)\geq d$ if and only if $\pi(f^{-1}(b))$ contains an open box $\prod_{i=1}^d(a_i,b_i)$ with $b_i-a_i>\Delta_\phi$ for some coordinate projection $\pi\colon \Gamma^n\to \Gamma^d$. This is a definable condition on $\phi$ and $b$. 
\end{proof}

\noindent
Corollary~\ref{uniform_definability_of_dimension} suggests considering the limits of $\dim_\phi$ and $\rk_\phi$ as $\phi$ increases in $\Psi$.
Note that $\bigcap_{\phi\in\Psi}\Delta_{\phi}=\{0\} = \Delta_{\infty}$. This allows us to compute the $\infty$-dimension of a definable set $X$ as a limit of its $\phi$-dimensions.

\begin{cor}[Continuity of dimension]
Suppose $X\subseteq\Gamma^n$ is definable. Then:
\[
\lim_{\phi\to\infty}\dim_{\phi}(X) \ = \ \dim_{\infty}(X).
\]
\end{cor}
\begin{proof}
Let $d \coloneqq \dim_\infty(X)$. Since $\dim_\phi(X)$ is increasing in $\phi$, it is enough to show that $\dim_\phi(X) \geq d$ for some $\phi \in \Psi$. By Corollary~\ref{agreement_topological_dimension}, there is a coordinate projection $\pi\colon \Gamma^n\to \Gamma^d$ and an open box $\prod_{i=1}^d(a_i,b_i)$ contained in $\pi(X)$. Take $\phi$ large enough with $b_i-a_i>\Delta_\phi$ for $i = 1,\ldots,d$. For this $\phi$, we have $\dim_\phi X \geq d$. 
\end{proof}

\noindent
This continuity of dimension can fail for the ranks of individual elements:

\begin{example}
Let $(\Gamma,\psi)$ be the standard model, and let $\alpha=(1,1/2,1/3,1/4,\ldots)$ be an element in an immediate elementary extension $(\Gamma^*,\psi^*)$ of $(\Gamma,\psi)$ (cf.~\cite[Example 2]{GehretACTlogNIP}). Then inside $(\Gamma^*,\psi^*)$ we have $\alpha\not\in\cl_{\infty}\Gamma$, although $\alpha\in\bigcap_{\phi\in\Psi}\cl_{\phi}\Gamma$. Thus, 
$\lim_{\phi\to\infty}\operatorname{rk}_\phi(\alpha|\Gamma)=0<1=\operatorname{rk}_{\infty}(\alpha|\Gamma)$.
\end{example}

\noindent
In connection with~\cite{fornasierohalupczok}, we record that our dimension $\dim_\phi$ is \emph{local} for each $\phi$ (see Remark~\ref{remark_phitop} for the connection with the $\phi$-topology defined in Remark~\ref{remark_phitopdef}):
\begin{cor}\label{cor_localprop}
A definable set $X\subseteq \Gamma^n$ has $\dim_{\phi} X\leq d$ if for every $x \in X$, there is a box $U = \prod_{i=1}^n(a_i,b_i)$ with $b_i-x_i, x_i-a_i>\Delta_\phi$ for each $i$ such that $\dim_\phi(X \cap U) \leq d$.
\end{cor}
\begin{proof}
This is an immediate consequence of~\cite[Theorem~3.1]{fornasierohalupczok} when $\phi = \infty$, but the proof goes through for all $\phi$.
Alternatively, here is an outline of a proof avoiding the use of the Erd\H{o}s--Rado Theorem and other machinery from \cite{fornasiero-dimmatroiddense} in the proof of \cite[Theorem~3.1]{fornasierohalupczok}.
Pass to a monster model.
For $x \in X$, assume that we have a box $U$ as in the statement. We then use the Small Sets Theorem recursively to shrink $U$ so that its endpoints are suitably independent from $x$ in the sense of \cite[Lemma~3.69]{fornasiero-dimmatroiddense}, which verifies the hypotheses of that lemma and so yields $\dim_{\phi} X \leq d$.
\end{proof}

\begin{verySmallSetsTheorem}\label{corverysmallsets}
For $\phi\in\Psi_{\infty}$ and $X\subseteq\Gamma^n$ definable in $(\Gamma,\psi)$ the following are equivalent:
\begin{enumerate}
\item\label{vsmallequiv1} $\dim_{\phi} X\leq 0$.
\item\label{vsmallequiv2} $X$ is contained in a product $\prod_{i=1}^n X_i$, where each $X_i\subseteq \Gamma$ is a definable $\phi$-small set.
\item\label{vsmallequiv3} $\bar{X}$ is a finite union of discrete sets.
\item\label{vsmallequiv4} $X$ can be covered by finitely many affine maps $\Psi^{m}\times\Delta_{\phi}\to\Gamma^n$.
\item\label{vsmallequiv5} 
$\bar{X}$ is d-finite.
\item\label{vsmallequiv6} $X$ is $\Psi\Delta_{\phi}$-internal.
\item\label{vsmallequiv7} $X$ is $\Psi\Delta_{\phi}$-coanalysable.
\end{enumerate}
Furthermore, if $(\Gamma,\psi)$ is the standard model, then the above are additionally equivalent to:
\begin{enumerate}
\setcounter{enumi}{7}
\item\label{vsmallequiv8} $|\bar{X}|\leq\aleph_0$.
\item\label{vsmallequiv9} $|\bar{X}|<2^{\aleph_0}$.
\end{enumerate}
If $\phi<\infty$, then the above are additionally equivalent to:
\begin{enumerate}\setcounter{enumi}{9}
\item\label{vsmallequiv10} $\bar{X}$ is closed and discrete, and thus $\bar{X}'=\varnothing$.
\end{enumerate}
Finally, if $\phi=\infty$, then the above are additionally equivalent to:
\begin{enumerate}\setcounter{enumi}{10}
\item\label{vsmallequiv11} $X$ is dp-finite.
\item\label{vsmallequiv12} The dp-rank of $X$ is less than $\aleph_0$.
\end{enumerate}
\end{verySmallSetsTheorem}
\begin{proof}
For \eqref{vsmallequiv1}$\Rightarrow$\eqref{vsmallequiv2}, let $X_i$ be the coordinate projection of $X$ onto the $i$th coordinate. Then $X \subseteq \prod_{i=1}^n X_i$ and $\dim_\phi X_i \leq \dim_\phi X \leq 0$ by Fact~\ref{fact_vdD89}\eqref{dim_funct_fact_def_bij}.
Suppose \eqref{vsmallequiv2} holds and for each $i \in \{1,\ldots,n\}$, take affine maps $F_{i,0},\ldots,F_{i,n_i} \colon \Psi^{m}\times\Delta_{\phi}\to\Gamma$ whose images cover $X_i$ (we arrange that $m$ is the same for each $F_{i,j}$). For $\bm{j} = (j_1,\ldots,j_n)\in \N^n$ with $j_i\leq n_i$ for each $i$, we let $F_{\bm{j}}\colon \Psi^{m}\times\Delta_{\phi}\to\Gamma^n$ be the affine map $(F_{1,j_1},F_{2,j_2},\ldots,F_{n,j_n})$. Then $X$ is covered by the images of the (finitely many) $F_{\bm{j}}$, establishing \eqref{vsmallequiv4}.

We have \eqref{vsmallequiv4}$\Rightarrow$\eqref{vsmallequiv6} by the  definition of internality (see Lemma~\ref{lem_dsetPsiint}) 
and \eqref{vsmallequiv6}$\Rightarrow$\eqref{vsmallequiv7} since internality always implies coanalysability (see Remark~\ref{rem_coanal}). Using Lemma~\ref{lem_coanaldphi0}, we obtain \eqref{vsmallequiv7}$\Rightarrow$\eqref{vsmallequiv1}, since  $\dim_{\phi}\Psi\Delta_{\phi}=0$. 

Take $X_1,\ldots,X_n \subseteq \Gamma$ as in \eqref{vsmallequiv2}, so the Small Sets Theorem gives natural numbers $m_1,\ldots,m_n$ such that $\bar{X}_i^{(m_i)} = \varnothing$ for each $i$. Using Lemma~\ref{CBproduct1} and induction on $n$, we get $\bar{X}^{(m)} = \varnothing$ for $m = m_1+\cdots+m_n$, so \eqref{vsmallequiv2}$\Rightarrow$\eqref{vsmallequiv5}, and we always have \eqref{vsmallequiv5}$\Rightarrow$\eqref{vsmallequiv3} (see Lemma~\ref{CBsmall_implies_finite_union_discrete}). 

Let $X\subseteq \Gamma^n$ and suppose that $\bar{X}$ is a finite union of discrete sets. 
To establish \eqref{vsmallequiv3}$\Rightarrow$\eqref{vsmallequiv1}, we need to show that $\dim_\phi X \leq 0$. We proceed by induction on $n$ (the case $n=1$ being the Small Sets Theorem). Since $\bar{X}$ is nowhere dense (Lemma~\ref{discrete_implies_nowhere_dense_under_mild_assumptions}),
 the Kinda Small Sets Theorem~\ref{smallishsets} tells us that $X$ can be covered by finitely many sets of the form
\[
\{\alpha \in \Gamma^n:q\cdot \alpha \in Y\}
\]
where $q \in \Q^n$ is a tuple of rational numbers, not all zero, and $Y\subseteq \Gamma$ is definable and $\phi$-small. Using \ref{D2}, we may assume that $X$ is contained in one of these sets. For $\gamma \in Y$, let $X_\gamma\coloneqq\{\alpha \in X:q\cdot \alpha = \gamma\}$, so $X_\gamma$ is the intersection of $X$ with the hyperplane $\{\alpha \in \Gamma^n:q\cdot \alpha =\gamma\}$, which is homeomorphic to $\Gamma^{n-1}$ via some coordinate projection  map $\pi\colon \Gamma^n\to \Gamma^{n-1}$. Clearly, $\overline{X_\gamma}$ is also a finite union of discrete sets, so $\pi(\overline{X_\gamma})= \overbar{\pi(X_\gamma)}$ is a finite union of discrete sets as well. Our induction hypothesis and Fact~\ref{fact_vdD89}\eqref{dim_funct_fact_def_bij} give $\dim_\phi X_\gamma = \dim_\phi \pi(X_\gamma)\leq 0$. Finally, using that  $X = \bigcup_{\gamma \in Y}X_\gamma$, that $\dim_\phi Y\leq 0$, and Fact~\ref{fact_vdD89}\eqref{dim_funct_fact_map_fibers}, we get that $\dim_\phi X \leq 0$. 

In the case that $\phi<\infty$, we obtain \eqref{vsmallequiv2}$\Rightarrow$\eqref{vsmallequiv10}, using that a product of closed discrete sets is closed and discrete. Then \eqref{vsmallequiv10}$\Rightarrow$\eqref{vsmallequiv5} trivially. If $\phi = \infty$, then \eqref{vsmallequiv2}$\Rightarrow$\eqref{vsmallequiv11} by Fact~\ref{dpfacts}, \eqref{vsmallequiv11}$\Rightarrow$\eqref{vsmallequiv12} is trivial, and \eqref{vsmallequiv12}$\Rightarrow$\eqref{vsmallequiv2} again by Fact~\ref{dpfacts}, taking $X_i$ to be the coordinate projection of $X$ onto the $i$th coordinate. Finally, the equivalence \eqref{vsmallequiv1}$\Leftrightarrow$\eqref{vsmallequiv2} and basic properties of cardinality give \eqref{vsmallequiv1}$\Leftrightarrow$\eqref{vsmallequiv8}$\Leftrightarrow$\eqref{vsmallequiv9} in the standard model.
\end{proof}

\begin{remark}\label{dim0_definable_in_pair}
Let $X \subseteq \Gamma^n$ be definable with $\dim_\infty X\leq 0$. Then $X$ is already definable in the reduct $(\Gamma,\Psi)$ by \eqref{vsmallequiv4} above, since the structure $(\Psi;<)$ is purely stably embedded. 
\end{remark}
\begin{remark}
The equivalence \eqref{vsmallequiv8}$\Leftrightarrow$\eqref{smallequiv9} above asserts that ``the continuum hypothesis holds for definable sets in the standard model''.
\end{remark}

\begin{cor}[d-minimality]\label{corfornasierodmin}
The theory $T_{\log}$ is d-minimal in the stronger sense defined by Fornasiero in~\cite[Definition 9.1]{fornasiero-dimmatroiddense}, i.e.:
\begin{enumerate}
\item If $X \subseteq \Gamma$ is definable with empty interior, then $X$ is a finite union of discrete sets. 
\item If $X \subseteq \Gamma^n$ is definable and discrete, then $\pi_1(X)$ has empty interior, where $\pi_1$ is the projection onto the first coordinate.
\item If $X \subseteq \Gamma^2$ and $U \subseteq \pi_1(X)$ are definable, $U$ is open and nonempty, and $X_a$ has nonempty interior for each $a \in U$, then $X$ has nonempty interior.
\end{enumerate}
\end{cor}
\begin{proof}
The first condition holds by the Small Sets Theorem. For $X$ as in the second condition, the Very Small Sets Theorem~\ref{corverysmallsets} gives $\dim_\infty\pi_1 X\leq \dim_\infty X \leq 0$, so $\pi_1 X$ has empty interior by the Small Sets Theorem. For $X,U$ as in the third condition, we have $\dim_\infty U = 1 = \dim_{\infty} X_a$ for all $a \in U$ by the Small Sets Theorem. Then $\dim_\infty X = 2$, so $X$ has nonempty interior by the Kinda Small Sets Theorem~\ref{smallishsets}. 
\end{proof}

\begin{cor}\label{phi_small_is_finite_sk0}
Let $X\subseteq\Gamma^n$ be definable and suppose $\phi=s^k0$ for some $k$. Then $\dim_\phi X \leq 0$ if and only if $\bar{X}$ is finite.
\end{cor}
\begin{proof}
The implication \eqref{vsmallequiv10}$\Rightarrow$\eqref{vsmallequiv1} of the Very Small Sets Theorem tells us that $\dim_\phi X \leq 0$ if $\bar{X}$ is finite. For the other implication, it suffices by 
\eqref{vsmallequiv1}$\Rightarrow$\eqref{vsmallequiv4} of the Very Small Sets Theorem to consider the case $X=\operatorname{image}(F)$, where $F\colon\Psi^m\times\Delta_{\phi}\to\Gamma^n$ is an affine map. Define $G\colon(\Psi^{\leq \phi})^m\to\Gamma^n$ by $G(\alpha)\coloneqq F(\alpha,0)$. Then $\operatorname{image}(G)$ is finite, and $\overline{\operatorname{image}(G)}=\bar{X}$.
\end{proof}

\begin{remark}\label{remark_phitop}
In connection with the $\phi$-topology defined in Remark~\ref{remark_phitopdef}, we remark the following.
\begin{itemize}
    \item We can replace \eqref{smallequiv4} in the Small Sets Theorem by the statement that $X$ has nonempty interior in the $\phi$-topology, and similarly with \eqref{ksmallequiv3} in the Kinda Small Sets Theorem. Hence, a definable $X \subseteq \Gamma^n$ has nonempty interior in the $\phi$-topology if and only if $\bar{X}$ has nonempty interior in $\bar{\Gamma}^n$, from which a version of Corollary~\ref{agreement_topological_dimension} for the $\phi$-topology follows. This is also needed to invoke \cite[Proposition~2.1]{fornasierohalupczok} in the proof of Corollary~\ref{dim_equals_dim_closure} below.
    \item Regarding Corollary~\ref{cor_localprop}, every neighbourhood of $x\in \Gamma^n$ in the $\phi$-topology contains a box $U = \prod_{i=1}^n(a_i,b_i)$ with $b_i-x_i, x_i-a_i>\Delta_\phi$ for each $i$. (Although for $\phi<\infty$ such a box is not open in the $\phi$-topology, $x$ nevertheless belongs to its interior.) Hence, the ``local property'' expressed in that remark refers to the $\phi$-topology.
    \item The statements on local linearity (Corollary~\ref{cor_locallin}) and d-minimality (Corollary~\ref{corfornasierodmin}) are given only in the $\phi=\infty$ case for simplicity, although there exist uniform versions of these statements as well using the $\phi$-topology (alternatively, working on the quotients $\bar{\Gamma}^n$).
\end{itemize}
\end{remark}

\begin{cor}\label{dim_equals_dim_closure}
Let $X\subseteq\Gamma^n$ be definable.
Then $X$ has the same $\phi$-dimension as its topological closure in the $\phi$-topology.
Hence, $X$ has the same $\phi$-dimension as its topological closure in the topology coming from the order topology on~$\Gamma$. \end{cor}
\begin{proof}
This is \cite[Proposition~2.1]{fornasierohalupczok}, which is available by Remark~\ref{remark_phitop}.
\end{proof}

\subsection{Failure of elimination of imaginaries}\label{subsec_EIfailure} First observe the following lemma:

\begin{lemma}
Let $X \subseteq \Gamma$ be a $\varnothing$-definable set and let $E$ be a $\varnothing$-definable equivalence relation on $X$. Suppose that in the standard model $E$ has uncountably many uncountable equivalence classes. Then there is no $\varnothing$-definable function $f\colon X\to\Gamma_{\infty}^{n}$ such that $f(x) = f(y) \Leftrightarrow xEy$ for $x,y \in X$.
\end{lemma}
\begin{proof}
Suppose towards contradiction that there is such a function $f$. Let 
\[
Z\ =\ \{z \in f(X): \dim_\infty f^{-1}(z) = 1\}.
\]
By our assumption on $E$ and the equivalence \eqref{smallequiv1}$\Leftrightarrow$\eqref{smallequiv11} in the Small Sets Theorem, $\dim_\infty(Z)\geq 1$. But then $\dim_{\infty}(X)\geq 2$ by Fact~\ref{fact_vdD89}\eqref{dim_funct_fact_map_fibers}, contradicting that $X \subseteq \Gamma$. 
\end{proof}

\noindent
This lemma can be used to show that $T_{\log}$ does not have elimination of imaginaries in the $1$-sorted language $\mathcal{L}_{\log}$; see~\cite[Lemma 8.4.7]{tent2012course}. Specifically, we consider the analogue of the so-called \emph{RV sort}, which corresponds to the equivalence relation $\sim_{\psi}$ on $\Gamma^{\neq}$ defined via:
\[
x \ \sim_{\psi} \ y\quad \iff \quad \psi(x) \ < \ \psi(x-y).
\]
\begin{cor}\label{RV_sort_not_EI}
The imaginary $\Gamma^{\neq}/\!\sim_{\psi}$ is not eliminated. \end{cor}

\section{Proof of the Small Sets Theorem}\label{sec_sstproof}

\noindent
In this section we prove the Small Sets Theorem as explained in Subsection~\ref{subsec_roadmap} and indicated in the diagram below. To avoid any confusion, we emphasize that the Small Sets Theorem and the Dimension Theorem are proved simultaneously in three steps:
\begin{itemize}
\item[(1st)] the equivalences between ~\eqref{smallequiv1}--\eqref{smallequiv7} and~\eqref{smallequiv13} of the Small Sets Theorem are established,
\item[(2nd)] the Dimension Theorem is then proved as explained in Subsection~\ref{subsec_dimthmproof}, and
\item[(3rd)] the remaining equivalences of the Small Sets Theorem are then established, \emph{possibly with the help of the Dimension Theorem}.
\end{itemize}

\noindent
\emph{In this section, we adopt the same global conventions as in Section~\ref{sec_dimensions_SST}.}

\subsection{Roadmap of the proof}\label{subsec_roadmap}
We provide the proof in stages, with each stage increasing the set of properties proved to be equivalent to \eqref{smallequiv1}.

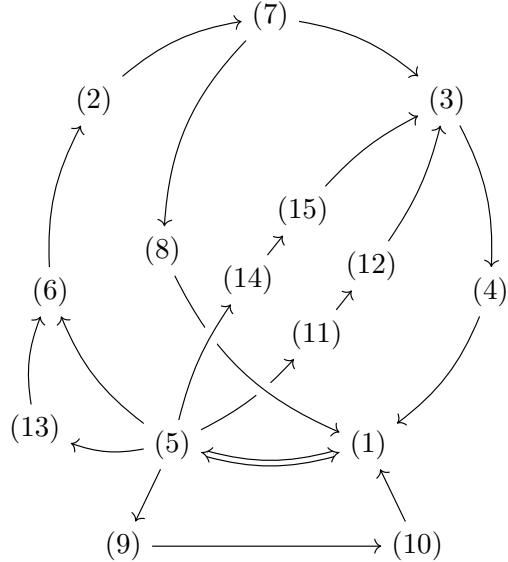
\begin{center}
\begin{tikzpicture}[commutative diagrams/every diagram]
\node (C7) at (0,3)
{\eqref{smallequiv7}};

\node (C3) at (3*.782,3*.623)
{\eqref{smallequiv3}};

\node (C2) at (-3*.782,3*.623)
{\eqref{smallequiv2}};

\node (C4) at (3*.975,-3*.223)
{\eqref{smallequiv4}};

\node (C6) at (-3*.975,-3*.223) 
{\eqref{smallequiv6}};

\node (C1) at (3*.434,-3*.901)
{\eqref{smallequiv1}};

\node (C5) at (-3*.434,-3*.901)
{\eqref{smallequiv5}};

\node (C10) at (3*1.5*.434,-3*1.5*.901)
{\eqref{smallequiv10}};

\node (C9) at (-3*1.5*.434,-3*1.5*.901)
{\eqref{smallequiv9}};

\node (C12) at (-1.2*.434+1.8*.782+.3*1.524,-1.2*.901+1.8*.623-.3*1.216)
{\eqref{smallequiv12}};

\node (C11) at (-1.8*.434+1.2*.782+.3*1.524,-1.8*.901+1.2*.623-.3*1.216)
{\eqref{smallequiv11}};

\node (C15) at (-1.2*.434+1.8*.782-.3*1.524,-1.2*.901+1.8*.623+.3*1.216)
{\eqref{smallequiv15}};

\node (C14) at (-1.8*.434+1.2*.782-.3*1.524,-1.8*.901+1.2*.623+.3*1.216)
{\eqref{smallequiv14}};

\node (C8) at (-1.5*.434-1.5*.975+.678,-3*.223+.541)
{\eqref{smallequiv8}};

\node (C13) at (-1.5*.434-1.5*.975-1.5*.678,-1.5*.901-1.5*.223-1.5*.541)
{\eqref{smallequiv13}};

\path[commutative diagrams/.cd, every arrow, every label]
(C1) edge[commutative diagrams/Leftrightarrow, bend left = 15] 
(C5)

(C15) edge[commutative diagrams/bend left = 10] 
(C3)

(C3) edge[commutative diagrams/bend left = 15] (C4)

(C2) edge[commutative diagrams/bend left = 15] 
(C7)

(C5) edge[commutative diagrams/bend left = 15] 
(C6)

(C6) edge[commutative diagrams/bend left = 15] 
(C2)

(C4) edge[commutative diagrams/bend left = 15] 
(C1)

(C5) edge[commutative diagrams/bend left = 15] 
(C13)

(C13) edge[commutative diagrams/bend left = 15] (C6)

(C5) edge 
(C9)

(C9) edge 
(C10)

(C10) edge 
(C1)

(C7) edge[commutative diagrams/bend left = 15] 
(C3)

(C7) edge[commutative diagrams/bend right=20] (C8)

(C14) edge (C15)

(C5) edge[commutative diagrams/bend right = 10] (C11)

(C8) edge[commutative diagrams/crossing over, bend right=20] 
(C1)

(C11) edge (C12)

(C12) edge[commutative diagrams/bend right = 10] (C3)

(C5) edge[commutative diagrams/crossing over, bend left = 10] 
(C14);
\end{tikzpicture}
\end{center}

\emph{Stage (I): Incorporating~\eqref{smallequiv5}.} We start with \eqref{smallequiv1}$\Leftrightarrow$\eqref{smallequiv5}, which is just Lemma~\ref{lem_dphi0equivdphiset}.

\emph{Stage (II): Incorporating \eqref{smallequiv2}, \eqref{smallequiv3}, \eqref{smallequiv4}, \eqref{smallequiv6},
\eqref{smallequiv7},
and \eqref{smallequiv13} (Subsections~\ref{finding_open_boxes_subsection},~\ref{Psi_functions_CB_subsection}, and~\ref{subsec_equilat}).}
\begin{itemize}
\item \eqref{smallequiv4}$\Rightarrow$\eqref{smallequiv1} Consider the contrapositive: suppose $\operatorname{dim}_{\phi}X=1$. Then there is some $\alpha$ in an elementary extension of $(\Gamma,\psi)$ such that $\alpha\in X^*$ and $\operatorname{rk}_{\phi}(\alpha|\Gamma)=1$. Since $X\in\operatorname{tp}(\alpha|\Gamma)$, it follows from Proposition~\ref{prop_simpleextinterior} that $X$ contains an interval $(a,b)$ such that $b-a>\Delta_{\phi}$.
\item \eqref{smallequiv5}$\Rightarrow$\eqref{smallequiv6} (case $\phi=\infty$) This is Corollary~\ref{cor_dsetCBrankemptyint}, which follows from explicitly computing the derived set of the image of an affine map $\Psi^n\to\Gamma$ in Proposition~\ref{prop_derivePsifun}. To put this in more topological terms: Lemma~\ref{lem_primdsetCBrank} shows that the image of an affine map $\Psi^n\to\Gamma$ is d-finite, so \eqref{smallequiv6} follows from the fact that the d-finite sets form an ideal (see Lemma~\ref{CB_small_ideal_properties}).
\item \eqref{smallequiv5}$\Rightarrow$\eqref{smallequiv13}$\Rightarrow$\eqref{smallequiv6} (case $\phi<\infty$) The step \eqref{smallequiv5}$\Rightarrow$\eqref{smallequiv13} is Corollary~\ref{cor_dphibarcloseddisc}. The direction \eqref{smallequiv13}$\Rightarrow$\eqref{smallequiv6} is trivial.
\item \eqref{smallequiv6}$\Rightarrow$\eqref{smallequiv2} This direction is always true for arbitrary subsets of arbitrary topological spaces, i.e., the ideal of d-finite sets is always contained in the ideal generated by the discrete sets (see Lemma~\ref{CBsmall_implies_finite_union_discrete}). Note that in complete generality, the converse can fail (see Examples~\ref{counterexample_discrete_not_CB_finite} and~\ref{second_counterexample_discrete_not_CB_finite}).
\item  \eqref{smallequiv2}$\Rightarrow$\eqref{smallequiv7} This direction is true for arbitrary subsets of arbitrary topological spaces that are $T_1$ and have no isolated points, i.e., the ideal generated by discrete sets is always contained in the ideal of nowhere dense sets under these assumptions; see Lemma~\ref{discrete_implies_nowhere_dense_under_mild_assumptions}. Note that the converse will fail for definable subsets of $\Gamma^n$ for $n>1$.
\item  \eqref{smallequiv7}$\Rightarrow$\eqref{smallequiv3} This direction is always true; see Lemma~\ref{nowhere_dense_implies_empty_interior}.
\item \eqref{smallequiv3}$\Rightarrow$\eqref{smallequiv4} This step is trivial. 
\end{itemize}

Note: At this point in the proof of the Small Sets Theorem we can prove the Dimension Theorem (Corollaries~\ref{dimension_theorem_uniqueness} and~\ref{dimension_theorem_existence}), and thus we may freely use the fact that each $\dim_{\phi}$ is a dimension function in the rest of the proof of the Small Sets Theorem.

\emph{Stage (III): Incorporating \eqref{smallequiv8} (definably meager).}
We now say what it means for an arbitrary subset of $\bar{\Gamma}^n$ to be \emph{definably meager}; recall that we are working in a $1$-sorted setting, so we give a definition purely in terms of definable subsets of $\Gamma^n$.

\begin{definition}
Suppose $Z\subseteq\bar{\Gamma}^n$. We say that $Z$ is \textbf{definably meager} in $\bar{\Gamma}^n$ if there exist a definable directed set $(I;\leq)$ and a definable family $Y\subseteq I\times\Gamma^n$ such that for every $a,b\in I$:
\begin{enumerate}
\item $\overline{Y_a}$ is nowhere dense in $\bar{\Gamma}^n$,
\item if $a\leq b$, then $\overline{Y_a}\subseteq\overline{Y_b}$, and
\item $Z\subseteq\bigcup_a\overline{Y_a}$.
\end{enumerate}
In this case, we say that $(I;\leq)$ and $Y$ \textbf{witness} that $Z$ is definably meager.
\end{definition}

\begin{itemize}
\item \eqref{smallequiv7}$\Rightarrow$\eqref{smallequiv8} This direction is trivial. Indeed, if $\bar{X}$ is nowhere dense, then the singleton directed set $(\{\ast\};\leq)$ and the definable family $Y\coloneqq\{\ast\}\times X$ witnesses that $\bar{X}$ is definably meager.
\item \eqref{smallequiv8}$\Rightarrow$\eqref{smallequiv1} Suppose $\bar{X}$ is definably meager, witnessed by $(I;\leq)$ and $Y\subseteq I\times\Gamma$. In particular, each $Y_a$ is $\phi$-small by \eqref{smallequiv7}$\Rightarrow$\eqref{smallequiv1}. By replacing $Y_a$ with $\bigcup_{b\leq a}Y_b$ (which does not change the set $\overline{Y_a}$, so $Y_a$ is still $\phi$-small), we may assume that $Y_a\subseteq Y_b$ for $a\leq b\in I$.
It now follows from Corollary~\ref{dim_directed_union_fact} that $W\coloneqq\bigcup_{a\in I}Y_a$ also satisfies $\dim_{\phi}W\leq 0$. Thus $\bar{W}$ is nowhere dense by \eqref{smallequiv1}$\Rightarrow$\eqref{smallequiv7}, and so $\bar{X}\subseteq \bar{W}$ is also nowhere dense.
\end{itemize}

\emph{Stage (IV): Incorporating \eqref{smallequiv9} and \eqref{smallequiv10} (internality and coanalysability; Subsection~\ref{interality_coanalysability_subsection}).}
\begin{itemize}
\item \eqref{smallequiv5}$\Rightarrow$\eqref{smallequiv9} This is clear since \eqref{smallequiv5} is expressing a precise form of \emph{internality}; see Lemma~\ref{lem_dsetPsiint}.
\item \eqref{smallequiv9}$\Rightarrow$\eqref{smallequiv10} This is because internality always implies coanalysability; see Remark~\ref{rem_coanal}.
\item \eqref{smallequiv10}$\Rightarrow$\eqref{smallequiv1} If $X$ is $\Psi\Delta_{\phi}$-coanalysable, then since $\dim_{\phi}\Psi\Delta_{\phi}=0$, it follows from Lemma~\ref{lem_coanaldphi0} that $\dim_{\phi}X\leq 0$. Note that this step uses the fact that $\dim_{\phi}$ is a dimension function, which we established after Stage (II) above.
\end{itemize}

\emph{Stage (V): Incorporating \eqref{smallequiv11} and \eqref{smallequiv12} (cardinality in the standard model).}
\begin{itemize}
\item \eqref{smallequiv5}$\Rightarrow$\eqref{smallequiv11} If $X$ satisfies \eqref{smallequiv5}, then $X\subseteq Y+\Delta_{\phi}$ for some $\infty$-small definable set $Y$. 
Since $Y$ is countable in the standard model and $\bar{X}\subseteq \bar{Y}$, it follows that $|\bar{X}|\leq\aleph_0$.
\item \eqref{smallequiv11}$\Rightarrow$\eqref{smallequiv12} This is trivial.
\item \eqref{smallequiv12}$\Rightarrow$\eqref{smallequiv3} This follows easily from the observation that in the linear order $\Gamma_{\log}/\Delta_{\phi}$, each interval has size continuum.
\end{itemize}

\emph{Stage (VI): Incorporating \eqref{smallequiv14} and \eqref{smallequiv15} (dp-rank; Subsection~\ref{dp_rank_subsection}).}
\begin{itemize}
\item \eqref{smallequiv5}$\Rightarrow$\eqref{smallequiv14} This is by Fact~\ref{dpfacts}, which uses that $\operatorname{dp}(\Psi)=1$.
\item \eqref{smallequiv14}$\Rightarrow$\eqref{smallequiv15} This is trivial.
\item \eqref{smallequiv15}$\Rightarrow$\eqref{smallequiv3} This is Lemma~\ref{dpint}.
\end{itemize}

\subsection{Finding open boxes}\label{finding_open_boxes_subsection}

Below, we show that any definable subset of $\Gamma$ of $\phi$-dimension 1 contains an open interval of width $>\Delta_\phi$, yielding \eqref{smallequiv4}$\Rightarrow$\eqref{smallequiv1} of the Small Sets Theorem. This is the unary case of the more general Proposition~\ref{prop_simpleextinterior}, which we need to establish the analogous direction \eqref{ksmallequiv3}$\Rightarrow$\eqref{ksmallequiv1} of the Kinda Small Sets Theorem. It is not clear to us whether this $n$-ary version can be deduced from the unary version.

\begin{prop}\label{prop_simpleextinterior}
Let $\alpha = (\alpha_1,\ldots,\alpha_n)$ be a tuple in an elementary extension of $(\Gamma,\psi)$ with $\rk_\phi(\alpha|\Gamma) = n$. Then  every set $X$ in $\operatorname{tp}(\alpha|\Gamma)$ contains an open box $\prod_{i=1}^n(a_i,b_i)$ where $b_i-a_i>\Delta_{\phi}$ for each $i$.
\end{prop}
\begin{proof}
We begin by replacing $\Gamma$ with $\operatorname{cl}_\phi(\Gamma)\subseteq \Gamma\langle \alpha\rangle$; then $\Gamma$ is still a model of $T_{\log}$ by Lemma~\ref{lem_elemsub} and we still have $\rk_\phi(\alpha|\Gamma)=n$, but we have arranged that $\Psi_{\Gamma\langle \alpha\rangle} = \Psi$, so 
\[
\Gamma \langle \alpha\rangle\ \cong\ \Gamma \oplus \bigoplus_{i=1}^n\Q\alpha_i.
\]
Take $\delta$ in a further elementary extension $(\Gamma^*,\psi^*)$ of $\Gamma\langle \alpha\rangle$ with 
\[
\Delta_\phi^*\ <\ \delta\ <\ \Gamma\langle \alpha\rangle^{>\Delta_\phi^*}. 
\]
Note that if $\phi<\infty$, then $\psi^*(\delta) = \phi$ and that if $\phi = \infty$, then $\psi^*(\delta)> \Psi$. By compactness, it is enough to show that if $X$ is in $\operatorname{tp}(\alpha|\Gamma)$, then $X^*$ contains the box $\prod_{i=1}^n(\alpha_i-\delta,\alpha_i+\delta)$. That is, we must show that for any $\varepsilon = (\varepsilon_1,\ldots,\varepsilon_n) \in (\Gamma^*)^n$ with $|\varepsilon_i|<\delta$ for each $i$, we have $\operatorname{tp}(\alpha|\Gamma) = \operatorname{tp}(\alpha+\varepsilon|\Gamma)$. Fix such an $\varepsilon$; then quantifier elimination allows us to further reduce the problem to finding an $\calL_{\log}$-isomorphism $\iota\colon\Gamma\langle \alpha\rangle \to \Gamma\langle \alpha+\varepsilon\rangle$ over $\Gamma$ with $\iota(\alpha)=\alpha+\varepsilon$. 

Note that for any $\beta \in \Gamma\langle \alpha\rangle \setminus \Gamma$, we have $\psi^*(\beta)\leq \phi$. Thus,  $|\beta| \in \Gamma\langle \alpha\rangle^{>\Delta_\phi^*}$, so $[\varepsilon_i]\leq [\delta]< [\beta]$. It follows that for $\gamma \in \Gamma$ and $q \in \Q^n$, the element $\gamma+ q\cdot \alpha$ is positive if and only if $\gamma+ q\cdot (\alpha+\varepsilon)$ is positive, so we have an ordered abelian group isomorphism 
\[
\iota \colon\Gamma\oplus\bigoplus_{i=1}^n\Q\alpha_i\to\Gamma\oplus\bigoplus_{i=1}^n\Q(\alpha_i+\varepsilon_i)
\]
over $\Gamma$ with $\iota(\alpha)=\alpha+\varepsilon$.

It remains to prove that $\iota$ commutes with $\psi$, $s$, and $p$. Consider an arbitrary element $\beta \in \Gamma\langle \alpha\rangle$; we can assume that $\beta \not\in \Gamma$. Since $\Psi_{\Gamma\langle\alpha\rangle} = \Psi$, it is enough to show that $\psi^*(\iota(\beta)) = \psi^*(\beta)$, and likewise for $s$ and $p$. 
Take $q \in \Q^n$ distinct from the zero tuple with $\iota(\beta) = \beta + q\cdot \varepsilon$.  Since $[\varepsilon_i]<[\beta]$ for each $i$, we have $[\beta + q\cdot\varepsilon] = [\beta]$; in particular, $\psi^*(\beta+q\cdot\varepsilon) = \psi^*(\beta)$. The Integral Identity (Fact~\ref{fact_identities}) gives $\textstyle{\int}\beta = \beta - s\beta \in \Gamma\langle \alpha\rangle\setminus \Gamma$, so we have
\[
(\textstyle{\int} \beta+q\cdot\varepsilon)'\ =\ \textstyle{\int} \beta+q\cdot\varepsilon+\psi^*(\textstyle{\int} \beta+q\cdot\varepsilon)\ =\ \textstyle{\int} \beta+q\cdot\varepsilon+\psi^*(\textstyle{\int} \beta)\ =\ \beta+q\cdot\varepsilon.
\]
It follows that $\textstyle{\int} \beta+q\cdot\varepsilon = \textstyle{\int}( \beta+q\cdot\varepsilon)$, and we again use the Integral Identity to conclude that $s\beta = \beta - \textstyle{\int}\beta = \beta +q\cdot\varepsilon -  \textstyle{\int}( \beta+q\cdot\varepsilon) = s(\beta+q\cdot\varepsilon)$. Finally for $p$, we note that since $\beta \not\in \Gamma  \supseteq \Psi = \psi^*((\Gamma\oplus\Q\alpha)^{\neq})$ and the same for $\beta+q\cdot\varepsilon$, we have $p(\beta)=\infty=p(\beta+q\cdot\varepsilon)$.
\end{proof}

\subsection{Derived sets and \texorpdfstring{$\Psi$}{Ψ}-functions} \label{Psi_functions_CB_subsection}
We establish here the direction \eqref{smallequiv5}$\Rightarrow$\eqref{smallequiv6} of the Small Sets Theorem in the case $\phi=\infty$.

\medskip\noindent
For $\alpha,\beta\in\Gamma$ with $\beta\neq 0$, set $\alpha \prec_{\psi} \beta$ if $\psi(\alpha)>\psi(\beta)$, so if $\alpha\prec_{\psi}\beta$, then $[\alpha]<[\beta]$. 
\begin{lemma}\label{lem_arbsmalldiff}
For every $\varepsilon \in \Gamma^>$, there exist $\delta_0,\delta_1\in \Psi^{\geq\psi(\varepsilon)}$ such that $0<\delta_0-\delta_1 \prec_{\psi} \varepsilon$.
\end{lemma}
\begin{proof}
Set $\delta_0\coloneqq s\psi(\varepsilon)$ and $\delta_1\coloneqq\psi(\varepsilon)$.
Then $\delta_0-\delta_1 = s\psi(\varepsilon)-\psi(\varepsilon) = -\textstyle{\int}\psi(\varepsilon)$
by the Integral Identity (Fact~\ref{fact_identities}).
It remains to note that $0<-\textstyle{\int}\psi(\varepsilon) \prec_{\psi} \varepsilon$ by \cite[Lemma~9.2.18(iii)]{ADAMTT}.
\end{proof}

\begin{def*}
A \textbf{$\Psi$-function} is $F\colon \Psi^I\to\Gamma$ defined by $F\alpha=\sum_{i\in I}q_i\alpha_i + \beta$ for $\alpha=(\alpha_i)_{i\in I} \in \Psi^I$, where $I\subseteq \N$ is finite, $(q_i)_{i\in I} \in (\Q^{\times})^I$, and $\beta\in\Gamma$. 
\end{def*}
\noindent
We allow $I=\0$, in which case $F$ has constant value $\beta$.
Note that each $\infty$-small subset of $\Gamma$ is contained in finite unions of images of $\Psi$-functions.
Studying the limit points of such images will establish the part of the Small Sets Theorem claimed above.

\medskip\noindent
In the remainder of this subsection, $F\colon \Psi^I \to \Gamma$ is a $\Psi$-function given by $F\alpha=\sum_{i \in I}q_i\alpha_i + \beta$ as in the definition above.
We set $\|F\| \coloneqq \sum_{i\in I}q_i$. For $J \subseteq I$, we define the $\Psi$-function $F_J \colon \Psi^J \to \Gamma$ by
\[
F_J\alpha \ \coloneqq \ \sum_{j\in J}q_j\alpha_j+\beta,
\]
for $\alpha=(\alpha_j)_{j\in J}\in\Psi^J$,
in which case $F\alpha = F_J(\alpha_j)_{j\in J}+F_{I\setminus J}(\alpha_i)_{i\in I\setminus J} - \beta$ for all $\alpha \in \Psi^I$.
When convenient, we implicitly regard $F$ as a function $\Psi^J \to \Gamma$, for $J \supseteq I$, by setting $q_j=0$ for $j\in J \setminus I$.
In particular, for $J \subseteq I$, we write $F_J\alpha$ instead of $F_J(\alpha_j)_{j\in J}$, for $\alpha\in\Psi^I$.

\medskip\noindent
Recall that $\Psi$ is $\Q$-linearly independent \cite[Lemma~6.8]{GehretACTlog}, which yields the following uniqueness property. 

\begin{lemma}\label{key_uniqueness_lemma}
If $G\colon \Psi^J \to \Gamma$ is a $\Psi$-function with $G\alpha=\sum_{j\in J}\tilde{q}_j\alpha_j + \tilde{\beta}$ such that 
$F\alpha=G\alpha$
for all $\alpha \in \Psi^{I\cup J}$, then $I=J$, $q_i=\tilde{q}_i$ for all $i \in I$, and $\beta=\tilde{\beta}$. 
\end{lemma}

\begin{lemma}\label{lem_arbsmallPsifun}
Suppose that $\|F\|=0$ and $I\neq\0$.
Then for every $\varepsilon\in\Gamma^>$ there exists $\alpha\in\Psi^I$ such that
$0 < |F\alpha-\beta| \prec_{\psi} \varepsilon$.
\end{lemma}
\begin{proof}
Pick $i_0 \in I$ and set $I_0\coloneqq I\setminus \{i_0\}$.
Let $\varepsilon>0$, so Lemma~\ref{lem_arbsmalldiff} gives $\delta_0,\delta_1\in\Psi$ such that $0<\delta_0-\delta_1\prec_{\psi}\varepsilon$. Then $\alpha\in \Psi^{I}$ defined by
$\alpha_{i_0} \coloneqq \delta_0$ and $\alpha_i \coloneqq\delta_1$ for $i \in I_0$
satisfies
\[
F\alpha-\beta \ = \ q_{i_0}\alpha_{i_0}+\textstyle\sum_{i\in I_0}q_i\alpha_i \ = \ q_{i_0}(\delta_0+q_{i_0}^{-1}\sum_{i\in I_0}q_i \delta_1) \ = \ q_{i_0}(\delta_0-\delta_1),
\]
and thus $0 < |F\alpha-\beta| \prec_{\psi} \varepsilon$.
\end{proof}

\noindent
In fact, the assumptions on $F$ in the lemma above characterize when $\beta$ is a limit point of $\image(F)$ in $\Gamma$.
\begin{cor}
We have $\beta \in \image(F)' \Leftrightarrow \|F\|=0\ \text{and}\ I\neq\0$.
\end{cor}
\begin{proof}
If $\|F\|=0$ and $I\neq\0$, then
Lemma~\ref{lem_arbsmallPsifun} shows that $\beta\in\image(F)'$.
Conversely, if $I=\0$,
then $\image(F)=\{\beta\}$ has no limit points, and if $\|F\|\neq 0$, then \cite[Lemma~6.4]{GehretACTlog} gives $\psi(F\alpha-\beta)=s0$ for all $\alpha\in\Psi^I$.
\end{proof}

\begin{prop}\label{prop_derivePsifun}
We have
\[
\image(F)' \ = \ \bigcup_{\substack{\0\neq J\subseteq I \\ \|F_J\|=0}} \image(F_{I\setminus J}).
\]
\end{prop}
\begin{proof}
First, suppose that $\0 \neq J \subseteq I$ and $\|F_J\|=0$.
Let $(\alpha_i)_{i\in I\setminus J}\in\Psi^{I\setminus J}$ and $\varepsilon\in\Gamma^>$ be arbitrary.
Lemma~\ref{lem_arbsmallPsifun} gives $(\alpha_j)_{j \in J}\in\Psi^J$ so that $0<|F_J\alpha-\beta|\prec_{\psi}\varepsilon$,
where $\alpha\coloneqq(\alpha_i)_{i \in I}$.
Note that
\[
F_J\alpha \ = \ F_J\alpha+F_{I\setminus J}\alpha-F_{I\setminus J}\alpha \ = \  F\alpha+\beta-F_{I\setminus J}\alpha,
\]
so $0 < |F_{I\setminus J}\alpha-F\alpha| \prec_{\psi} \varepsilon$.
Hence $F_{I\setminus J}\alpha \in \image(F)'$.

For the reverse inclusion, take an elementary extension $(\Gamma^*,\psi^*)$ of $(\Gamma,\psi)$ containing an element $\delta\in(\Psi^*)^{>\Psi}$. By QE and UA, we arrange that $\Gamma^*=\Gamma\langle\delta\rangle$.
The proof of \cite[Lemma~4.11]{GehretACTlog} shows
\[
\Gamma\langle\delta\rangle \ = \ \Gamma\oplus\bigoplus_{k\in\Z}\Q s^k\delta
\]
as an internal direct sum of $\Q$-linear subspaces, where $s^k\delta\coloneqq p^{-k}\delta$ for $k<0$ in $\Z$.
Fix $\gamma\in\Gamma$ and assume that for every $\varepsilon \in \Gamma^>$, there exists $\alpha \in \Psi^I$ such that $0<|F\alpha-\gamma|<\varepsilon$.
Set $\varepsilon(\delta)\coloneqq -\textstyle{\int}\delta$, which satisfies $0<\varepsilon(\delta)<\Gamma^>$ since $\textstyle{\int}\Psi$ is cofinal in $\Gamma^<$ by \cite[Lemma~9.2.15]{ADAMTT}.
By elementarity, this yields $\alpha^*\in (\Psi^*)^I$ such that
\[
0 \ < \ |F\alpha^*-\gamma| \ < \ \varepsilon(\delta).
\]
It follows that $F\alpha^*\not\in\Gamma$, so $J\coloneqq\{i\in I: \alpha^*_i>\Psi\} \neq \0$.
Note that $I \setminus J = \{ i \in I : \alpha^*_i \in \Psi \}$. This gives for each $j\in J$ a $k_j\in\Z$ such that $\alpha_j^*=s^{k_j}\delta$.
Setting $\tilde{\beta}\coloneqq \gamma-F_{I\setminus J}\alpha^* \in \Gamma$,
we have \[
0 \ < \ \big|\sum_{j\in J}q_js^{k_j}\delta-\tilde{\beta}\big| \ < \ \varepsilon(\delta).
\]
By combining some of the coefficients $q_j \in \Q^{\times}$ and shrinking $I$ if necessary, we arrange that the integers $k_j$ are distinct.
Note that we still have $J\neq\0$.
Then
\[
\Psi^* \ \ni \ \psi\big(\sum_{j\in J}q_js^{k_j}\delta-\tilde{\beta}\big) \ \geq \ \psi(\varepsilon(\delta))\ =\ s\delta\ > \ \Psi.
\]
Now, \cite[Theorem 6.6]{GehretACTlog} completely characterizes the values of $\psi(\sum_{j\in J}q_js^{k_j}\delta-\tilde{\beta})$.
Since $\psi(\tilde{\beta})\in\Psi_{\infty}$, this characterization forces either $s^{k+1}\delta<\psi(\tilde{\beta})$ or $s^{k+1}\delta<s(q^{-1}\tilde{\beta})$, where $k\coloneqq\min_{j\in J}\{k_j\}$ and $q\coloneqq\sum_{j\in J}q_j$.
The latter is impossible, since $s(q^{-1}\tilde{\beta})\in\Psi<s^{k+1}\delta$.
Hence $s^{k+1}\delta<\psi(\tilde{\beta})$, in which case $\Psi<s^{k+1}\delta$ gives $\psi(\tilde{\beta})=\infty$.
That is, $\tilde{\beta}=0$, as desired.
That $q=0$ also follows from the case distinctions in \cite[Theorem 6.6]{GehretACTlog}.
\end{proof}

\begin{lemma}\label{lem_primdsetCBrank}
If $|I|=n \geq 1$, then $\image(F)^{(n)}=\0$.
\end{lemma}
\begin{proof}
Note that $\image(F)'=\{\beta\}'=\0$ when $I=\0$.
By Proposition~\ref{prop_derivePsifun}, we have
\[
\image(F)' \ = \ \bigcup_{\substack{\0\neq J\subseteq I \\ \|F_J\|=0}} \image(F_{I\setminus J}).
\]
By induction on $n$, for each nonempty $J\subseteq I$ we have $\image(F_{I\setminus J})^{(n-1)}=\0$. Then
\[
\image(F)^{(n)} \ = \ \bigcup_{\substack{\0\neq J\subseteq I \\ \|F_J\|=0}} \image(F_{I\setminus J})^{(n-1)} \ = \ \0,
\]
as desired (see Lemma~\ref{CB_der_properties}).
\end{proof}

\noindent
Note that Lemma~\ref{lem_primdsetCBrank} only provides an upper bound.
For example, the set $X=\Psi+\dots+\Psi$ of $n$-fold sums satisfies $X'=\0$.

\medskip\noindent
The following gives  \eqref{smallequiv5}$\Rightarrow$\eqref{smallequiv6} in the Small Sets Theorem in the case $\phi=\infty$:

\begin{cor}\label{cor_dsetCBrankemptyint}
If $X$ is $\infty$-small, then $X^{(n)}=\0$ for some $n \in \N$.
\end{cor}
\begin{proof}
This is immediate from Lemma~\ref{lem_primdsetCBrank} and the fact that d-finite sets form an ideal (see Lemma~\ref{CB_small_ideal_properties}).
\end{proof}

\noindent
We have now completed enough of the Small Sets Theorem (i.e., Stage (II) in case $\phi=\infty$) to get that $\dim_{\infty}$ is a dimension function; see Corollary~\ref{dimension_theorem_existence}. 
This yields some additional uniformity in families, for which we define a \textbf{parametrized $\Psi$-function} to be $F \colon \Psi^n \times \Gamma^m \to \Gamma$ defined by
\[
F(\alpha,\beta)\ =\ \sum_{i=0}^{n-1} q_i\alpha_i + \sum_{j=0}^{m-1} r_j\beta_j,
\]
for $\alpha=(\alpha_0,\dots,\alpha_{n-1})\in\Psi^n$ and $\beta=(\beta_0,\dots,\beta_{m-1})\in\Gamma^m$, where $q_i\in \Q^{\times}$ for $i=0,\dots,n-1$ and $r_j\in \Q$ for $j=0,\dots,m-1$.
\begin{cor}\label{cor_dimfibre}
If $X \subseteq \Gamma^{m+1}$ is $A$-definable, then there exist $a \in A^{\ell}$ and parametrized $\Psi$-functions $F_i \colon \Psi^{n_i}\times\Gamma^\ell\times\Gamma^m \to \Gamma$ for $i=0,\dots,k-1$ such that for all $\delta\in\Gamma^{m}$ with $\dim_{\infty} X_{\delta}=0$, we have
\[
X_{\delta}\ \subseteq\ \bigcup_{i=0}^{k-1} \image(F_i(\cdot,a,\delta)).
\]
\end{cor}
\begin{proof}
By Lemma~\ref{lem_existentialmat}, $\cl_{\infty}$ is an existential matroid on a monster model $(\Gamma^*,\psi^*) \succcurlyeq (\Gamma,\psi)$, and $\dim_{\infty}$ is a dimension function.
Hence, $\dim_{\infty}Y=\sup\{ \rk_{\infty}(y|A) : y \in Y^* \}$ for $A$-definable $Y\subseteq \Gamma^n$ by \cite[Theorem~4.3 and Remark~3.45]{fornasiero-dimmatroiddense} (see also \cite[Definition~3.29]{fornasiero-dimmatroiddense}).
Compactness arguments as in \cite[Lemmas~2.3 and 2.7]{angel2016bounded} then yield the desired result, a stronger form of that latter lemma.
\end{proof}

\begin{cor}
If $X \subseteq \Gamma^{m+1}$ is definable, then there is $N \in \N$ such that for all $\delta\in\Gamma^{m}$ with $\dim_{\infty} X_{\delta}=0$, we have $X_{\delta}^{(N)}=\0$.
\end{cor}
\begin{proof}
Note that $\image(F(\cdot,\beta))^{(n)}=\0$ for all $\beta\in\Gamma$, where $F\colon \Psi^n \times \Gamma \to \Gamma$ is a parametrized $\Psi$-function, so this result follows from the previous corollary.
\end{proof}

\subsection{Equilateral sets and quotients}\label{subsec_equilat}
In this subsection we establish \eqref{smallequiv5}$\Rightarrow$\eqref{smallequiv13} of the Small Sets Theorem, which follows from purely valuation-theoretic arguments.
For that, let $(G, S; v)$ be an ordered abelian group with a convex (surjective) valuation $v \colon G \to S_{\infty}$, where $S$ is an ordered set and we extend $S$ to $S_{\infty} \coloneqq S \cup \{\infty\}$ in the usual way.
Note that if $S$ has no greatest element, then the valuation topology on $G$ agrees with the order topology on $G$.
On the other hand, if $S$ has a maximum, then the valuation topology on $G$ is discrete.
See for instance \cite[Chapter~2]{ADAMTT} for basic definitions and facts.

\medskip\noindent
Let $s \in S$ and $X \subseteq G$.
We say $X$ is \textbf{$s$-equilateral} if $v(a-b)=s$ for all distinct $a,b \in X$.
We say $X$ is \textbf{anti-equilateral} if for every $s$, $X$ contains no infinite $s$-equilateral subset.
To see how this property passes to quotients, let $\Delta_{s}\coloneqq\{ a \in G : va>s \}$, which is a convex subgroup of $G$.
Then $v$ induces a convex valuation $v\colon G/\Delta_{s} \to S^{\leq s}_{\infty}$ on the quotient ordered abelian group $G/\Delta_{s}$, defined by $v(a+\Delta_{s})=va$ for $a \in G \setminus \Delta_{s}$.
Let $\bar{X}$ denote the image of $X$ under the quotient map $G \to G/\Delta_{s}$.

\begin{lemma}\label{lem_equilatquot}
If $X \subseteq G$ is anti-equilateral, then $\bar{X} \subseteq G/\Delta_{s}$ is anti-equilateral.
\end{lemma}
\begin{proof}
Let $X \subseteq G$ and suppose that $\bar{X} \subseteq G/\Delta_{s}$ is an infinite $s'$-equilateral set, where $s' \in S^{\leq s}$.
Take $Y \subseteq X$ so that the quotient map $G \to G/\Delta_{s}$ restricts to a bijection $Y \to \bar{X}$.
Then $Y$ is infinite and it is easy to check that $Y$ is $s'$-equilateral, for if $y_1,y_2 \in Y$ with $y_1\neq y_2$, then $y_1-y_2 \notin \Delta_{s}$ and so
\[
v(y_1-y_2)\ =\ v\big((y_1-y_2)+\Delta_{s}\big)\ =\ v\big((y_1+\Delta_{s})-(y_2+\Delta_{s})\big)\ =\ s'. \qedhere
\]
\end{proof}

\begin{lemma}\label{lem_equilatcloseddisc}
Suppose that $S$ has a maximum $s=\max S$ and $X\subseteq G$ is anti-equilateral.
Then $X$ is closed and discrete in the order topology on~$G$.
\end{lemma}
\begin{proof}
We can assume that $G$ is not discrete.
To show that $X$ is discrete, let $x \in X$.
Then we have an interval $(a,b)$ in $G$ with $x \in (a,b)$ and $v(b-a)=s$.
It follows that $(a,b)$ is an infinite $s$-equilateral set, so $(a,b)\cap X$ is finite.
By shrinking $(a,b)$ further, we arrange that $(a,b)\cap X = \{x\}$, as desired.
Similarly, if $x \in \cl(X)$, then $x \in X$.
\end{proof}

\noindent
The previous lemma only uses that $X$ contains no infinite $s$-equilateral subset.
Now we apply these lemmas to the setting $(\Gamma,\psi)\models T_{\log}$, construed as a structure $(\Gamma, \Psi; \psi)$ in the notation above.

\begin{lemma}
Every $\infty$-small $X \subseteq \Gamma$ is anti-equilateral. \end{lemma}
\begin{proof}
Let $F\colon \Psi^n \to \Gamma$ be a $\Psi$-function and fix $\phi\in\Psi$.
By the Pigeonhole Principle, it suffices to prove that $\image(F)$ contains no infinite $\phi$-equilateral subset. Let $Y \subseteq \Psi^n$ be infinite and, for $k=0,\dots,n-1$, let $\pi_k \colon Y \to \Psi$ be projection onto coordinate $k$.
If $\pi_0^{-1}(p\phi)$ is finite, replacing $Y$ by $Y \setminus \pi_0^{-1}(p\phi)$ arranges that $\pi_0^{-1}(p\phi)=\0$ while keeping $Y$ infinite.
If $\pi_0^{-1}(p\phi)$ is infinite, replacing $Y$ by $\pi_0^{-1}(p\phi)$ arranges that $\pi_0^{-1}(p\phi)=Y$ while keeping $Y$ infinite.
Repeating this procedure ensures that for each $k=0,\dots,n-1$, either $\pi_k^{-1}(p\phi)=\0$ or $\pi_k^{-1}(p\phi)=Y$.
Let $I=\{k\in\{0,\dots,n-1\} : \pi_k^{-1}(p\phi)=\0 \}$, so $\psi(F(\alpha)-F(\beta)) = \psi(F_I(\alpha)-F_I(\beta))$ for all $\alpha,\beta\in Y$.
For fixed $\alpha,\beta \in Y$, by combining the coefficients of the equal $\alpha_i$ and $\beta_j$, $i,j\in I$, we see that $F_I(\alpha)-F_I(\beta)$ is equal to a $\Q$-linear combination of distinct such $\alpha_i$, $\beta_j$ with coefficients summing to $0$.
Then by \cite[Lemma~6.4]{GehretACTlog}, $\psi(F_I(\alpha)-F_I(\beta))\in \{ s\alpha_i, s\beta_i : i \in I\}\cup\{\infty\}$.
In particular, $\psi(F(\alpha)-F(\beta))\neq\phi$ for all $\alpha,\beta \in Y$.
Now, if $Y\subseteq\Psi^n$ is such that $F(Y) \subseteq \image(F)$ is $\phi$-equilateral, then the above shows that $F(Y)$ is finite.
\end{proof}

\noindent
Note that the quotient map $\pi_{\phi}\colon\Gamma \to \Gamma/\Delta_{\phi}$ is injective on $\Psi^{<s\phi}$ by the Successor Identity (Fact~\ref{fact_identities}).
In particular, the value set $\Psi^{<s\phi}$ of $\Gamma/\Delta_{\phi}$ is order isomorphic to $\bar{\Psi}$ via $\pi_\phi$, and we consider $\Gamma/\Delta_{\phi}$ as a convexly valued ordered abelian group with $\psi \colon (\Gamma/\Delta_{\phi})^{\neq} \to \bar{\Psi}$. 
Moreover, $(\Gamma/\Delta_{\phi}, \psi)$ is an asymptotic couple by \cite[Lemma~9.2.24]{ADAMTT}, but that is not needed here.
\begin{cor}\label{cor_dphibarcloseddisc}
If $X\subseteq\Gamma$ is $\phi$-small and $\phi\in\Psi$, then $\bar{X}\subseteq\Gamma/\Delta_{\phi}$ is closed and discrete. \end{cor}
\begin{proof}
Let $Y\subseteq \Gamma$ be $\infty$-small, so $Y$ is anti-equilateral. Note that $\psi((\Gamma/\Delta_{\phi})^{\neq})$ has a maximum $\bar{\phi}$.
Then by Lemmas~\ref{lem_equilatquot} and~\ref{lem_equilatcloseddisc}, $\bar{Y}$ is closed and discrete in $\Gamma/\Delta_{\phi}$.
It remains to note that for every $\phi$-small $X \subseteq\Gamma$, there is an $\infty$-small $Y\subseteq\Gamma$ with $\bar{Y}=\bar{X}$.
\end{proof}

\subsection{Internality and coanalysability}\label{interality_coanalysability_subsection}
We now investigate how previous properties are connected to the model-theoretic notion of \emph{internality} to a definable set, namely internality to the distinguished family of sets~$\Psi\cup\Delta_{\phi}$.

\begin{def*}
Let $X \subseteq \Gamma^{n}_{\infty}$ be definable.
We say that $X$ is \textbf{$\Psi\Delta_{\phi}$-internal} if $X \subseteq \image(f)$ for some definable $f \colon (\Psi\cup\Delta_{\phi})^m \to \Gamma^n_{\infty}$.
\end{def*}

\noindent
Note that the collection of $\Psi\Delta_{\phi}$-internal definable subsets of $\Gamma^n_{\infty}$ (for a fixed $n$) form an ideal of definable sets on $\Gamma^n_{\infty}$.
Hence:
\begin{lemma}\label{lem_dsetPsiint}
If $X\subseteq\Gamma^n$ can be covered by finitely many affine maps $\Psi^m\times\Delta_{\phi}\to\Gamma^n$, then $X$ is $\Psi\Delta_{\phi}$-internal.
\end{lemma}

\noindent
In case $\phi=\infty$, \cite[Theorem~6.3]{GehretACTlog} shows that each definable $F \colon \Psi \to \Gamma_{\infty}$ is given piecewise by so-called $s$-functions.
Now we generalize this notion in order to characterize definable functions $F \colon \Psi^n \to \Gamma_{\infty}$ in Proposition~\ref{characterization_of_functions_from_Psin}, answering a question Hieronymi asked the first-named author.
This proposition is not used except in Corollary~\ref{cor_deffunctiondset}. In particular, our proof of the Small Sets Theorem does not require it. Nevertheless, it gives a more precise description of $\Psi$-internal sets.

\begin{def*}
We call $F \colon \Psi^m \to \Gamma_{\infty}$ a \textbf{generalized $s$-function} if for $\alpha=(\alpha_0,\dots,\alpha_{m-1})\in \Psi^m$
\[
F(\alpha)\ = \ \sum_{i=0}^{m-1} \sum_{j=0}^{n-1} q_{i,j}s^{k_j}(\alpha_i) + \beta,
\]
where $\beta \in \Gamma_{\infty}$, $q_{i,j}\in \Q$, and $k_j\in \Z$ for $i=0,\dots,m-1$ and $j=0,\dots,n-1$.
Here, $s^{k}\coloneqq p^{-k}$ for $k \in \Z^{<}$.
\end{def*}

\begin{prop}\label{characterization_of_functions_from_Psin}
Let $F \colon \Psi^n \to \Gamma_{\infty}$ be definable.
Then $F$ is given piecewise by (finitely many) generalized $s$-functions.
\end{prop}
\begin{proof}
Let $(\Gamma^*, \psi^*)$ be a $|\Gamma|^+$-saturated elementary extension of $(\Gamma,\psi)$ and $(\alpha_0^*, \dots, \alpha_{n-1}^*) \in (\Psi^*)^n$.
By universal axiomatization and quantifier elimination, we have $F(\alpha_0^*,\dots,\alpha_{n-1}^*) \in \Gamma\langle \alpha_0^*,\dots,\alpha_{n-1}^*\rangle$.
Let $i_0=\min\{i<n : \alpha_i^*\notin\Psi\}$.
Then either $\alpha_{i_0}^*>\Psi$ or there is a nonempty downward closed $B\subsetneq\Psi$ with $s(B)\subseteq B$ and $B<\alpha_{i_0}^*<\Psi^{>B}$.
In the former case appealing to \cite[Lemma~4.11]{GehretACTlog} and its proof, and in the latter to \cite[Lemma~4.12]{GehretACTlog} and its proof, we have $\Gamma\langle \alpha_{i_0}^*\rangle = \Gamma\oplus\bigoplus_{k\in\Z} \Q s^k(\alpha_{i_0}^*)$ as $\Q$-linear subspaces of $\Gamma^*$.
Note that in applying \cite[Lemma~4.12]{GehretACTlog}, condition 2.\ is satisfied because $B<\alpha_{i_0}^*<\Psi^{>B}$, $s(B)\subseteq B$, and $s$ is increasing on the downward closure of $\Psi$ in $\Gamma$ by \cite[Corollary~3.6]{GehretACTlog}.
Now take $i_1 = \min\{i<n: \alpha_{i}^* \notin \Psi \cup \{s^k(\alpha_{i_0}^*) : k \in \Z\} \}$ and proceed by induction to obtain $\{i_0,\dots,i_{m-1}\}\subseteq n$ with
\[
\Gamma\langle \alpha_0^*, \dots, \alpha_{n-1}^*\rangle\ =\ \Gamma\oplus\bigoplus_{j=0}^{m-1}\bigoplus_{k\in\Z} \Q s^k(\alpha_{i_j}^*).
\]
Then compactness yields a covering of $F \colon \Psi^n \to \Gamma_{\infty}$ by finitely many generalized $s$-functions.
It follows that $F$ is given piecewise by these generalized $s$-functions. \end{proof}

\noindent
One could give a proof without using compactness by doing induction on $n$ and terms and using \cite[Theorem~6.6 and Corollary~6.7]{GehretACTlog} in a fixed model of~$T_{\log}$.

\begin{lemma}\label{lem_gensfunprimdset}
Let $F \colon \Psi^m \to \Gamma_\infty$ be a generalized $s$-function. Then $\image(F)\subseteq\image(G)\cup\{\infty\}$ for some $\Psi$-function $G\colon \Psi^I\to \Gamma$.
\end{lemma} 
\begin{cor}\label{cor_deffunctiondset}
The image of every definable $F \colon \Psi^n \to \Gamma$ is $\infty$-small.
\end{cor}

\noindent
Note that this corollary gives a direct proof of \eqref{smallequiv9}$\Rightarrow$\eqref{smallequiv5} for $\phi = \infty$ in the Small Sets Theorem.
Next we connect internality to the related model-theoretic notion of being coanalysable relative to a definable set, from \cite{herwighrushovskimacpherson}.
We follow the presentation of \cite{adh-dimension}.
As for internality, we focus on coanalysability relative to~$\Psi\Delta_{\phi}$.

\begin{def*}
Let $X\subseteq\Gamma^n$ be defined by an $\calL_{\log}\cup\{\phi,\gamma\}$-formula, where $\gamma\in\Gamma^m$.
We say that $X$ is \textbf{$\Psi\Delta_{\phi}$-coanalysable} if for all extensions
$(\Gamma_1,\psi_1)\preccurlyeq(\Gamma_2,\psi_2)\models \Th_{\calL_{\log}\cup\{\phi,\gamma\}}(\Gamma,\psi)$, if the $\Psi$-set and the subgroup $\Delta_{\phi}$ do not grow, then neither does the interpretation of~$X$.
\end{def*}

\noindent
This definition is the relevant specific instance of a general model-theoretic definition, and is justified by \cite[Proposition~2.4]{herwighrushovskimacpherson}, which provides equivalent formulations in that general context; see also \cite[Section~6]{adh-dimension} for an exposition.
What we need here is the following observation.
\begin{remark}\label{rem_coanal}
If $X \subseteq \Gamma^{n}_{\infty}$ is definable and
    $\Psi\Delta_{\phi}$-internal, then $X$ is $\Psi\Delta_{\phi}$-coanalysable.
\end{remark}

\noindent
In fact, using the appropriate general definitions, internality always implies coanalysability, but the converse fails in some settings, including in the differential field $\mathbb{T}$ of logarithmic-exponential transseries, where internality and coanalysability are taken relative to the constant field \cite{adh-dimension}.
Conversely, we show that $\Psi\Delta_{\phi}$-internality and $\Psi\Delta_{\phi}$-coanalzyability are equivalent in models of $T_{\log}$.
For this, we connect $\Psi\Delta_{\phi}$-coanalysability to $\dim_{\phi}$.
\begin{lemma}\label{lem_coanaldphi0}
If $X \subseteq \Gamma^n_{\infty}$ is definable and $\Psi\Delta_{\phi}$-coanalysable, then $\dim_{\phi}X\leq 0$.
\end{lemma}
\begin{proof}
Since $\dim_{\phi}\Psi\Delta_{\phi}=0$ and $\dim_{\phi}$ is a dimension function, this is just \cite[Corollary~2.9]{angel2016bounded}.
\end{proof}

\subsection{Dp-rank}\label{dp_rank_subsection}

Let $X$ be a definable set and $\kappa$ a cardinal. Then the dp-rank of $X$ is less than $\kappa$ (written $\dpr(X)<\kappa$) if, in every elementary extension $(\Gamma^*,\psi^*)$ of $(\Gamma,\psi)$, there does not exist a collection of formulas $(\phi_{\alpha}(x,y_{\alpha}))_{\alpha<\kappa}$, an array of elements $(b_{\alpha,i})_{\alpha<\kappa,i<\omega}$ from $\Gamma^*$ and a family of elements $(a_{\eta})_{\eta\in\omega^{\kappa}}$ from $X^*$ such that
\[
\phi_{\alpha}(a_{\eta},b_{\alpha,j}) \quad \iff \quad j \ = \ \eta(i)\quad\text{for all $\alpha,j$.}
\]
The dp-rank of $X$ is equal to $\kappa$ (written $\dpr(X) = \kappa$) if $\dpr(X) < \kappa^+$ but $\dpr(X) \not< \kappa$. We say that $X$ is \textbf{dp-finite} if $\dpr(X)=n$ for some $n$. Note that $\dpr(X) <\aleph_0$ does not imply dp-finite in general, as $\dpr(X)$ may be less than $\aleph_0$ but not less than any $n$.

\begin{fact}\label{dpfacts}
Let $Y$ be another definable set and let $f$ be a definable function.
\begin{enumerate}
\item If $X \subseteq Y$, then $\dpr(X)\leq \dpr(Y)$.
\item If $X,Y$ are dp-finite, then so is $X\times Y$. 
\item $\dpr(f(X)) \leq \dpr(X)$. 
\item $\dpr(\Psi) = 1$.
\end{enumerate}
\end{fact}
\begin{proof}
(1) and (3) are straightforward and (2) holds by~\cite{KOU13}. For (4), we know by Fact~\ref{fact_Psistabembed} that $\Psi$ is stably embedded as a model of $\Th(\N,<)$, and this theory has dp-rank one~\cite{DGL11}.
\end{proof}

\begin{lemma}\label{dpint}
Suppose that $X \subseteq \Gamma$ has nonempty interior. Then $\dpr(X) = \aleph_0$.
\end{lemma}
\begin{proof}
First, since $T_{\log}$ is countable and NIP, we have $\dpr(X)<\aleph_1$; see~\cite[Remark 2.3]{KOU13}. Since $X$ has nonempty interior, it contains a coset of $\Delta_\phi$ for some sufficiently large $\phi \in \Psi$. By Fact~\ref{dpfacts}, it is enough to show that $\dpr(\Delta_\phi) = \aleph_0$. Arguing as in~\cite[Theorem 3.6 (8)]{GehretKaplan-ACTlogdistal}, one can easily find for each $\varepsilon \in \Delta_\phi^>$ an infinite definable discrete subset of $(0,\varepsilon)$, so $\dpr(X) \not<\aleph_0$ by~\cite[Theorem 2.11]{DG17}.
\end{proof}

\section{D-minimality criterion and applications}\label{sec_dminimality_criterion}

\noindent
In this section, we give a general criterion for when the expansion of a topological theory by a collection of unary functions is d-minimal. This criterion systematizes a number of earlier results, dating back to van den Dries's 1985 proof that the real field with a predicate for the integer powers of two is d-minimal~\cite{vdd-realspowerstwo}. Instead of working with the predicate directly, van den Dries instead considers the real field expanded by the function $\lambda$ that takes $x>0$ to the largest power of two less than or equal to $x$. It is shown by a straightforward induction that terms in this extended language are given locally by semialgebraic functions, off of a finite union of discrete sets. D-minimality follows from this fact and a quantifier elimination result. 

\medskip\noindent
D-minimality for $(\R,2^\Z)$ was generalized by Miller, who showed that $(\tilde{\R},\alpha^\Z)$ is d-minimal for any polynomially bounded o-minimal expansion $\tilde{\R}$ of the real field with field of exponents $\Q$ and any $\alpha>0$~\cite{miller-tameness}. The general structure of the proof is essentially the same, though some new facts about valuation theory for o-minimal structures are needed for the quantifier elimination result. This was further extended by Friedman and Miller, who showed that d-minimality for $(\tilde{\R},\alpha^\Z)$ is preserved after adding all subsets of all cartesian powers of $\alpha^\Z$~\cite{FM01}.

\medskip\noindent
This process of describing terms in an expansion locally by functions definable in a base theory is relatively straightforward over the reals, but takes a bit more care when working over arbitrary topological structures. Here, we isolate sufficient conditions that allow us
to prepare terms in this way. These ``preparations,'' made explicit here, were used in~\cite{kaplan2023dichotomy} to show that the expansion of a power bounded o-minimal field $\bm{R}$ by a monomial group $\mathfrak{M}$ (that is, the image of a section of a $T$-convex valuation) is d-minimal.

\medskip\noindent
Below is a list of structures and theories that can be easily shown to be d-minimal, using our criterion along with known quantifier elimination results. We note that our criterion does not (to our knowledge) establish d-minimality in the stronger sense defined by Fornasiero in~\cite[Definition 9.1]{fornasiero-dimmatroiddense}. In particular, we see no way to prove Corollary~\ref{corfornasierodmin} without the Kinda Small and Very Small Sets Theorems.
\begin{enumerate}
\item Our criterion generalizes the methods used to establish d-minimality for $(\R,2^\Z)$, $(\tilde{\R},\alpha^\Z)$, and $(\bm{R},\mathfrak{M})$ as mentioned above, so these examples fall under our framework.
\item The field $(\Q_p,p^\Z)$ of $p$-adics with a predicate for the powers of $p$. This can be construed as a two-sorted structure $(\Q_p,\Z; v,\pi)$, where $v\colon \Q_p^\times \to \Z$ is the $p$-adic valuation and where $\pi\colon \Z\to \Q_p$ is the cross-section $z\mapsto p^z$. Then d-minimality follows from our criterion below, along with Ax and Kochen's quantifier elimination~\cite{AK66}. D-minimality was shown by Scowcroft, also using an induction on complexity of terms~\cite{scowcroft-moredefpadic}.
\item The expansion of an o-minimal structure by an \emph{iteration sequence} (a predicate for the iterates of a sufficiently fast-growing definable function). Quantifier elimination and d-minimality were established by Miller and Tyne~\cite{MT06}. 
\item A \emph{tame pair} of o-minimal fields: an o-minimal field $\bm{R}$ expanded by a predicate for a proper elementary substructure $\bm{A} \prec \bm{R}$ that is Dedekind complete in $\bm{R}$. Van den Dries and Lewenberg showed that these expansions eliminate quantifiers when considered with the standard part map $R\to A$, which takes $x$ in the convex hull of $A$ to the closest element in $A$. D-minimality follows from quantifier elimination and our criterion. To our knowledge, d-minimality has not been previously observed, though it follows for tame pairs of real closed ordered fields by d-minimality of the differential field $\mathbb{T}$ of logarithmic-exponential transseries, which is an expansion of such a pair; see~\cite[Corollary~16.6.11]{ADAMTT}. \item D-minimality for $T_{\log}$ can quickly be established using our criterion and the quantifier elimination result in~\cite{GehretACTlog}, as opposed to using the Small Sets Theorem. We describe how in Subsection~\ref{subsec_asympdmin}.
\item In Subsection~\ref{subsec_valueddmin}, we establish d-minimality for henselian valued fields of residue characteristic zero, equipped with a section of the valuation and a lift of the residue field. This makes use of the quantifier elimination in~\cite{vdD14}.
\item In Subsection~\ref{subsec_powersdmin}, we show how one can use our criterion to show that $(\tilde{\R},\alpha^\Z)^{\#}$, the expansion of a polynomially bounded o-minimal expansion $\tilde{\R}$ of the real field with field of exponents $\Q$ by all subsets of all cartesian powers of $\alpha^\Z$ for some $\alpha>0$, is d-minimal. As mentioned above, this was originally shown by Friedman and Miller, albeit by very different methods~\cite{FM01}. We approach this by showing that the quantifier elimination result of Miller for $(\tilde{\R},\alpha^\Z)$ can be extended to this expansion by considering this structure as a two-sorted structure.
\end{enumerate}

\subsection{The d-minimality criterion}\label{d_minimality_criterion_statement_subsection}
For the remainder of this section, $\calL$ is a multi-sorted language in the sorts $\mathcal{S}$, $s$ ranges over $S$, $\bm{M}$ is an $\calL$-structure, and $T$ is a complete $\calL$-theory. 

\medskip\noindent
We also fix a family $\chi=(\chi_s(x_s;y_s))_{s\in\mathcal{S}}$ of partitioned $\mathcal{L}$-formulas $\chi_s(x_s;y_s)$ where $x_s$ is a variable in the sort $s$, and $y_s$ is a multivariable.

\begin{definition}
We say that $\bm{M}$ is a \textbf{topological $\mathcal{L}$-structure} (with respect to $\chi$) if for each $s\in\mathcal{S}$, the family:
\[
\{\chi_s(M_s;a):a\in M_{y_s}\}
\]
is a basis for a topology on $M_s$. 
We say that $T$ is a \textbf{topological $\calL$-theory} (with respect to $\chi$) if $\bm{M}$ is a topological $\calL$-structure (with respect to $\chi$) for every model $\bm{M}$ of $T$.
\end{definition}

\noindent
Note that if $\bm{M}$ is a topological $\calL$-structure, then $\operatorname{Th}_{\calL}(\bm{M})$ is a topological $\calL$-theory. 
For a tuple of sorts $\bm{s}$, we let $M_{\bm{s}}$ denote the corresponding cartesian product. When $\bm{M}$ is a topological $\calL$-structure, we construe such a cartesian product $M_{\bm{s}}$ as a topological space using the topology generated by the $\chi_s$ and the product topology. 
If $X\subseteq M_{\bm{s}}$ is definable, then $\operatorname{int}(X)$ and $\operatorname{br}(X)$ are also definable.

\medskip\noindent
\emph{For the rest of this section $\bm{M}$ is a topological $\calL$-structure and $T$ is a topological $\calL$-theory.} 

\medskip\noindent
We say that $\bm{M}$ is $T_1$ at $s$ if the topology on $M_s$ is $T_1$, and we say that $T$ is $T_1$ at $s$ if every model $\bm{M}$ of $T$ is $T_1$ at $s$. Note that if $\bm{M}$ is $T_1$ at $s$, then so is $\operatorname{Th}_{\calL}(\bm{M})$.

\begin{definition}\label{d_minimality_def_new}
Suppose $\bm{M}$ and $T$ are $T_1$ at $s$. We say that $\bm{M}$ is \textbf{d-minimal} at $s$ if every definable subset of $M_s$ either has interior or is the union of finitely many discrete subsets of $M_s$. We say that $T$ is \textbf{d-minimal} at $s$ if $\bm{M}$ is d-minimal at $s$ for every model $\bm{M}$ of $T$.
\end{definition}

\noindent
We fix a distinguished sort $s_0 \in \calS$, and we write $M_0$ in place of $M_{s_0}$. We make the following assumptions on our topological $\calL$-theory $T$:
\begin{enumerate}[label=(\Roman*)]
\item\label{assumptionI} $T$ is $T_1$ and d-minimal at $s_0$. \item\label{assumptionII} For every model $\bm{M} \models T$, every $s \in \calS$ and every $\calL$-definable function $g \colon M_0\to M_s$, the set $\operatorname{cl}(\mathsf{Discont}(g))$ is a finite union of discrete sets. 
\end{enumerate}

\begin{remark}
By assuming that $T$ (and not just $\bm{M}$) is d-minimal at $s_0$, we obtain a uniform version of \ref{assumptionII}: For every finite tuple of sorts $\bm{s}$, every $s \in \calS$, and every $\calL$-definable function $g\colon M_{\bm{s}}\times M_0\to M_s$, there is $N$ such that $\operatorname{cl}(\mathsf{Discont}(g_x))$ is a union of $N$ discrete sets for all $x \in M_{\bm{s}}$. This uses that being a union of at most $N$ discrete sets is a definable condition;
see Remark~\ref{CB_remark}.
\end{remark}

\begin{remark}
Note that \ref{assumptionII} does not follow from \ref{assumptionI} in general, even though the set of discontinuities is definable. Indeed, Johnson notes that in the real field with the Sorgenfrey (lower limit) topology, every definable set with empty interior is finite, but $x \mapsto -x$ is nowhere continuous~\cite[Remark 1.16]{Jo24}.
\end{remark}

\noindent
Consider $\calL(\fr{F})\coloneqq \calL\cup\fr{F}$, where $\fr{F}$ is a set of new unary function symbols. Let $T(\fr{F})$ be a complete $\calL(\fr{F})$-theory extending $T$. Given $\bm{M}\models T$ we denote by $(\bm{M},\fr{F})$ an expansion of $\bm{M}$ to a model of $T(\fr{F})$.

\begin{prop}[d-minimality criterion]\label{d-minimality_criterion_ms}
Suppose the following conditions hold:
\begin{enumerate}[label=(\Alph*)]
\item\label{conditionA} For every $\calL(\fr{F})$-formula $\varphi(t)$ where $t$ is a unary  variable in the sort $s_0$, there exist an $\calL$-formula $\varphi^*(x_1,\ldots,x_n)$ and $\calL(\fr{F})$-terms $\tau_1(t),\ldots,\tau_n(t)$ such that:
\[
T(\fr{F}) \ \vdash \ \varphi(t)\leftrightarrow \varphi^*(\tau_1(t),\ldots,\tau_n(t)).
\]
\item\label{conditionB} Every new function symbol $\fr{f}\colon M_{s'}\to M_{s}$ in $\fr{F}$ is locally constant off of a finite union of strongly discrete sets; see Definition~\ref{def:stronglydiscrete}.\end{enumerate}
Then $T(\fr{F})$ is d-minimal at $s_0$.
\end{prop}

\noindent
Our criterion for d-minimality ensures that $\calL(\fr{F})$-terms are given locally by $\calL$-definable functions off of a finite union of discrete sets. The following definition makes this more precise:

\begin{definition}
Suppose $\tau\colon M_0\to M_s$ is an $\calL(\fr{F})$-term, with $s \in \calS$. We define a \textbf{preparation} of $\tau$ to be a triple $(\bm{s},B,f)$ consisting of a finite tuple of sorts $\bm{s}$, an $\calL(\fr{F})$-definable set $B\subseteq M_{\bm{s}}\times M_0$, and an $\calL$-definable function $f \colon M_{\bm{s}}\times M_0\to M_s$ such that:
\begin{enumerate}
\item $B_x$ is an open subset of $M_0$ for each $x\in M_{\bm{s}}$ and $B_x\cap B_{x'}=\0$ for $x\neq x'$,
\item $M_0\setminus \bigcup_{x}B_x$ is a finite union of discrete sets,
\item $\tau(t)=f_x(t)$ for each $x$ and for all $t\in B_x$.
\end{enumerate}
Moreover, if a preparation $(\bm{s},B,f)$ of $\tau$ additionally satisfies:
\begin{enumerate}
\setcounter{enumi}{3}
\item $f_x|_{B_x}\colon B_x\to M_s$ is continuous for each $x$,
\end{enumerate}
then we say that $(\bm{s},B,f)$ is a \textbf{continuous} preparation of $\tau$.
\end{definition}

\begin{lemma}\label{lem:contprep_ms}
Let $\tau\colon M_0\to M_s$ be an $\calL(\fr{F})$-term. 
If $\tau$ has a preparation, then $\tau$ has a continuous preparation.
\end{lemma}
\begin{proof}
Let $(\bm{s},B,f)$ be a preparation of $\tau$.
For each $x$ define $D_x\coloneqq\operatorname{cl}(\mathsf{Discont}(f_x|_{B_x}))$. Then there is $N \in \N$ such that $D_x\subseteq B_x$ is a union of $N$ discrete sets for each $x$.
It follows from Lemma~\ref{discrete_gluing_lemma} that $D\coloneqq\bigcup_{x}D_x$ is a union of $N$ discrete sets. Next, define $B^*\subseteq B$ by declaring for each $x\in M_{\bm{s}}$:
\[
B_x^* \ \coloneqq \ B_x\setminus D_x \ = \ B_x\setminus D.
\]
We claim that $(\bm{s},B^*,f)$ is a continuous preparation. The definition of $B_x^*$ ensures that (3) and (4) are satisfied.
Since each $D_x$ is closed, it follows that $B_x^*$ is open, which is (1). 
Finally, for (2) note that:
\[
\textstyle
M_0\setminus\bigcup_xB_x^* \ =\  (M_0\setminus\bigcup_xB_x)\cup D,
\]
which is a finite union of discrete sets.
\end{proof}

\begin{lemma}
If $\tau$ is a variable, then $(\0,M_0,\tau)$ is a continuous preparation of $\tau$, where $\0$ is the empty tuple of sorts.
\end{lemma}

\begin{lemma}
Let $\bm{s} = (s_1,\ldots,s_n)$ be a tuple from $\calS$, let $s \in \calS$, let $\sigma\colon M_{\bm{s}}\to M_s$ be an $\calL$-term and let $\tau_i\colon M_0\to M_{s_i}$ be an $\calL(\fr{F})$-term for each $i = 1,\ldots,n$. If each $\tau_i$ has a preparation, then the $\calL(\fr{F})$-term $\tau\coloneqq\sigma(\tau_1,\ldots,\tau_n)$ has a preparation.
\end{lemma}
\begin{proof}
For each $i=1,\ldots,n$, let $(\bm{s}_i,B_i,g_i)$ be a preparation of $\tau_i$. Next, set $\bar{\bm{s}} = (\bm{s}_1,\ldots,\bm{s}_n)$ and define $B\subseteq M_{\bar{\bm{s}}}\times M_0$ and $g\colon M_{\bar{\bm{s}}}\times M_0\to M_s$ by declaring for all $x=(x_1,\ldots,x_n)\in M_{\bm{s}_1}\times\cdots\times M_{\bm{s}_n}$:
\[
B_x \ \coloneqq  \ B_{1,x_1}\cap\cdots\cap B_{n,x_n},\qquad g(x,t) \ \coloneqq  \ \sigma(g_1(x_1,t),\ldots,g_n(x_n,t)).
\]
We claim that $(\bar{\bm{s}},B,g)$ is a preparation of $\tau$. Conditions (1) and (3) are clear. For (2), note that:
\[
\textstyle M_0\setminus (\bigcup_xB_x) \ =\  \bigcup_{1\leq i\leq n}(M_0\setminus \bigcup_{x_i}B_{i,x_i}),
\]
which is a finite union of discrete sets.
\end{proof}

\begin{lemma}
Suppose $\tau =\fr{f}(\sigma)$ for some $\calL(\fr{F})$-term $\sigma\colon M_0\to M_{s'}$ and some $\fr{f}\colon M_{s'}\to M_{s}$ in $\fr{F}$. Suppose also that:
\begin{enumerate}
\item $\sigma$ has a preparation, and
\item $\fr{f}$  is locally constant off of a finite union of strongly discrete sets.
\end{enumerate}
Then $\tau$ has a preparation.
\end{lemma}
\begin{proof}
Let $(\bm{s},C,h)$ be a preparation of $\sigma$, which by Lemma~\ref{lem:contprep_ms} we may assume to be continuous. Define $B\subseteq M_{\bm{s}}\times M_{s}\times M_0$ by:
\[
B_{x,y} \ \coloneqq  \ \operatorname{int}\{t\in C_x:\fr{f}(h_x(t))=y\} \ \subseteq \ C_x
\]
and define $f\colon M_{\bm{s}}\times M_{s}\times M_0\to M_{s}$ by:
\[
f(x,y,t) \ \coloneqq  \ y.
\]
We claim that $((\bm{s},s),B,f)$ is a preparation of $\tau$.

(1) By definition, each $B_{x,y}$ is open. Next suppose $t\in B_{x,y}\cap B_{x',y'}$. Then $t\in C_x\cap C_{x'}$, and thus $x=x'$ since $(\bm{s},C,h)$ is a preparation. By definition of $B_{x,y}$, we have $y = \fr{f}(h_x(t)) = y'$. 

(2) As $T$ is d-minimal at $s_0$, we find $N \in \N$ such that $\operatorname{br}(h_x^{-1}(z))$ is a union of $N$ discrete sets for each $(x,z)\in M_{\bm{s}}\times M_{s'}$. Then Proposition~\ref{key_discreteness_prop2_new} gives $N' \in \N$ such that $C_x\setminus\bigcup_yB_{x,y}$ is a union of $N'$ discrete sets. Thus by Lemma~\ref{discrete_gluing_lemma} we have that $\bigcup_x(C_x\setminus \bigcup_yB_{x,y})$ is a union of $N'$ discrete sets as well. It follows that
\[
\textstyle M_0\setminus\bigcup_{x,y}B_{x,y} \ = \ (M_0\setminus\bigcup_{x}C_x) \ \cup \ \bigcup_x(C_x\setminus \bigcup_yB_{x,y})
\]
is a finite union of discrete sets.

(3) Let $x,y$ and $t\in B_{x,y}$ be arbitrary. Using that $(\bm{s},C,h)$ is a preparation for $\sigma$ and $t\in B_{x,y} \subseteq C_x$, we see that
\[
\tau(t) \ = \ \fr{f}(\sigma(t)) \ =\ \fr{f}(h(x,t))\ = \ y \ = \ f_{x,y}(t).\qedhere
\]
\end{proof}

\begin{cor}\label{all_terms_preps_ms}
If each $\fr{f}\in\fr{F}$  is locally constant off of a finite union of strongly discrete sets, then every $\calL(\fr{F})$-term has a (continuous) preparation.
\end{cor}

\noindent
We can now put together the material on preparations to prove Proposition~\ref{d-minimality_criterion_ms}. 

\begin{proof}[Proof of Proposition~\ref{d-minimality_criterion_ms}]
Let $D\subseteq M_0$ be $\calL(\fr{F})$-definable. By removing the interior of $D$ (which is open and $\calL(\fr{F})$-definable), we may assume that $D$ has empty interior. We will show that $D$ is a finite union of discrete sets. 
By \ref{conditionA}, there exist an $n$-ary $\calL$-definable relation $R(x_1,\ldots,x_n)$ and $\calL(\fr{F})$-terms $\tau_1(t),\ldots,\tau_n(t)$ such that $D$ is of the form:
\[
D \ = \ \{t\in M_0:R(\tau_1(t),\ldots,\tau_n(t))\}.
\]
Next, for $i=1,\ldots,n$, by \ref{conditionB} and Corollary~\ref{all_terms_preps_ms} we take preparations $(\bm{s}_{i},B_{i},f_{i})$ of each term $\tau_{i}$. Set $\bar{\bm{s}}\coloneqq (\bm{s}_1,\ldots,\bm{s}_n)$ and define  $B\subseteq M_{\bar{\bm{s}}}\times M_0$ by declaring for $x=(x_1,\ldots,x_n) \in M_{\bm{s}_1}\times \cdots \times M_{\bm{s}_n}$: 
\[
B_x \ \coloneqq  \ B_{1,x_1}\cap\cdots\cap B_{n,x_n}\ \subseteq M_0.
\]
Then each $B_x$ is open and 
\[
\textstyle M_0\setminus\bigcup_{x}B_x \ = \ \bigcup_{1\leq i\leq n} (M_0 \setminus \bigcup_{x_{i}\in M_{\bm{s}_i}}(B_{i})_{x_{i}})
\]
is a finite union of $\calL(\fr{F})$-definable discrete sets.
We have 
\begin{align*}
D \cap B_x\ &= \ \{ t \in B_x : R_{i}\big(f_{0}(x_{0},t),\dots,f_{n-1}(x_{n-1},t)\big)\ \} \\
&= \ B_x\cap 
\underbrace{\{ t \in M_0 : R_{i}\big(f_{0}(x_{0},t),\dots,f_{n-1}(x_{n-1},t)\big) \}}_{\eqqcolon C_x}.
\end{align*}
Suppose $z\in\operatorname{int}(C_x)$, so there is an open $U$ such that $z\in U\subseteq C_x$. Since $D\cap B_x$ does not have interior, it must be the case that $U\cap B_x=\varnothing$ (since $B_x$ is open), thus $z\not\in D\cap B_x$. In particular, $D\cap B_x=B_x\cap \operatorname{br}(C_x)$. By our assumption that $T$ is d-minimal at $s_0$, $\operatorname{br}(C_x)$ is a union of $N$ discrete sets, for some $N \in \N$ not depending on $x$, so $D \cap B_x$ is a union of $N$ discrete sets. 

By Lemma~\ref{discrete_gluing_lemma}, the set $D \cap \bigcup_{x \in M_{\bar{\bm{s}}}} B_x$ is a union of $N$ discrete sets.
Finally,
\[
\textstyle D\ =\ (D \cap \bigcup_{x \in M_{\bar{\bm{s}}}} B_x) \cup (D \cap (M_0 \setminus \bigcup_{x \in M_{\bar{\bm{s}}}} B_x))
\]
is a finite union of discrete sets.
\end{proof}

\subsection{Asymptotic couple}\label{subsec_asympdmin}

Here, we describe how our d-minimality criterion can be used to quickly show that $T_{\log}$ is d-minimal. Let $(\Gamma,\psi) \models T_{\log}$. We take as a base theory $T$ the theory of divisible ordered abelian groups with an infinite element, construed in the language 
\[
\calL\ \coloneqq\ \{0,\infty,+,-,<,\delta_1,\delta_2,\delta_3,\ldots\}.
\]
Then $T$ is o-minimal, so assumptions \ref{assumptionI} and \ref{assumptionII} hold. 
Moreover $\calL_{\log} = \calL \cup \{\psi,s,p\}$, so $T_{\log}$ is an $\calL_{\log}$-theory with quantifier elimination and a universal axiomatization. In particular, $T_{\log}$ satisfies condition \ref{conditionA} of Proposition~\ref{d-minimality_criterion_ms}. Lemma~\ref{loc_const_primitives} gives that $\psi$, $s$, and $p$ are all locally constant off of a discrete set  (indeed, a strongly discrete set; see Lemma~\ref{stronglydiscreteexamples}), so condition \ref{conditionB} holds as well.

\begin{remark}
Using Corollary~\ref{all_terms_preps_ms}, one can easily prove the unary version of Corollary~\ref{cor_locallin}. That is, one can show that any definable function $f \colon \Gamma\to \Gamma_\infty$ is locally linear off of a finite union of discrete sets.
\end{remark}

\subsection{Valued fields}\label{subsec_valueddmin}
Let $K$ be a field and let $v$ be a henselian valuation on $K$ of residue characteristic zero. We let $\bm{K}$ be the three-sorted structure $(K,\k,\Gamma; v,\res)$, where $v \colon K^\times \to \Gamma$ and the residue map $\res\colon K\to \k$ are assumed to be surjective (and where the residue map sends everything with negative valuation to 0). We let $\L$ be the natural language of this three-sorted structure, but Morleyized on the sorts $\k$ and $\Gamma$ (so $\L$ also includes an additional relation symbol for every $\varnothing$-definable subset of any cartesian power of $\k$ or $\Gamma$).

\medskip\noindent
We view $\bm{K}$ as a topological $\L$-structure, where the sort $K$ is equipped with the valuation topology and the sorts $\Gamma$ and $\k$ are equipped with the discrete topology. 
We let $T$ be the $\L$-theory of $\bm{K}$. Then $T$ is a topological $\L$-theory with quantifier elimination. The models of $T$ are precisely those  henselian valued fields with residue field elementarily equivalent to $\k$ and value group elementarily equivalent to $\Gamma$.

\begin{prop}
The theory $T$ satisfies assumptions \ref{assumptionI} and \ref{assumptionII}, where the distinguished sort is the home sort. 
\end{prop}
\begin{proof}
The valuation topology on $K$ is $T_1$, and by~\cite[Proposition 5.1]{Fl11}, every definable subset of $K$ either has interior or is finite. For \ref{assumptionII}, Theorem 5.1.1 in~\cite{CHRK22} gives that every $\L$-definable function $K\to K$ is continuous on a dense open (hence cofinite) set. To verify that \ref{assumptionII} holds for $\L$-definable functions $g\colon K\to \k$, we fix such a $g$ and consider the $\L$-formula $\phi(x,y)$ defining $g$, so $x$ is a variable of sort $K$ and $y$ is a variable of sort $\k$. Let $\tilde{\L}$ be the language containing sorts for $\k$ and $\Gamma$ (as a field and an ordered abelian group, both Morleyized as above) together with a sort $\k/(\k^\times)^n$ and a corresponding quotient map $\pi_n\colon \k\to \k/(\k^\times)^n$ for each $n>1$. By~\cite[Theorem 5.15]{ACGZ}, $\phi(x,y)$ is equivalent to a formula $\psi(\sigma_1(x),\ldots,\sigma_m(x),y)$ where $\psi$ is an $\tilde{\L}$-formula and each $\sigma_i$ is of the form $vP(x)$, $\res(P(x)/Q(x))$, or $\res^n(P(x))$, where $P,Q \in K[X]$ and where $\res^n\colon K\to \k/(\k^\times)^n$ is the map
\[
\res^n(a)\ \coloneqq\ \begin{cases}
0& \text{if}\ va \not\in n\Gamma,\\
\pi^n\circ \res(a/b^n)& \text{if}\ va = nvb\ \text{(this does not depend on choice of $b$).}
\end{cases}
\]
Since $v$, $\res$, and the maps $\res^n$ are locally constant away from 0, each $\sigma_i$ is locally constant off of a finite set. It follows that $g$ is also locally constant (in particular, continuous) off of a finite set. 
The same argument works for $\L$-definable functions $h\colon K\to \Gamma$.
\end{proof}

\noindent
We now extend $\L$ to $\L^*\coloneqq \L\cup\{s,\ell\}$ and $T$ to an $\L^*$-theory $T^*$ with  axioms stating that $s$ is a section of the valuation and $\ell$ is a lift of the residue field. Note that not every model of $T$ admits an expansion to a model of $T^*$, but if $\bm{K}\models T$ is $\aleph_1$-saturated, then $\bm{K}$ does admit such an expansion.

\begin{cor}
$T^*$ is d-minimal (on the home sort). 
\end{cor}
\begin{proof}
As a consequence of van den Dries's AKE theorem with lift and section~\cite[Section~5.3]{vdD14}, the theory $T^*$ eliminates quantifiers in the language $\L^*$; see also~\cite[Theorem 2.2]{Ke24}. Condition \ref{conditionA} of Proposition~\ref{d-minimality_criterion_ms} follows. For condition \ref{conditionB}, we note that both $s$ and $\ell$ are locally constant everywhere, since $\Gamma$ and $\k$ have the discrete topology. \end{proof}

\subsection{The structure \texorpdfstring{$(\tilde{\R},\alpha^\Z)^{\#}$}{(R, αZ)\#}}\label{subsec_powersdmin}

Let $\tilde{\R}$ be a polynomially bounded expansion of the reals with field of exponents $\Q$, and let $\alpha>1$. 

\medskip\noindent
It is easy to see how our criterion can be used to show that $(\tilde{\R},\alpha^\Z)$ is d-minimal (using the quantifier elimination result of Miller--van den Dries). Less easy is to see how d-minimality of $(\tilde{\R},\alpha^\Z)^{\#}$, the expansion of $\tilde{\R}$ by all subsets of all cartesian powers of $\alpha^\Z$, follows from our criterion, as the criterion only allows for new functions, not new predicates. Here, we hint at how the multi-sorted setup can be used as a workaround.

\medskip\noindent
Let $(\tilde{\R},\Z^{\#})$ be a 2-sorted structure with underlying sets $\R$ and $\Z$. We construe this structure in the language $\calL$ that contains:
\begin{enumerate}
\item For each $n$, a function symbol for each function $\R^n\to \R$ that is definable without parameters in $\tilde{\R}$;
\item the language of ordered abelian groups $\{0,+,<\}$ on $\Z$, along with a relation symbol for each subset of $\Z^n$.
\end{enumerate}
Then $(\tilde{\R},\Z^{\#})$ eliminates quantifiers in $\calL$ and any substructure of $\tilde{\R}$ is an elementary substructure.

\medskip\noindent
Now expand $(\tilde{\R},\Z^{\#})$ by the map $\lambda\colon \R^>\to \Z$ sending $r \in \R^>$ to the least $z \in \Z$ with $\alpha^z\leq r<\alpha^{z+1}$. We extend $\lambda$ to all of $\R$ by defining $\lambda(r)\coloneqq 0$ for $r\leq 0$. Let $\calL^* \coloneqq \calL\cup\{\lambda\}$. Note that $\lambda$ is locally constant off of the set $\{0\} \cup \alpha^\Z$, which is a union of two strongly discrete sets. Thus, to prove d-minimality, it is enough to show that $(\tilde{\R},\Z^{\#};\lambda)$ has QE. This can be shown via a saturated embedding test, following the proof sketch given in~\cite[Section~8.6]{miller-tameness}.

\section{Final comments}\label{section_final_comments}

\subsection{}The reader may have noticed that the pregeometries $\cl_{\phi}$ are reminiscent of the so-called \emph{small closure for lovely pairs}~\cite[Definition 4.5]{berenstein2010lovely} and \emph{dense pairs}~\cite[Definition 8.27]{fornasiero-dimmatroiddense}. We are unaware of any existing generalized pair framework that fits our situation exactly, although the following observations are in order:
\begin{itemize}
\item $(\Gamma,\Psi)$ is a proper reduct of $(\Gamma,\psi)$: using the embedding lemma~\cite[4.6]{GehretACTlog} one can change the $\psi$-value of a suitable archimedean class without changing the underlying $\Psi$-set.
\item $\Psi$ is a predicate that names an indiscernible sequence over $\Gamma$ (in the language of ordered abelian groups with $\infty$); thus $(\Gamma,\Psi)$ (and so definable sets of $\infty$-dimension zero, by Remark~\ref{dim0_definable_in_pair}) may be further analyzed in the tradition of~\cite{baldwin2000stability}.
\item $(\Gamma,\chi)$ is a proper reduct of $(\Gamma,\psi)$, where $\chi(\alpha)\coloneqq \int\psi(\alpha)$ for $\alpha \neq 0$ is the \emph{contraction map} induced by the logarithm; this is proved in~\cite[Section~7.3]{GehretACTlogNIP}. The first-order theory of $(\Gamma,\chi)$ is studied in~\cite{bautista2019precontraction}.
\item $(\Gamma,\chi,\Psi)$ is interdefinable with $(\Gamma,\psi)$: given $\alpha\in\Gamma^{\neq}$, we may define $\psi(\alpha)$ to be the unique $\phi\in\Psi$ such that $\phi-s\phi=\chi(\alpha)$; here we use the restriction $s\colon\Psi\to\Psi$ of the successor function, which is definable in the reduct.
\item We believe that $\chi$ (and thus $\psi$) is not definable in any monadic expansion of $\Gamma$.
\end{itemize}

\subsection{} The reader may have also noticed that for definable sets $X\subseteq\Gamma^n$ with $\dim_{\infty}(X)\leq 0$ (i.e., the ``$\infty$-very small sets'') there are additional invariants available for a finer analysis. We record a few of them here:
\begin{itemize}
\item the least $m$ such that $X$ is a union of $m$ discrete sets;
\item the least $m$ such that $X^{(m)}=\varnothing$;
\item the least $m$ such that $X$ can be covered by finitely many affine maps of the form $\Psi^m\to\Gamma^n$;
\item the least $m$ such that there exists a definable surjection $\Psi^m\to X$;
\item the dp-rank of $X$;
\item when $(\Gamma,\psi)$ is the standard model, we may associate to $X$ the (growth rate of the) function $\N\to\N:k\mapsto |\pi_{s^k0}(X)|$, which we expect is always eventually equal to a polynomial with rational coefficients. For example, with the $X\subseteq\Gamma$ described in the beginning of Subsection~\ref{subsec_sststatement}, we have $|\pi_{s^k0}(X)|={k-1\choose 2}+{k-1\choose 1}+{k-1\choose 0}=\frac{1}{2}k^2-\frac{1}{2}k+1$.
\end{itemize}
We make no general claims about these invariants, although some basic relations between these quantities can be read off from our results in Section~\ref{sec_sstproof}.
\subsection{} Finally, we address some possible extensions and limitations of this work, in particular as it compares to analogous results for (the asymptotic couple of) the differential field $\mathbb{T}$ of \emph{logarithmic-exponential transseries} by Aschenbrenner, van den Dries, and van der Hoeven.
Note that their main~\cite[Theorem 0.1]{adh-revisitingclosedac} can be viewed as a ``Small Sets Theorem'' for that setting.

\medskip\noindent
It is clear that $(\Gamma_{\log},\psi)$ is part of $\mathbb{T}_{\log}^{\operatorname{eq}}$. However, the induced structure on $\Gamma_{\log}$ is more than that of the $1$-sorted asymptotic couple $(\Gamma_{\log},\psi)$. In fact, the following binary map is also definable in $\mathbb{T}_{\log}^{\operatorname{eq}}$:
\[
\operatorname{sc} \colon \R\times\Gamma_{\log}\to\Gamma_{\log},\quad \operatorname{sc}(r,\gamma) \ = \ r\gamma
\]
and consequently we also should consider the $2$-sorted strict expansion $\mathbf{\Gamma}_{\mathbb{T}_{\log}}=((\Gamma_{\log},\psi),\R;\operatorname{sc})$; this is for the same reason as given in~\cite{adh-revisitingclosedac}, namely that the constant power ``map'' on $\mathbb{T}_{\log}$ induces a scalar multiplication by the constant/residue field $\R$ on the value group $\Gamma_{\log}$.
Of course, the $2$-sorted setting introduces new types of discrete sets that are not captured by the Small Sets Theorem, such as $\R s0\coloneqq\operatorname{sc}(\R,s0)$, but we expect them to be ``orthogonal'' to the discrete set $\Psi$ in a relevant sense.

\medskip\noindent
As in~\cite{adh-revisitingclosedac}, we also know that for nonzero differential polynomials $G(Y)\in\mathbb{T}_{\log}\{Y\}$ the subset $\{vy:y\in\mathbb{T}_{\log}^{\times}:G(y)=0\}$ of $\Gamma_{\log}$ is discrete since~\cite[Corollaries 14.3.10, 14.3.11]{ADAMTT} also apply to $\mathbb{T}_{\log}$. For example, the differential polynomial $G(Y)=x(YY''-(Y')^2)+YY'$ yields the discrete set $\R s0$ in $\Gamma_{\log}$; this is because the nonzero zeros of $G(Y)$ satisfy $(xY^{\dagger})'=0$, i.e., $Y^{\dagger}\in\R x^{-1}$. 
Conversely, certain discrete sets like $\Psi$ cannot occur of the form ``$v(Z(G)^{\neq})$'' by~\cite[Corollary 13.4.5]{ADAMTT}.

\medskip\noindent
We do not know whether $\mathbf{\Gamma}_{\mathbb{T}_{\log}}$ is purely stably embedded in $\mathbb{T}_{\log}^{\operatorname{eq}}$: a positive answer would require an elimination theory for $\mathbb{T}_{\log}$, and a negative answer would require a counterexample. This was answered in the negative for $\mathbb{T}$ in~\cite[p.\ 532]{adh-revisitingclosedac}.
However, their counterexample produces an injection $c\mapsto [v(e^{e^{cx}})]:\R^>\to[\Gamma_{\mathbb{T}}^{\neq}]$ into the archimedean classes of $\Gamma_{\mathbb{T}}$; such a phenomenon is impossible in our setting since $[\Gamma_{\log}^{\neq}]\cong\Psi$ is countable.

\appendix

\section{Topology}\label{sec_appendix}

\noindent
Our main theorems concern the agreement (or disagreement) of various ideals of definable sets. In this appendix, we recall some basic facts we need about those ideals that arise from topology. In particular, we observe that under mild topological assumptions (e.g., if $X$ is $T_1$ with no isolated points), we get a linear inclusion of ideals:
\[
\text{finite} \ \subseteq \ \text{d-finite} \ \subseteq \ \text{finite unions of discrete sets} \ \subseteq \ 
\text{nowhere dense}. \]
See Lemmas~\ref{finite_is_dfinite_T1},~\ref{CBsmall_implies_finite_union_discrete}, and~\ref{discrete_implies_nowhere_dense_under_mild_assumptions}.
Furthermore, in Subsection~\ref{border_loc_const_functs_subsection} we establish Proposition~\ref{key_discreteness_prop2_new}, which is needed for the d-minimality criterion~\ref{d-minimality_criterion_ms}.

\medskip\noindent
\emph{Throughout, $X$ and $Y$ are topological spaces and we let $A,B$ range over subsets of $X$.}

\subsection{Preliminaries}
We denote the \textbf{interior} and \textbf{closure} of $A$ in $X$ by $\operatorname{int}(A)$ and $\operatorname{cl}(A)$. We may use subscripts if we want to emphasize the ambient space. For example:

\begin{lemma}\label{cont_map_bigger_closure}
If $f\colon X\to Y$ is continuous, then for any $A\subseteq X$ we have $\operatorname{cl}_X(A)\subseteq f^{-1}(\operatorname{cl}_Y(f(A)))$, equivalently, $f(\cl_X(A))\subseteq\cl_Y(f(A))$.
\end{lemma}

\noindent
We occasionally make the following assumptions about our topological space:

\begin{definition}
We say that $X$ is $T_1$ if every singleton is closed, and we say that $X$ has \textbf{no isolated points} if every singleton is not open.
\end{definition}

\subsection{Sets with empty interior}
We say a set $A$ has \textbf{empty interior in $X$} if $\operatorname{int}_X(A)=\varnothing$. When the ambient space $X$ is understood from context, we just say $A$ has \textbf{empty interior}. Conversely, we say $A$ \textbf{has interior} if $A$ does not have empty interior.

\medskip\noindent
The sets with empty interior in general do not form an ideal; moreover, the ideal that they generate is often the improper ideal (e.g., consider $\Q$ and $\R\setminus\Q$ in $\R$ with the usual topology.)

\subsection{The ideal of nowhere dense sets}
We say $A$ is \textbf{nowhere dense in $X$} if $\operatorname{int}_X(\operatorname{cl}_X(A))=\varnothing$. When the ambient space $X$ is understood from context, we just say $A$ is \textbf{nowhere dense}.

\medskip\noindent
The nowhere dense subsets of $X$ form an ideal of subsets of $X$, which is closed under taking closures:

\begin{lemma}[Nowhere dense sets form an ideal]\label{nowheredense_ideal_properties}
\hfill\begin{enumerate}
\item\label{nowheredense_ideal_union} A finite union of nowhere dense sets is nowhere dense.
\item\label{nowheredense_downward_closed} If $A\subseteq B$ and $B$ is nowhere dense, then $A$ is nowhere dense.
\item\label{nowheredense_closure} $A$ is nowhere dense iff $\operatorname{cl}(A)$ is nowhere dense.
\end{enumerate}
\end{lemma}

\noindent
The following is obvious:
\begin{lemma}\label{nowhere_dense_implies_empty_interior}
If $A$ is nowhere dense, then $A$ has empty interior.
\end{lemma}

\subsection{The ideal generated by discrete sets}
 Recall that a set $A$ is \textbf{discrete} (in $X$) if $A$ is a discrete space when equipped with the subspace topology.

\medskip\noindent
In general, the discrete subsets of $X$ do not form an ideal (consider $\{0\}\cup\{1/n:n\geq 1\}$ in $\R$ with the usual topology), although really we are interested in the ideal generated by the discrete sets, i.e., the sets that are finite unions of discrete sets.

\medskip\noindent
We first observe that under mild topological assumptions, the ideal generated by discrete sets is sandwiched between the ideal of finite sets and the ideal of nowhere dense sets:

\begin{lemma}\label{discrete_implies_nowhere_dense_under_mild_assumptions}
Suppose $X$ is $T_1$ and has no isolated points. Then for every $D\subseteq X$ we have:
\[
\text{$D$ is finite} \quad \Rightarrow \quad \text{$D$ is discrete} \quad \Rightarrow \quad \text{$D$ is nowhere dense}.
\]
\end{lemma}
\begin{proof}
The first implication is immediate from $T_1$. Now assume $D$ is discrete and assume towards a contradiction that $\operatorname{int}(\operatorname{cl}(D))\neq\varnothing$. Then there exists a nonempty open set $U\subseteq \cl(D)$. Note that $U\cap D\neq\varnothing$ and fix $d\in U\cap D$. Since $D$ is discrete, we can choose an open set $V\subseteq U$ such that $V\cap D=\{d\}$.
Next, next note that $W\coloneqq V\setminus\{d\}$ is open (by $T_1$) and nonempty (since $X$ has no isolated points). Moreover, we have $W\subseteq \operatorname{cl}(D)$, however $W\cap D=\varnothing$, a contradiction.
\end{proof}

\noindent
In general both directions $\Leftarrow$ can fail. Indeed, the second is the distinction between the Very Small Sets and the Kinda Small Sets when $n>1$.

\medskip\noindent
Sets in the ideal generated by discrete sets can be glued together in the following obvious way:

\begin{lemma}\label{discrete_gluing_lemma}
Suppose:
\begin{enumerate}
\item $(B_i)_{i\in I}$ is a disjoint family of open sets,
\item $(D_i)_{i\in I}$ is a family with $D_i \subseteq B_i$ for each $i \in I$,
\item there is $N \in \N$ such that each $D_i$ is a union of $N$ discrete sets.
\end{enumerate}
Then $\bigsqcup_{i\in I}D_i$ is a union of $N$ discrete sets. 
\end{lemma}

\noindent
When applying this lemma, the following concept is useful:

\begin{definition}\label{def:stronglydiscrete}
We say a set $D\subseteq X$ is \textbf{strongly discrete} if there exists a family $(V_d)_{d\in D}$ of pairwise disjoint open sets such that $V_d\cap D=\{d\}$ for each $d\in D$.
\end{definition}

\noindent
In many topological spaces of interest, we get \emph{strongly discrete} for free:

\begin{lemma}\label{stronglydiscreteexamples}
Suppose $X$ is either (i) a metric space equipped with the metric topology, or (ii) a valued abelian group equipped with the valuation topology. Then for every $D\subseteq X$ we have: $D$ is discrete iff $D$ is strongly discrete.
\end{lemma}

\noindent
Note the above lemma can fail, even for hausdorff spaces. Indeed, in the  \emph{Niemytzki's tangent disk topology}~\cite[II.82]{steen1978counterexamples} on $X=\R\times[0,+\infty)$ the axis $\R\times \{0\}$ is discrete but not strongly discrete.

\subsection{The ideal of d-finite sets}

Given $x\in X$ we say that $x$ is a \textbf{limit point} of $A$ (in $X$) if for every neighbourhood $U$ of $x$:
\[
(U\setminus\{x\})\cap A \ \neq \ \varnothing.
\]
Define the \textbf{derived set} of $A$ (in $X$) to be:
\[
A' \ \coloneqq \ \{x\in X:\text{$x$ is a limit point of $A$}\}.
\]
Note that the derived set $A'$ of $A$ depends on the ambient space $X$.
For $n\geq 0$ we define $A^{(n)}$ recursively by setting $A^{(0)}\coloneqq A$ and $A^{(n+1)}\coloneqq (A^{(n)})'$.

\medskip\noindent
Here are some basic facts about derived sets:

\begin{lemma}\label{CB_der_properties}
Given $A,B$ and a finite family $(A_i)_{i\in I}$ of subsets of $X$ we have for every $n\geq 0$:
\begin{enumerate}
\item\label{CB_der_commute_union} $(\bigcup_{i\in I}A_i)^{(n)}=\bigcup_{i\in I}A_i^{(n)}$.
\item\label{CB_der_closure_nth} $A^{(n)}\subseteq A\cup A'=\operatorname{cl}(A)$; in particular, $\operatorname{cl}(A)^{(n)}=A^{(n)}\cup A^{(n+1)}=\cl(A^{(n)})$.
\item\label{CB_der_monotone} If $A\subseteq B$, then $A^{(n)}\subseteq B^{(n)}$.
\item $A$ is discrete and closed if and only if $A'=\varnothing$.
\end{enumerate}
\end{lemma}

\noindent
In particular, if $A$ is closed, then  $(A^{(n)})_{n\geq 0}$ is a decreasing sequence of closed sets. Moreover, if $X$ is $T_1$, then $A'=\cl(A)'$, hence $(A^{(n)})_{n\geq 1}$ is a decreasing sequence of closed sets. However, this behavior need not be typical in general:

\begin{example}Suppose $X=\{a,b\}$ with the indiscrete topology. Then for $A=\{a\}$ we have:
\[
A^{(n)} \ \coloneqq \ \begin{cases}
\{a\} & \text{if $n$ even}, \\
\{b\} & \text{if $n$ odd}.
\end{cases}
\]
Thus in general $A^{(n)}$ need not be closed, the sequence $(A^{(n)})_{n\geq 0}$ need not be eventually decreasing, and $\operatorname{cl}(A)^{(n)}$ need not equal $A^{(n)}$.
\end{example}
\noindent
We introduce the following terminology:

\begin{definition}
We say that $A$ is \textbf{d-finite in $X$} if $A^{(n)}=\varnothing$ for some $n\geq 0$. When the ambient space $X$ is understood from context, we just say $A$ is \textbf{d-finite}.
\end{definition}

\noindent
The d-finite subsets of $X$ form an ideal of subsets of $X$, which is closed under taking closures:

\begin{lemma}[d-finite sets form an ideal]\label{CB_small_ideal_properties}
\hfill\begin{enumerate}
\item\label{CB_finite_ideal_union} A finite union of d-finite sets is d-finite.
\item\label{CB_finite_downward_closed} If $A\subseteq B$ and $B$ is d-finite, then $A$ is d-finite.
\item\label{CB_finite_closure} $A$ is d-finite iff $\operatorname{cl}(A)$ is d-finite.
\end{enumerate}
\end{lemma}
\begin{proof}
\eqref{CB_finite_ideal_union} follows from~\ref{CB_der_properties}\eqref{CB_der_commute_union}, \eqref{CB_finite_downward_closed} follows from~\ref{CB_der_properties}\eqref{CB_der_monotone}, and \eqref{CB_finite_closure} follows from~\ref{CB_der_properties}\eqref{CB_der_closure_nth}.
\end{proof}

\noindent
The ideal of d-finite sets contains the ideal of finite sets under the $T_1$ assumption:
\begin{lemma}\label{finite_is_dfinite_T1}
If $X$ is $T_1$ and $A$ is finite, then $A'=\varnothing$, so $A$ is d-finite.
\end{lemma}

\noindent
The ideal of d-finite sets is always contained in the ideal generated by discrete sets:

\begin{lemma}\label{CBsmall_implies_finite_union_discrete}
For every $A$, the difference $A\setminus A'$ is discrete; hence, if $A$ is d-finite, then $A$ is a finite union of discrete sets.
\end{lemma}

\noindent
In general the converse can fail in a $T_1$ space:

\begin{example}\label{counterexample_discrete_not_CB_finite}
Consider the ordinal $X\coloneqq\omega_1$ equipped with the order topology, and the subset $D\coloneqq\{\alpha+1:\alpha\in \omega_1\}\subseteq X$ of successor ordinals. Then $D$ is discrete, although $D$ is not d-finite.
\end{example}

\noindent
Here is a $T_1$ example without any isolated points:

\begin{example}\label{second_counterexample_discrete_not_CB_finite}
Consider the product $X=\R\times\R$. Equip each point $(r,0)$ with the usual neighbourhood basis inherited from the standard topology on $\R^2$. Equip each point $(r,x)$ with $x\neq0$ with the neighbourhood basis of vertical intervals of the form $\{(r,t):x-\varepsilon<t<x+\varepsilon\}$, where $0<\varepsilon<|x|$. Consider the set $D\coloneqq\R\times\{1/n:n>0\}$. Then $D$ is discrete, although $D^{(n)}=\R\times\{0\}$ for every $n\geq 1$. Thus $D$ is not d-finite.
\end{example}

\begin{remark}[Relation to Cantor--Bendixson derivative]\label{CB_remark}
The \emph{Cantor--Bendixson (CB) derivative} is initially defined for topological spaces $X$ as $d_{\operatorname{CB}}(X)\coloneqq X'$ as in~\cite[(6.10)]{kechris2012classical}, although one often extends this definition to arbitrary subsets $A\subseteq X$ via $d_{\operatorname{CB}}(A)\coloneqq A\setminus\operatorname{isol}(A)$ (as done in~\cite[3.28]{fornasiero2021dminimalstructures} for hausdorff spaces), where $\operatorname{isol}(A)$ is the set of isolated points of $A$ in the subspace topology induced by $X$. 
In a $T_1$-space, we have $d_{\operatorname{CB}}^n(A) = \varnothing$ if and only if $A$ is a union of at most $n$ discrete sets, which shows that this property is definable in families; see~\cite[3.29]{fornasiero2021dminimalstructures}.

Following the terminology from~\cite{kechris2012classical}, we say that $A$ has \emph{finite CB-rank} if the sequence $(d_{\operatorname{CB}}^n(A))_{n\geq 0}$ is eventually constant, in which case we call the eventual value the \emph{perfect kernel} of $A$. 
It then follows that the property ``$A$ is d-finite'' is equivalent to ``$\cl(A)$ has finite CB-rank with empty perfect kernel''. 
Note that if $A$ is d-finite, then $A$ has finite CB-rank (with empty perfect kernel), although the converse need not hold as shown in Examples~\ref{counterexample_discrete_not_CB_finite} and~\ref{second_counterexample_discrete_not_CB_finite}. We could not find an existing name for the property \emph{d-finite} in the literature.
\end{remark}

\noindent
The following implies that the product of d-finite sets is again d-finite in the product topology:

\begin{lemma}\label{CBproduct1}
Suppose $A \subseteq X$ and $C\subseteq Y$ are closed. For $k \in \N$, we have 
\[
\textstyle
(A\times C)^{(k)}\ =\ \bigcup_{m+n=k}(A^{(m)}\times C^{(n)}),
\]
where $A\times C$ is considered as a subset of the product space $X\times Y$. In particular, if $A^{(m)} = C^{(n)}= \varnothing$ for some $m,n \in \N$ with $m+n>0$, then $(A\times C)^{(m+n-1)}= \varnothing$.
\end{lemma}
\begin{proof}
The $k=1$ case $(A\times C)'=(A'\times C)\cup (A\times C')$ is routine. The general case follows by induction on $k$.
\end{proof}

\subsection{Border and locally constant functions}\label{border_loc_const_functs_subsection}

\noindent
Given a set $A\subseteq X$, we let $\operatorname{br}A \coloneqq A \setminus \operatorname{int}A$ denote the \textbf{border} of $A$.
The \emph{border} of $A$ should not be confused with the \emph{boundary} of $A$ ($=\operatorname{cl}(A)\setminus\operatorname{int}(A)$) or the \emph{frontier} of $A$ ($=\operatorname{cl}(A)\setminus A$); we do not use boundary or frontier in this paper.

\medskip\noindent
We mainly use border to detect where a function is (not) locally constant: Suppose $\fr{f}\colon Y\to Z$ is an arbitrary function. Then we have:
\[
\textstyle\{y\in Y:\text{$\fr{f}$ is locally constant at $y$}\} \ = \ \bigsqcup_{z\in Z}\operatorname{int}\fr{f}^{-1}(z),
\]
which is always an open set; we also have:
\[
\textstyle\{y\in Y:\text{$\fr{f}$ is not locally constant at $y$}\} \ = \ \bigsqcup_{z\in Z}\operatorname{br}\fr{f}^{-1}(z),
\]
which is always a closed set, being the complement of the first set.

\medskip\noindent
We conclude with Proposition~\ref{key_discreteness_prop2_new} below, which is a technical membership criterion for this ideal that we use in the proof of Proposition~\ref{d-minimality_criterion_ms}. First a lemma:

\begin{lemma}\label{int_fr_discrete_fibers_new}
Suppose $h\colon X\to Y$ is continuous and $D\subseteq Y$ is arbitrary. Then:
\[
\textstyle\operatorname{int}h^{-1}D \ \supseteq \ \bigsqcup_{d\in D}\operatorname{int}h^{-1}(d)\quad\text{and thus}\quad \operatorname{br}h^{-1}D \ \subseteq \ \bigsqcup_{d\in D}\operatorname{br}h^{-1}(d).
\]
\end{lemma}
\begin{proof}
For $d\in D$ we have $h^{-1}(d)\subseteq h^{-1}(D)$, thus $\operatorname{int}h^{-1}(d)\subseteq \operatorname{int}h^{-1}(D)$, which yields the first inclusion. The second inclusion follows from taking complements inside $h^{-1}(D)$.
\end{proof}

\begin{prop}\label{key_discreteness_prop2_new}
Let $h\colon X\to Y$ be a continuous function, and $\fr{f}\colon Y\to Z$ be a function. Suppose:
\begin{enumerate}
\item $\fr{f}\colon Y\to Z$ is locally constant outside a union of $M$ strongly discrete sets, and 
\item there is $N \in \N$ such that $\operatorname{br}h^{-1}(y)$ is a union of $N$ discrete sets for every $y\in Y$.
\end{enumerate}
Then
\[
\textstyle X\setminus\bigsqcup_z\operatorname{int}X_z \ = \ \bigsqcup_z\operatorname{br}X_z 
\]
is the union of $MN$ discrete sets, where $X_z\coloneqq h^{-1}\fr{f}^{-1}(z)$, i.e., the composition $\fr{f}\circ h\colon X\to Z$ is locally constant outside a union of $MN$ discrete sets.
\end{prop}
\begin{proof}
Set $D\coloneqq\bigsqcup_{z\in Z}\operatorname{br}(\fr{f}^{-1}(z))$. By assumption, we may take a partition $(D^i)_{1\leq i\leq M}$ of $D$ into $M$ strongly discrete sets $D^i$.

For $z\in Z$ define $U_z\coloneqq\operatorname{int}\fr{f}^{-1}(z)$ and $D_z\coloneqq\operatorname{br}\fr{f}^{-1}(z)$. Note that our assumption on $\fr{f}$ says that $D=\bigsqcup_{z}D_z$, and hence each $D_z$ is a union of the $M$ discrete sets $D_z^i\coloneqq D_z\cap D^i$ where $i=1,\ldots,M$; hence $D^i=\bigsqcup_{z\in Z}D^i_z$. Thus we have a partition:
\[
\textstyle \fr{f}^{-1}(z) \ = \ U_z\sqcup D_z \ = \ U_z\sqcup\bigsqcup_{i=1}^MD_z^i.
\]
Pulling this back along the continuous function $h\colon X\to Y$ yields:
\[
\textstyle X_z \ = \ h^{-1}\fr{f}^{-1}(z) \ = \ h^{-1}(U_z)\sqcup\bigsqcup_{i=1}^Mh^{-1}(D_z^i) \ = \ \underbrace{\textstyle h^{-1}(U_z)\sqcup\bigsqcup_{i=1}^M\operatorname{int} h^{-1}(D_z^i)}_{\text{open}}\sqcup\bigsqcup_{i=1}^M\operatorname{br}h^{-1}(D_z^i).
\]
Hence we have $\operatorname{br}X_z\subseteq \bigsqcup_{i=1}^M\operatorname{br}h^{-1}(D^i_z)\subseteq D_z$ for each $z\in Z$. Thus by Lemma~\ref{int_fr_discrete_fibers_new} we have:
\[
\textstyle \bigsqcup_z\operatorname{br}X_z \ \subseteq \ \bigsqcup_{z\in Z}\bigsqcup_{i=1}^M\operatorname{br}h^{-1}(D^i_z) \ = \ \bigsqcup_{z\in Z}\bigsqcup_{i=1}^M\bigsqcup_{d\in D^i_z}\operatorname{br}h^{-1}(d) \ = \ \bigsqcup_{i=1}^M\bigsqcup_{d\in D^i}\operatorname{br}h^{-1}(d).
\]
It suffices to show for each $i$ that $\bigsqcup_{d\in D^i}\operatorname{br}h^{-1}(d)$ is a finite union of $N$ discrete sets. Take for each $d\in D^i$ an open $V_d\subseteq Y$ such that $V_d\cap D^i=\{d\}$; moreover, since $D^i$ is strongly discrete, we may assume that the family $(V_d)_{d\in D^i}$ is pairwise disjoint. Then $(h^{-1}V_d)_{d\in D^i}$ is a disjoint family of open sets with $\operatorname{br}h^{-1}(d)\subseteq h^{-1}V_d$ for each $d\in D^i$, so $\bigsqcup_{d\in D^i}\operatorname{br}h^{-1}(d)$ is a union of $N$ discrete sets by Lemma~\ref{discrete_gluing_lemma}.
\end{proof}

\section*{Acknowledgements}
\noindent
We are grateful to the following individuals for conversations around the topics of this paper: Matthias Aschenbrenner, Johannes Aspman, Gilles Bareilles, David Bradley-Williams, Monroe Eskew,  Philipp Hieronymi, Martin Hils, Nayoon Kim, Rufus Lawrence, Jana Lep\v{s}ov\'{a}, Shenyuan Ma, Jakub Mare\v{c}ek, Julia Millhouse, Benjamin Riff, Silvain Rideau-Kikuchi, and Ale\v{s} Wodecki.

\medskip\noindent
Gehret was partially supported from the European Union’s Horizon Europe research and innovation programme under grant agreement No.~101070568. 
This research was supported by the National Science Foundation under Award No.\ DMS-2103240.
This research was funded in whole or in part by the Austrian Science Fund (FWF) 10.55776/ESP450. Part of this research was conducted while E.\ Kaplan was hosted by the Max Planck Institute for Mathematics, and he thanks the MPIM for its support and hospitality. 
For open access purposes, the authors have applied a CC BY public copyright licence to any author accepted manuscript version arising from this submission.


\providecommand{\bysame}{\leavevmode\hbox to3em{\hrulefill}\thinspace}
\providecommand{\href}[2]{#2}

\end{document}